\documentclass[11pt]{article} 
 
\usepackage{amsmath}
\usepackage{amssymb}
\usepackage{amsthm}
\usepackage{geometry}
\usepackage[colorlinks=true, linkcolor=blue, citecolor=blue]{hyperref}
\usepackage{url}
\usepackage{enumitem}
\usepackage[utf8]{inputenc}
\usepackage[english]{babel}
\usepackage{longtable}
\usepackage[dvipsnames,table,longtable]{xcolor}
\usepackage{booktabs}
\usepackage[ruled,linesnumbered]{algorithm2e}
\usepackage{bbm}
\usepackage{mathtools}
\usepackage{graphicx}
\usepackage{fix-cm}
\usepackage{subcaption}
\usepackage{pgfplots}
\usepackage[title]{appendix}
\usepackage{multirow}
\usepackage{pdflscape} 
\usepackage{scalerel}
\usepackage{lipsum}
\usepackage{comment}
\usepackage{afterpage}
\usepackage{float}

\usepackage{tikz}
\usetikzlibrary{arrows}
\usetikzlibrary{plotmarks}
\usetikzlibrary{shadings,automata,shadows,calc,arrows}
\usetikzlibrary{decorations.pathmorphing,patterns}
\usetikzlibrary{decorations.pathreplacing}
\usetikzlibrary{decorations.markings}
\usetikzlibrary{calc} 
\usetikzlibrary{shapes,snakes}
\usetikzlibrary{positioning}
\usetikzlibrary{shapes.geometric}

\usepackage{authblk}

\usepackage{lipsum}
\makeatletter
\renewcommand\footnoterule{%
  \kern-3\p@
  \hrule\@width \textwidth
  \kern2.6\p@}
\makeatother

\makeatletter
\renewcommand*{\@fnsymbol}[1]{\ensuremath{\ifcase#1\or *\or \dagger\or \ddagger\or **\or
   \mathsection\or \mathparagraph\or \|\or  \dagger\dagger
   \or \ddagger\ddagger \else\@ctrerr\fi}}
\makeatother

\allowdisplaybreaks

\newcommand{\textBF}[1]{%
    \pdfliteral direct {2 Tr 0.5 w} 
     #1%
    \pdfliteral direct {0 Tr 0 w}%
}

\newcommand{\R}{\mathbb{R}}

\newcommand{\Z}{\mathbb{Z}}
\newcommand{\N}{\mathbb{N}}

\newcommand{\Zplus}{\Z_{+}}

\newcommand{\bsubeq}{\begin{subequations}}
\newcommand{\esubeq}{\end{subequations}}
\newcommand{\BI}{\begin{itemize}}
\newcommand{\EI}{\end{itemize}}
\newcommand{\I}{\item}
\newcommand{\BE}{\begin{enumerate}}
\newcommand{\EE}{\end{enumerate}}

\newcommand{\crev}{\color{black}}
\newcommand{\cRev}{\color{black}}

\newcommand{\rankfunc}{rank}
\newcommand{\doublefunc}{double}
\newcommand{\degr}{d}
\newcommand{\density}{D}
\newcommand{\neighb}{\mathcal{N}}
\newcommand{\indicate}{\mathbbm{I}}

\newcommand{\clique}{\mathcal{K}}

\newcommand{\instance}{(G, K)}
\newcommand{\graph}{G}
\newcommand{\digraph}{\overrightarrow{G}}
\newcommand{\cycles}{\mathcal{C}}

\newcommand{\EFfeas}{\mathcal{EF}}

\newcommand{\adjmatrix}{A}
\newcommand{\K}{K}
\newcommand{\vertexset}{\mathcal{V}}
\newcommand{\edgeset}{\mathcal{E}}
\newcommand{\arcset}{\mathcal{A}}
\newcommand{\LOrankSet}{\mathcal{R}^{LO}}
\newcommand{\DVOPrankSet}{\mathcal{R}^{DVOP}}

\newcommand{\IPVR}{\mathbb{(IP^{VR})}}

\newcommand{\vertexCP}{\mathbb{(CP^{VERTEX})}}
\newcommand{\rankCP}{\mathbb{(CP^{RANK})}}
\newcommand{\combCP}{\mathbb{(CP^{COMBINED})}}
\newcommand{\naiveMP}{(\mathbb{MP}\mathbbm{1})}
\newcommand{\naiveSP}{(\mathbb{SP}\mathbbm{1})}
\newcommand{\efMP}{(\mathbb{MP}\mathbbm{2})}
\newcommand{\efSP}{(\mathbb{SP}\mathbbm{2})}
\newcommand{\efExtensive}{\mathbb{(EF)}}
\newcommand{\OGcycles}{\mathbb{(CYCLES)}}
\newcommand{\OGranks}{\mathbb{(RANKS)}}
\newcommand{\CCG}{\mathbb{(CCG)}}
\newcommand{\IP}{\mathbb{(IP)}}
\newcommand{\precvar}{p}
\newcommand{\vrvar}{x}
\newcommand{\rankvar}{r}
\newcommand{\vertvar}{v}
\newcommand{\doublevar}{y}

\newcommand{\ipindicator}{z}
\newcommand{\cliquevar}{\kappa}
\newcommand{\witnessvar}{w}

\newcommand{\rankind}{r}
\newcommand{\vertind}{v}

\newcommand{\MinD}{\texttt{MIN DOUBLE}}
\newcommand{\MinN}{\texttt{MIN NODES}}

\newcommand{\iis}{\mathcal{IIS}}

\newcommand{\numdoub}{nodes}

\newcommand{\Lorder}{(}
\newcommand{\Rorder}{)}

\newtheorem{ex}{Example}
\newtheorem{definition}{Definition}

\newtheorem{lemma}{Lemma}
\newtheorem{corollary}{Corollary}
\newtheorem{conj}{Conjecture}
\newtheorem{prop}{Proposition}
\newtheorem{thm}{Theorem}

\newtheorem{fix}{Variable Fixing Rule}

\title{Integer Programming, Constraint Programming, and Hybrid Decomposition Approaches to Discretizable Distance Geometry Problems}

\author[1]{Moira MacNeil\thanks{m.macneil@mail.utoronto.ca}}
\author[1]{Merve Bodur\thanks{bodur@mie.utoronto.ca}}

\affil[1]{Department of Mechanical and Industrial Engineering, University of Toronto}

\begin{document}
	
	\maketitle
	
	\begin{abstract}
Given an integer dimension $\K$ and a simple, undirected graph $\graph$ with positive edge weights, the Distance Geometry Problem (DGP) aims to find a realization function mapping each vertex to a coordinate in $\R^\K$ such that the distance between pairs of vertex coordinates is equal to the corresponding edge weights in $\graph$. The so-called discretization assumptions reduce the search space of the realization to a finite discrete one, which can be explored via the branch-and-prune (BP) algorithm. Given a \emph{discretization vertex order} in $\graph$, the BP algorithm constructs a binary tree where the nodes at a layer provide all possible coordinates of the vertex corresponding to that layer. The focus of this paper is finding \emph{optimal BP trees} for a class of Discretizable DGPs. More specifically, we aim to find a discretization vertex order in $\graph$ that yields a BP tree with the least number of branches. We propose an integer programming formulation and three constraint programming formulations that all significantly outperform the state-of-the-art cutting plane algorithm for this problem. Moreover, motivated by the difficulty in solving instances with a large and low density input graph, we develop two hybrid decomposition algorithms, strengthened by a set of valid inequalities, which further improve the solvability of the problem.
\end{abstract}

\section{Introduction} \label{sec:intro}

Distance Geometry is the study of problems where we wish to determine positions in a geometric space of points while preserving some known distances between the points \cite{dgpbook,mucherino2012DG}. It has wide application areas, including astronomy, where we position stars relative to each other, {\cRev and }robotics, where the distances are arm lengths and we are trying to determine a set of positions within reach of a robot  \cite{dgpbook, liberti2014euclidean, mucherino2012DG}. In molecular geometry, Nuclear Magnetic Resonance spectroscopy is used to {\crev find interatomic distances of large molecules, which are often} proteins. {\crev This process gives measurements} in two dimensions, but finding {\crev the} three dimensional structure {\crev of such molecules} is key for determining their functional properties. In this case, we are positioning the atoms in three-dimensional Euclidean space  \cite{dgpbook}. In wireless sensor localization, the network has components with fixed positions, such as routers, and we wish to determine the positions of mobile wireless sensors, such as smartphones \cite{liberti2014euclidean}. A new variant of the Distance Geometry Problem (DGP), namely dynamical DGP, has stemmed from applications such as air traffic control, crowd simulation, multi-robot formation, and human motion retargeting, all of which involve a temporal aspect \cite{mucherinoHDR,mucherino2018application}. Other applications include  statics and graph rigidity, graph drawing, and clock synchronization \cite{dgpbook, liberti2014euclidean, mucherino2012DG}. 

The DGP can be represented on a graph where the vertices are the points we would like to position and weighted edges represent known distances between pairs of points, {\crev note that this graph need not be complete.} In addition to the graph, DGP takes as input an integer $\K$, the dimension of $\R$ into which the graph is positioned. Formally, we give the definition of \cite{dgpbook}.
\begin{definition}[Distance Geometry Problem] \label{def:DGP}
	Given, an integer $\K > 0$ and a simple, undirected graph $\graph = (\vertexset,\edgeset)$ with edge weights $w  : \edgeset \to (0, \infty)$,
	find a function $x : \vertexset \to  \R^\K$ such that for all $\{u,v\} \in \edgeset$: $
	\|x(u) - x(v)\| = w(u,v)$.
\end{definition}

\begin{figure}[h]
	\begin{subfigure}[b]{0.4\textwidth}
		\begin{center}
				\begin{tikzpicture}[scale=0.8,main_node/.style={scale=0.9,circle,fill=white!80,draw,inner sep=0pt, minimum size=16pt},	line width=1.2pt]
			\node[main_node] (v0) at (3,3) {$v_0$};
			\node[main_node] (v1) at (3,0) {$v_1$};
			\node[main_node] (v2) at (0,0) {$v_2$};
			\node[main_node] (v3) at (0,3) {$v_3$};
			\draw[-] (v0) -- node[right, blue] {\footnotesize $2\sqrt{2}$}  (v1);
			\draw[-] (v1) -- node[below, blue] {\footnotesize $\sqrt{13}$}  (v2);
			\draw[-] (v0) -- node[xshift=.6cm, yshift=0.9cm,blue] {\footnotesize 5}  (v2);
			\draw[-] (v0) -- node[above, blue] {\footnotesize 2}  (v3);
			\draw[-] (v1) -- node[xshift=.6cm, yshift=-0.9cm, blue] {\footnotesize 3}  (v3);
			\draw[-] (v2) -- node[left, blue] {\footnotesize 4}  (v3);
			\end{tikzpicture}
		\end{center}
		\caption{ A complete graph with four vertices.} \label{fig:complete}
	\end{subfigure}
~
	\begin{subfigure}[b]{0.5\textwidth}
		\begin{center}
			\begin{tikzpicture}[scale=0.7,main_node/.style={scale=0.9,circle,fill=white!80,draw,inner sep=0pt, minimum size=16pt},line width=1.2pt]
			\draw[help lines, color=gray!30, dashed] (-1.9, -0.9) grid (5.8,2.9);
			\draw[->,ultra thick] (-2,0)--(6,0) node[right]{$x_1$};
			\draw[->,ultra thick] (0,-1)--(0,3) node[above]{$x_2$};
			\node[main_node, label={[xshift=-.5cm, yshift=-0.05cm, blue]\footnotesize$(0,0)$}] (v0) at (0,0) {$v_0$};
			\node[main_node,  label={[xshift=0cm, yshift=-0.05cm, blue]\footnotesize$(2,2)$}] (v1) at (2,2) {$v_1$};
			\node[main_node,  label={[xshift=0cm, yshift=-0.05cm,blue]\footnotesize$(5,0)$}] (v2) at (5,0) {$v_2$};
			\node[main_node,  label={[xshift=0.2cm, yshift=0.0cm, blue]\footnotesize$(1.3,-0.55) $}] (v3) at (1.3,-0.55) {$v_3$};
			\end{tikzpicture}
			\caption{An embedding of the complete graph in $\R^2$.} \label{fig:completeEmbed}
		\end{center}
	\end{subfigure}	
	\caption{A realization for a complete graph in $\R^2$.} 
	\label{fig:DGPex}
\end{figure}
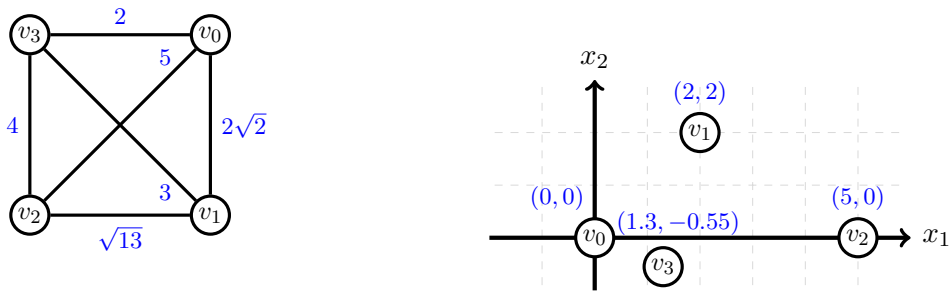

The function $x$  is called a \textit{realization} for $\graph$ or an \textit{embedding} of $\graph$ (see Figure \ref{fig:DGPex}). If $\graph$ is not connected, determining if it has a realization is equivalent to determining if its connected components have a realization so we assume $\graph$ is connected \cite{cassioli2015}. We use the Euclidean norm, however it can be any metric. We mention as well the interval DGP, where in Definition \ref{def:DGP},  the norm $\|x(u) - x(v)\|$ belongs to a given interval of weights, instead of being equal to a particular one \cite{Goncalves2017,Lavor2013}.

The DGP is $\mathcal{NP}$-Complete for $\K=1$ and $\mathcal{NP}$-Hard for $\K>1$ \cite{saxe1980embeddability}, and solution methods include nonlinear programming, semi-definite programming, and the geometric build-up methods  \cite{liberti2014euclidean,mucherinoHDR, mucherino2012DG}. In the special case where the distance between all pairs of vertices in $\graph$ are known, that is $\graph$ is complete, and we assume they yield a realization in $\R^\K$, this realization can be  found by solving a series of linear equations \cite{dgpbook}. 

In most applications, the instance is not a complete graph. The distance between some pairs of vertices is not available, and so this procedure does not apply. In such a case we would like to make use of combinatorial methods to solve the DGP, thus we must establish conditions under which the solution space of the DGP can be discretized. {\cRev For this, we assume that there exists a solution to the DGP for the instance. The solution set is then non-empty and thus is either finite or uncountable modulo translations, rotations, and reflections \cite{dgpbook, liberti2014euclidean}.} The solution to the system of linear equations is now the intersection of a line segment and a sphere  \cite{Lavor2012, dgpbook, liberti2014euclidean, mucherinoHDR}. 
The solution set is finite under the following conditions \cite{Lavor2012,dgpbook, liberti2014euclidean,mucherinoHDR, Mucherino2012}:
\BE[label=(\roman*)]
\I There is a realization for $\K$ vertices of the instance, and 
\I {\cRev For }every other vertex, $v \in \vertexset$, {\cRev there exist} edges $\{v, i\}, \{v, j\}, \{v, l\} \in \edgeset$, where the $i,j,l$ positions have already been fixed, so that we will be able to fix a position for $v$. 
\EE
These assumptions mean if the positions of $\K$ vertices {\crev are} fixed, the $(\K+1)^{\text{th}}$ vertex has at most two possible positions in relation to the previously fixed vertices {\crev since we are solving for the intersection of a line segment and a sphere. Similarly, }the $(\K+2)^{\text{th}}$ vertex has at most four possible positions relative to the previous vertices, and so on, so that the last vertex to be placed has $2^{|\vertexset| -\K}$ possible positions. Thus this search space induces a binary tree structure where each layer of the tree enumerates all possible positions for a fixed vertex, where a realization is a path in the tree from the root to a leaf \cite{mucherinoHDR}. 

If there are more edges in the graph than those that satisfy (i) and (ii), it is possible to prune positions for a vertex from the solution space. This leads to the notion of the branch-and-prune (BP) algorithm, which enumerates the possible positions of vertices one-by-one  and prunes a branch whenever there is an extra edge between the current vertex and previous vertices that is incompatible with the position \cite{dgpbook, liberti2014euclidean, mucherinoHDR}. We can think of an optimal BP search tree as the smallest search tree for a given instance \cite{Goncalves2017}. In fact, (i) and (ii) are satisfied, and the solution space is finite only if there exists a total order on the vertices satisfying the following definition \cite{Lavor2012, dgpbook, liberti2014euclidean, mucherinoHDR}: 

\begin{definition}[Discretizable DGP] \label{def:DDGP}
	Given, an integer $\K > 0$, and a simple, unweighted, undirected graph $\graph = (\vertexset,\edgeset)$, $\graph$ is an instance of Discretizable DGP if there exists a total order on $\vertexset$, $\Lorder v_0, v_1, \hdots, v_{|\vertexset|-1} \Rorder$, such that: 
	\BE[label=(\roman*)]
	\I $\graph[\{v_0, v_1, \hdots, v_{\K-1}\}]$ is a clique.
	\I For all $v_i \in \{ v_{\K}, v_{\K+1}, \hdots, v_{|\vertexset|-1} \}$, $v_i$ has
	\BE
	\I at least $\K$ adjacent predecessors
	\I a set of exactly $\K$ adjacent predecessors $\{v_{j_1}, v_{j_2}, \hdots, v_{j_\K} \}$ where $\graph[\{v_{j_1}, v_{j_2}, \hdots, v_{j_\K} \}]$ is a clique and the volume of the simplex formed by the realizations of $\{v_{j_1}, v_{j_2}, \hdots, v_{j_\K} \}$ is positive.
	\EE
	\EE	
\end{definition}
\noindent We let $\graph[\vertexset']$ be the subgraph of $\graph$ induced by $\vertexset' \subseteq \vertexset$ and a clique is a complete subgraph. We define an adjacent predecessor of a vertex $v \in \vertexset$ as $u \in \vertexset$ with $\{u,v\} \in \edgeset$ such that $u$ precedes $v$ in the order. {\crev We note that {\cRev such Discretizable DGPs (DDGPs)} are feasibility problems with no objective, wherein we wish to determine only if a vertex order exists for an instance.}

The focus of this paper is finding optimal BP trees for a class of DDGPs, namely the Discretization Vertex Ordering Problem (DVOP)\footnote{In the literature, somewhat confusingly, the Discretization Vertex Ordering Problem (DVOP) is sometimes referred to as the problem of finding an order for the DDGP \cite{cassioli2015}.}, {\crev that is, we wish to find a vertex order {\cRev with the smallest search tree} over all possible DVOP orders}. {\cRev This is an $\mathcal{NP}$-Complete problem \cite{omer2019}.} We present the DVOP in detail in Section \ref{sec:prelim}, following the convention of \cite{Lavor2012}, which distinguishes DVOP from DDGP and establishes the DVOP as a total order that does not verify the simplex-related conditions, (ii) (b), of the DDGP definition. The DDGP is $\mathcal{NP}$-Hard, and the DVOP is $\mathcal{NP}$-Complete \cite{Lavor2012, dgpbook, liberti2014euclidean}. However, if $\K$ is fixed, there exists a greedy algorithm to solve DVOP, given all possible initial cliques. Thus DVOP with fixed $\K$ is polynomial \cite{Lavor2012,liberti2014euclidean}.

The rest of the paper is organized as follows. In Section \ref{sec:prelim}, we present the DVOP, and explain in detail the problem of finding an optimal discretization order, $\MinD$. We also review two existing integer programming (IP) formulations, and a branch-and-cut procedure from the literature. In Section \ref{sec:models}, we introduce a novel IP formulation and three novel constraint programming (CP) formulations for  $\MinD$. In Section \ref{sec:hybrid}, we present two hybrid IP-CP decomposition algorithms, as well as some valid inequalities for the problem. Finally, in Section \ref{sec:results}, we present a computational study.

We note that an overview of our paper, namely the models/methods from the literature as well as our proposed models/methods are provided in Table \ref{table:MinDsummary} of Appendix \ref{app:summary}.

\section{Preliminaries} 
\label{sec:prelim}
\medskip
\subsection{Notation}
All sets are denoted calligraphically. Let $\graph = (\vertexset, \edgeset)$ be an undirected graph, where $\vertexset$ is the set of vertices and $\edgeset$ is the set of edges. The \emph{{\crev adjacency} matrix} of $\graph$ is denoted by $\adjmatrix$, i.e.,  $\adjmatrix_{v,u} = 1 $ if and only if edge $\{u,v\} \in \edgeset$. We define a directed graph $\digraph = (\vertexset, \arcset)$, where $\arcset$ is the set of directed arcs in $\digraph$, i.e., $\arcset = \{(u,v) \cup (v,u) : \{u,v\} \in \edgeset\}$. We adopt the convention of denoting an undirected edge, $\{u,v\}$, and a directed arc, $(u,v)$. Denote the \emph{neighbourhood} of a vertex $\vertind$ as $\neighb(\vertind)$,  i.e., $\neighb(\vertind) = \{ u \in \vertexset : \{u,v\} \in \edgeset\}$, thus $\vertind \notin \neighb(\vertind)$ and the \emph{degree} of $\vertind$ as $\degr(\vertind) = |\neighb(\vertind)|$. We let $\graph[\vertexset'] = (\vertexset', \edgeset')$ be the subgraph of $\graph$ \emph{induced} by $\vertexset' \subseteq \vertexset$, and thus $\edgeset' = \{\{u,v\}  \in \edgeset : u,v \in \vertexset' \}$. A \emph{clique}, $\clique$, in $\graph$ is a set of vertices $\{v_1, v_2, \hdots, v_{|\clique|}\} \subseteq \vertexset$ such that $\{v_i, v_j\} \in \edgeset$ for all $v_i, v_j \in \clique$ such that $v_i \neq v_j$. We define a directed cycle, $C$ as a subgraph of $\digraph$, $C = (\vertexset^C, \arcset^C)$, where $C$ forms a path where the first node is the same as the last node. We define an \emph{adjacent predecessor} of a vertex $v \in \vertexset$ as $u \in \vertexset$ with $\{u,v\} \in \edgeset$ such that $u$ precedes $v$ in a vertex order.

For $a, b \in \Zplus$, $a \leq b$,  we use the notation $[a] = \{0, 1, \hdots, {\crev a}\}$ and $[a, b] = \{a, a+1, \hdots, b\}$. If $a>b$, then  $[a,b]=\emptyset$, similarly if $a<0$, then  $[a]=\emptyset$. We use $\indicate(\cdot)$ as the \emph{indicator function}, which evaluates to $1$ if the boolean expression it operates on is true and $0$ if it is false.

Indices follow these conventions: indices start at $0$, so that the possible positions of a vertex order are $[|\vertexset|-1]$. We let $|\vertexset| = n$, and use $|\vertexset|$ in relation to vertices and $n$ in relation to ranks of a vertex order. 

\subsection{Problem Definition}

The \textit{Discretization Vertex Order Problem} (DVOP) \cite{Lavor2012} is the search for a total order of the vertices of a simple, connected, undirected graph $\graph$, given an integer dimension $\K$, that satisfies the following:
\BE[label=(\roman*)]
\item the first $\K$ vertices in the order form a clique in the input graph, $\graph$,  and
\item the following vertices each have at least $\K$ adjacent vertices in $\graph$ as predecessors in the order. 
\EE
We refer to a total order that satisfies (i) and (ii) as a \textit{DVOP order}, in this case we say the instance $\instance$ is \textit{feasible}, otherwise it is \textit{infeasible}.

In order to characterize the DVOP we introduce a general function $\rankfunc: \vertexset \to [n-1]$. Let $\LOrankSet$ be the set of $\rankfunc(\cdot)$ that give a  linear ordering, i.e., $\rankfunc(\cdot)$ is bijective  \cite{coudert2016}. Let $\DVOPrankSet \subseteq \LOrankSet$ be the set of $\rankfunc(\cdot)$ that give a linear ordering and satisfies (i) and (ii). We say $\rankfunc(\cdot)$ \emph{characterizes} a DVOP order if $\rankfunc(\cdot) \in \DVOPrankSet$. Note that (ii) implies that the vertex at rank $\K$ must be adjacent to all the vertices at ranks $[\K-1]$. Combined with (i) this implies that the first $\K+1$ vertices in the DVOP order induce a clique in $\graph$, we call this clique the first or initial clique. Formally, we have:
\BE[label=(\roman*)]
\item $\graph[\{\vertind \in \vertexset : \rankfunc(v) \leq \K\}]$ is a clique, and
\item $|\{ u \in \neighb(\vertind) : \rankfunc(u) \leq \rankfunc(\vertind)-1\}| \geq K$ for all $\vertind \in \vertexset$ with $rank(\vertind) \geq \K +1$.
\EE

\begin{figure}[h]
	\centering
	\begin{tikzpicture}[scale=0.7,main_node/.style={scale=0.9,circle,fill=white!80,draw,inner sep=0pt, minimum size=16pt},
	line width=1.2pt]
	\node[main_node] (v0) at (3,3) {$v_0$};
	\node[main_node] (v1) at (4.5,1.5) {$v_1$};
	\node[main_node] (v2) at (3,0) {$v_2$};
	\node[main_node] (v3) at (0,0) {$v_3$};
	\node[main_node] (v4) at (-1.5,1.5) {$v_4$};
	\node[main_node] (v5) at (0, 3) {$v_5$};
	\draw[-] (v0) -- (v1) -- (v2) -- (v3) -- (v4) -- (v2) -- (v0) -- (v5) -- (v3) -- (v1) -- (v5) -- (v2);
	\end{tikzpicture}
	\caption{\label{fig:DVOPgraph} A graph instance which is feasible for DVOP with $\K =2$.}
\end{figure}
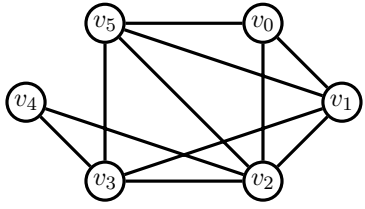

\begin{ex}
	Figure \ref{fig:DVOPgraph} shows a graph for which $\Lorder v_0, v_1, v_2, v_3, v_4, v_5\Rorder$ is a DVOP order for $\K =2$. That is for some $\rankfunc(\cdot) \in \DVOPrankSet$: $\rankfunc(v_0) = 0$, $\rankfunc(v_1) = 1$, $\rankfunc(v_2) = 2$ and so on. We can see the vertices $v_0, \text{ and } v_1$ are adjacent and thus form a 2-clique, the vertices $v_2, v_3,$  and $v_4$ each have two adjacent predecessors and $v_5$ has four adjacent predecessors in the order. Note that the first $\K +1$ vertices in the order, namely $v_0, v_1, v_2$, form a clique since again the vertex at position $\K$ must be adjacent to the $\K$ previous vertices.  
\end{ex}

However, there is no DVOP order for the $\K=3$ case for the graph in Figure \ref{fig:DVOPgraph}, since $v_4$ is not in a $4$-clique and cannot have $3$ adjacent predecessors as it has $2$ neighbours in $\graph$.


Recall, from Section \ref{sec:intro}, the solution space of DDGPs can be represented as a binary tree structure that may be searched using the branch and prune (BP) algorithm. Here the vertex order dictates the manner in which we search over the continuous $\R^\K$ space to find a solution to the DGP. Given a DVOP order, the BP algorithm solves the DDGP by fixing the coordinates of the first $\K$ vertices in the order and enumerating the possible realizations of the remaining vertices. In the BP search tree, branching on a vertex with exactly $\K$ predecessors in the order yields at most two child nodes which are called \textit{double vertices}. Otherwise if the vertex has more than $K$ predecessors it has at most one child in the BP tree and is a \textit{single vertex} \cite{Omer2017}.
In order to characterize the notion of a double vertex we define  function which we call $\doublefunc(\cdot)$ as
\[\doublefunc(\vertind) = 
\begin{cases}
1, & \text{if } \vertind \text{ has exactly } \K \text{ adjacent predecessors in the order} \\
0, & \text{otherwise}
\end{cases} \quad \forall \ \vertind \in \vertexset.
\]
The number of double vertices can be used as a measure of the size of the BP tree and defined by the recursion,
\[ \numdoub(\rankind)= \begin{cases} 
1, & \rankind \in [\K-1]\\
(\doublefunc(\rankfunc^{-1}(\rankind)) + 1)\cdot  \numdoub(\rankind -1 ), &  \rankind \in [\K, n-1] 
\end{cases}
\]
where $\numdoub(\rankind)$ is the number of nodes at level $\rankind$ of the BP tree \cite{Omer2017}. Notice that because each level of the tree gives all positions of a single vertex in the order, the levels of the tree are equivalent to the positions of the order. The maximum number of nodes in the BP tree is
the sum of $\numdoub(\cdot)$ over all positions. Note that we say the \textit{maximum} number of nodes in the BP tree, not the number of nodes, since $\graph$ may have  extra edges that allow us to prune positions and reduce the number of nodes in the tree.

Given this measure of BP tree size, \cite{Omer2017} defined two optimization problems: $\MinD$  and  $\MinN$, both of which try to find an optimal DVOP order but under different objective functions. $\MinD$ seeks to minimize the number of nodes which are doubles, while $\MinN$ seeks to minimize the maximum number of nodes in the BP tree and thus requires both a minimum number of double vertices and to fix the positions of these vertices as close to the end of the order as possible so as to have the smallest effect on  tree size, due to doubling the number of nodes. 
	\begin{figure}[h] 
	\begin{subfigure}{0.5\textwidth}
		\centering
		\begin{tikzpicture}[main_node/.style={scale=0.88, circle,fill=white!80,draw,inner sep=0pt, minimum size=16pt},
		line width=1.2pt]
		\node[main_node] (a) at (0,0) {$v_0$};
		\node[main_node] (b) at (0, -0.75) {$v_1$};
		\node[main_node] (c) at (0,-1.5) {$v_2$};
		\node[main_node] (d1) at (-2,-2) {$v_3$};
		\node[main_node] (d2) at (2,-2) {$v_3$};
		\node[main_node] (e1) at (-3,-2.5) {$v_4$};
		\node[main_node] (e2) at (-1,-2.5) {$v_4$};
		\node[main_node] (e3) at (1, -2.5) {$v_4$};
		\node[main_node] (e4) at (3,-2.5) {$v_4$};
		\node[main_node] (f1) at (-3.5, -3.5) {$v_5$};
		\node[main_node] (f2) at (-2.5, -3.5) {$v_5$};
		\node[main_node] (f3) at (-1.5, -3.5) {$v_5$};
		\node[main_node] (f4) at (-0.5, -3.5) {$v_5$};
		\node[main_node] (f5) at (3.5, -3.5) {$v_5$};
		\node[main_node] (f6) at (2.5, -3.5) {$v_5$};
		\node[main_node] (f7) at (1.5, -3.5) {$v_5$};
		\node[main_node] (f8) at (0.5, -3.5) {$v_5$};
		\draw[-] (a) -- (b) -- (c) -- (d1) -- (e1)--(f1);
		\draw[-] (e1) -- (f2);
		\draw[-] (d1) -- (e2) -- (f3);	
		\draw[-] (e2) -- (f4);
		\draw[-] (c) -- (d2) -- (e4) -- (f5);	
		\draw[-] (d2) -- (e3) -- (f8);	
		\draw[-] (e3) -- (f7);
		\draw[-] (e4) -- (f6);
		\end{tikzpicture}
		\caption{\label{fig:DVOPBP1}  A BP tree for DVOP order $\Lorder v_0, v_1, v_2, v_3, v_4, v_5\Rorder$ of the graph in Figure \ref{fig:DVOPgraph}.}
	\end{subfigure}
	~
	\begin{subfigure}{0.5\textwidth}
		\centering
		\begin{tikzpicture}[main_node/.style={scale=0.88,  circle,fill=white!80,draw,inner sep=0pt, minimum size=16pt},
		line width=1.2pt]
		\node[main_node] (a) at (0,0) {$v_3$};
		\node[main_node] (b) at (0, -0.75) {$v_5$};
		\node[main_node] (c) at (0,-1.5) {$v_2$};
		\node[main_node] (d1) at (-2,-2) {$v_1$};
		\node[main_node] (d2) at (2,-2) {$v_1$};
		\node[main_node] (e1) at (-2,-2.75) {$v_0$};
		\node[main_node] (e2) at (2,-2.75) {$v_0$};
		\node[main_node] (f1) at (-2,-3.5) {$v_4$};
		\node[main_node] (f2) at (2, -3.5) {$v_4$};
		\draw[-] (a) -- (b) -- (c) -- (d1) -- (e1)--(f1);
		\draw[-] (c) -- (d2) -- (e2)--(f2);
		\end{tikzpicture}
		\caption{\label{fig:DVOPBP2}  A BP tree for DVOP order $\Lorder v_3, v_5, v_2, v_1, v_0, v_4\Rorder$ of the graph in Figure \ref{fig:DVOPgraph}.}
	\end{subfigure}
	\caption{BP trees for two DVOP orders of the graph in Figure \ref{fig:DVOPgraph}.}
\end{figure}
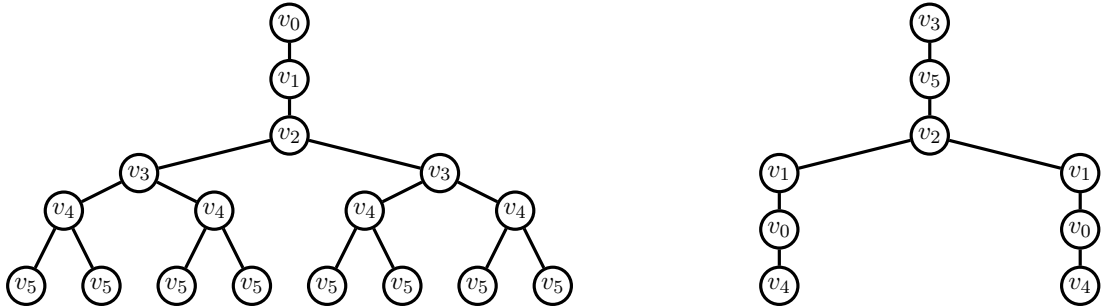
\begin{ex}
	The order previously given for the graph in  Figure \ref{fig:DVOPgraph},  $\Lorder v_0, v_1, v_2, v_3, v_4, v_5\Rorder$, has three double vertices, $v_2, v_3, v_4$, which all have exactly $\K$ adjacent predecessors in the order. However, a DVOP order that minimizes the number of doubles for this instance is $\Lorder v_3, v_5, v_2, v_1, v_0, v_4\Rorder$, which has two double vertices, $v_2, v_4$. Figure \ref{fig:DVOPBP1} shows the BP tree for the first order, with three doubles, it has at most $17$ nodes in total.
	Recall, that each node of the BP tree is a possible coordinate position for that vertex in  $\R^\K$. Figure \ref{fig:DVOPBP2} shows the BP tree for the second order, with two doubles, and at most $9$ nodes  almost half as many as the BP tree of first order. We also note that the tree in Figure \ref{fig:DVOPBP1} has width $8$, whereas the tree in Figure \ref{fig:DVOPBP2} has width $2$. Thus fewer doubles reduced the size of the BP tree both with respect to the width and the number of nodes, giving a clear motivation to solve the aforementioned optimization problems. 
\end{ex}

The general framework for $\MinD$ is 
\bsubeq \label{form:generalMinD}
\begin{alignat}{3}
\min & \ \sum_{\rankind \in [n-1]} \doublefunc(\rankfunc^{-1}(\rankind)) \\
\text{s.t.} & \ \rankfunc(\cdot) \in \DVOPrankSet \\
&\hspace*{-0.12cm} \sum_{u \in \neighb(\vertind)} \hspace*{-0.12cm} \indicate{(\rankfunc(u) \leq \rankfunc(\vertind)-1)}  &  \begin{cases}
= \K & \hspace*{-0.26cm} \Rightarrow \doublefunc(\vertind) = 1 \\
\geq \K+1 & \hspace*{-0.26cm} \Rightarrow   \doublefunc(\vertind) = 0 
\end{cases} \ \,
\forall \ \vertind \in \vertexset : \rankfunc(\vertind) \geq K \label{generalMinD:linking}
\end{alignat}
\esubeq
Constraints \eqref{generalMinD:linking} can be simplified since at optimality the objective will ensure there are the fewest doubles possible, thus we need only to enforce that  $\vertind$ has at least $\K+1$ adjacent predecessors if it is not a double. Otherwise, we need to enforce that $\vertind$ has a least $\K$ adjacent predecessors. Thus, we continue with the following reduced version of \eqref{form:generalMinD}:
\bsubeq \label{form:MinDouble}
\begin{alignat}{2}
\min \ & \sum_{\rankind \in [n-1]} \doublefunc(\rankfunc^{-1}(\rankind)) \\
\text{s.t.} \ & \rankfunc(\cdot) \in \DVOPrankSet  \label{rankConstr} \\
&\sum_{u \in \neighb(\vertind)} \indicate{(\rankfunc(u) \leq \rankfunc(\vertind)-1)}  \geq \K+1 - \doublefunc(\vertind) 
&& 
\quad \forall \ \vertind \in \vertexset \text{ s.t. } \rankfunc(\vertind) \geq K \label{predConstr}
\end{alignat}
\esubeq
Note that constraints \eqref{predConstr} imply that every vertex $\vertind$ with $\rankfunc(\vertind) \geq K$ has at least $\K$ adjacent predecessors since $\doublefunc(\vertind) \in \{0,1\}$. Moreover, if $\vertind$ also has $\leq \K$, i.e., $ = \K$, adjacent predecessors, then constraints \eqref{predConstr} makes it a double. Otherwise, constraint \eqref{predConstr} for $\vertind$ becomes redundant, so the value of $\doublefunc(\vertind)$ can be either $0$ or $1$. But due to the objective, which is minimizing the number of doubles, we will have $\doublefunc(\vertind) = 0$, which is what is desired by definition of $\doublefunc(\cdot)$.



We are able to extend the general formulation for $\MinD$ to one for $\MinN$ by simply changing the objective function to
$\sum_{\vertind \in \vertexset} \numdoub(\rankfunc(\vertind))$
since the number of nodes in the tree, $\numdoub(\rankfunc(\vertind))$, at each rank follows by definition from whether the vertex at that position in the order is a double or not, $\doublefunc(\vertind)$.
Furthermore, we consider the multi-objective case of minimizing both the number of doubles and the maximum number of nodes. We define the multi-objective problem as follows.
\begin{equation}
\label{form:generalMultiObj}
\min \  \left \{  \ \sum_{\vertind \in \vertexset} \numdoub(\rankfunc(\vertind)) , \sum_{\rankind \in [n-1]} \doublefunc(\rankfunc^{-1}(\rankind)) \right \} \ \text{s.t.} \  \eqref{rankConstr}-\eqref{predConstr}
\end{equation}

\begin{figure}[h]
\centering
\begin{tikzpicture}[scale = 0.8]
  \pgfplotsset{
   width=10cm,height=4.5cm
  }

  \begin{axis}[axis x line=bottom,
    axis y line=left,
    axis on top,
    domain=0:1,
    ytick distance=1,
    ymin=1.75,ymax=3.25,
    xmin=11,xmax=25,
    xlabel=$\MinN$ objective,
    ylabel=$\MinD$ objective,
  ]
    \addplot[only marks, mark options={scale=2}]
    coordinates{ 
      (24,3)
      (14,2)
      (20,3)
      (16,2)
    }; 
    
    \addplot[color=blue,mark=square*,mark options={scale=2}] coordinates {
	(12,2)
};

    \label{plot_one}

    \addlegendimage{/pgfplots/refstyle=plot_one}
  \end{axis}
\end{tikzpicture}  
\caption{Set of feasible solutions in the objective space of the multi-objective problem \eqref{form:generalMultiObj} for the example in Figure \ref{fig:DVOPgraph}, where the non-dominated point is marked with a square.}
\label{fig:pareto}
\end{figure}
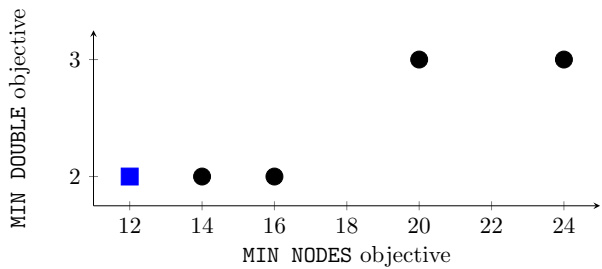

\begin{ex}
Consider the graph in Figure \ref{fig:DVOPgraph}. For this example, the multi-objective problem has 180 feasible (DVOP) solutions in the decision space, which yield the set of five images in the objective space as shown in Figure \ref{fig:pareto}. Its Pareto frontier consists of the square point. 
\end{ex}
{\crev We conjecture that the $\MinN-\MinD$ multi-objective problem has a feasible ideal point. If Conjecture \ref{conj:pareto} is true we can solve $\MinD$ and retrieve from its solution an optimal solution to $\MinN$. 
\begin{conj}
\label{conj:pareto}
	The Pareto frontier of the $\MinN-\MinD$ multi-objective problem \eqref{form:generalMultiObj} consists of a single non-dominated point.
\end{conj}
}

For the remainder of this paper, we study solution methods to $\MinD$, but these methods can be easily extended to $\MinN$ {\crev by a couple simple modifications to the proposed models {\cRev (more explicitly by defining additional variables for $\numdoub(\cdot)$  values and linking them with the variables describing the DVOP order, leaving the subproblems and cuts for the decompositions models intact).}} {\cRev An example of an IP model for $\MinN$ can be found in Appendix  \ref{appendix:IPminNodes}}.

\subsection{Existing Mathematical Models} \label{sec:OGmodels}

Prior to this work, \cite{Omer2017} present two IP formulations and one branch-and-cut procedure for $\MinD$. These are summarized below, while full details can be found in Appendix \ref{appendix:existing}.
\BI
\I The cycles formulation $\OGcycles$: They introduce three sets of binary variables to  indicate precedence between vertices, the double vertices, and the initial clique, respectively. The constraints break $2$-cycles and $3$-cycles in the precedence variables, select the first clique, and link precedence and the initial clique variables to the doubles.
\I The rank formulation $\OGranks$: They replace the cycle breaking constraints in $\OGcycles$ with an adaptation of the \cite{Miller1960} formulation for the travelling salesman problem, for which they introduce rank variables, and link the precedence variables with the ranks.
\I Cycle cut generation $\CCG$: They break cycles in the precedence variables iteratively within a branch-and-cut procedure.
\EI

\section{Mathematical Models} \label{sec:models}
{\crev To formulate the mathematical models we define binary variables $\doublevar_\rankind= 1$ if the vertex at position $\rankind \in [n-1]$ is a double, $0$ otherwise. It is also possible to define this variable with respect to the vertex index $\vertind \in \vertexset$ as in the existing IP formulations. When the rank-based double variables are used, we fix the first $\K$ positions to be single vertices, as by definition they cannot be doubles. We also fix position $\K$ to be a double, since we have an initial clique of size $\K + 1$, which implies that the $(\K+1)^{\text{th}}$ vertex will always be adjacent to exactly $\K$ vertices and thus is always double. For clarity we write these as constraints in the formulations, however they are implemented as bounds and they may be omitted by projecting out the $\doublevar_\rankind$ for $\rankind \in [\K]$ and adding one to the objective function. }
\medskip

\subsection{Integer Programming Formulation}
Our first IP formulation extends the DVOP formulation presented by \cite{Lavor2012} to incorporate constraints for \MinD. Define binary variables $ \vrvar_{\vertind \rankind}=1$ if the vertex $\vertind \in \vertexset$ is at rank $\rankind \in [n-1]$, $0$ otherwise. We also introduce binary indicator variables $\ipindicator_{\vertind\rankind}$ for all $\vertind \in \vertexset, \ \rankind \in [n-1]$ to express logical constraints. Then, the formulation is:
\bsubeq \label{form:MDIPvr}
\begin{alignat}{2}
\IP: \min \ &  \sum_{\rankind \in [n-1]} \doublevar_\rankind && \label{MDIPvr:obj} \\
\text{s.t.} \ &\sum_{\rankind \in [n-1]}\vrvar_{\vertind\rankind} = 1 && \forall \  \vertind \in \vertexset \label{MDIPvr:vertex1rank}\\
&\sum_{\vertind \in V} \vrvar_{\vertind\rankind} = 1 && \forall \ \rankind \in [n-1] \label{MDIPvr:rank1vertex}\\
& \sum_{u \in \neighb(\vertind)} \sum_{j \in [\rankind -1]} \vrvar_{uj} \geq \rankind\vrvar_{\vertind\rankind} \quad &&  \forall \ \vertind \in \vertexset, \rankind \in [1,\K] \label{MDIPvr:initalClique}\\
&\sum_{u \in \neighb(\vertind)} \sum_{j \in [\rankind -1]} \vrvar_{uj} \geq \K \vrvar_{\vertind \rankind} &&  \forall \ \vertind \in \vertexset, \rankind \in [\K+1, n-1] \label{MDIPvr:Kpred}\\
& \doublevar_\rankind = 0 && \forall \ \rankind \in [\K-1] \label{MDIPvr:cliqueDoubles} \\
&  \doublevar_\K =1 && \label{MDIPvr:firstDouble}\\
& \sum_{u \in \neighb(\vertvar)} \sum_{j \in [\rankind -1]}\vrvar_{uj}    \geq (K+1)\ipindicator_{\vertind\rankind} \quad &&   \forall \ \vertind \in \vertexset, \rankind \in [\K, n-1] \label{MDIPvr:logical1}\\ 
&  \vrvar_{\vertind\rankind} -  \doublevar_\rankind    \leq \ipindicator_{\vertind\rankind} \quad &&   \forall \ \vertind \in \vertexset, \rankind \in [\K, n-1]\label{MDIPvr:logical2}\\ 
&\vrvar \in \{0,1\}^{|\vertexset| \times n},  \doublevar   \in \{0,1\}^n, \ipindicator \in  \{0,1\}^{|\vertexset| \times n} \quad && 
\end{alignat}
\esubeq
Constraints \eqref{MDIPvr:vertex1rank} and \eqref{MDIPvr:rank1vertex} ensure that there is a bijection from the vertices to the ranks, that is each vertex has exactly one rank and vice versa. Constraints \eqref{MDIPvr:initalClique} ensure that we have an initial clique of size $\K +1$, since each vertex in the clique will be adjacent to all its predecessors. Constraints \eqref{MDIPvr:Kpred} ensure that each vertex after the initial clique has at least $\K$ adjacent predecessors. {\crev Constraints \eqref{MDIPvr:cliqueDoubles} and \eqref{MDIPvr:firstDouble} fix the double value of positions $[\K]$.}
Constraints \eqref{MDIPvr:logical1} and  \eqref{MDIPvr:logical2} link the $\vrvar_{\vertind\rankind} $ and $ \doublevar_\vertind$ variables, where we want to enforce that if vertex $\vertind$ is a non-double, then for any rank $\rankind$, either $\vertind$ is not in position $\rankind$ or it has at least $\K + 1$ adjacent predecessors, that is
$    
\doublevar_\rankind = 0 \land  \vrvar_{\vertind\rankind}  = 1 \implies  \sum_{u \in \neighb(\vertvar)} \sum_{j \in [\rankind -1]} \vrvar_{uj} \geq K+1.
$
The left-hand-side can be rewritten giving
$
\vrvar_{\vertind\rankind} - \doublevar_\rankind \geq 1 \implies  \sum_{u \in \neighb(\vertvar)} \sum_{j \in [\rankind -1]} \vrvar_{uj} \geq K+1.
$
We then use the indicator variable $\ipindicator_{\vertind\rankind}$ taking the value of $1$ if the left-hand-side of the inequality holds   to rewrite the implication as constraints \eqref{MDIPvr:logical1} and  \eqref{MDIPvr:logical2}. 

For any feasible solution to $\IP$, $\rankfunc(v) = \sum_{r \in [\rankind -1]} \rankind \vrvar_{vr}, \forall \ \vertind \in \vertexset $ characterizes the order, while $\doublefunc(\rankfunc^{-1}(\rankind))= \doublevar_\rankind, \forall \ \rankind \in [n-1]$ provides the information about the doubles.

\subsection{Constraint Programming Formulations}
The effectiveness of Constraint Programming (CP) for solving permutation-based problems can be leveraged to solve $\MinD$. We present three CP formulations for $\MinD$.

We begin with a natural translation of $\IP$ into CP to define the primal CP formulation. Define integer variables $\rankvar_\vertind \in [n-1]$ denoting the rank of  vertex $\vertind \in \vertexset$. 
\bsubeq\label{form:MDCPrank}
\begin{alignat}{2}
\hspace*{-0.6cm}  \rankCP:  & \min \ \sum_{\vertind \in \vertexset} \doublevar_\vertind + 1&& \label{MDCPrank:obj}\\
& \text{s.t.} \ \text{AllDifferent}(\rankvar_0, \rankvar_1,..., \rankvar_{n-1}) \label{MDCPrank:allDiff}\\ 
& \rankvar_i \geq \K+1 \lor \rankvar_j \geq \K+1  &&\forall \ i,j  \in \vertexset : i\neq j, \{i,j\} \notin \edgeset \label{MDCPrank:nonEdges} \\
& \rankvar_v \geq \K+1 \Rightarrow \hspace*{-0.25cm} \sum_{u \in \neighb(\vertind)} \hspace*{-0.1cm}  \indicate(\rankvar_u \leq \rankvar_\vertind -1 ) \geq \K +   (1-\doublevar_{\vertind}) \ \ && \forall \ \vertind \in \vertexset \label{MDCPrank:pred}\\ 
&  \doublevar \in \{0,1\}^n, \rankvar \in [n-1]^{|\vertexset|}
\end{alignat}
\esubeq
Objective \eqref{MDCPrank:obj} minimizes the number of vertices which are double and adds one for the vertex in position $\K$ which is always a double by definition. Constraint  \eqref{MDCPrank:allDiff} is the CP global constraint AllDifferent, which forces every variable in the set to have a unique value, because this constraint acts on $n$ variables, all of which have the same domain which also has $n$ values, \eqref{MDCPvertex:allDiff}  ensures that each vertex and rank map one to one. Constraints \eqref{MDCPrank:nonEdges} constrain that we have an initial clique of size $\K+1$, by having only pairs of vertices which are adjacent in the first $K+1$ ranks. Logical constraints \eqref{MDCPrank:pred} ensure that all vertices with ranks outside the initial clique have $\K$ adjacent predecessors if they are double and at least $\K+1 $ adjacent predecessors otherwise. These constraints use CP cardinality clause constraints, which specify a particular  number of boolean variables must be true. However, constraints \eqref{MDCPrank:pred} do not use boolean variables directly, instead each $(\rankvar_u \leq \rankvar_\vertind -1 )$ is a predicate taking a boolean value. 

The function which characterizes the order for $\rankCP$ is $\rankfunc(\vertind) = \rankvar_\vertind,  \forall \ \vertind \in \vertexset$ and the double function is $\doublefunc(\vertind) = \doublevar_\vertind,  \forall \ \vertind \in \vertexset$.

The second CP formulation for $\MinD$ is the dual formulation to \eqref{form:MDCPrank}. We define integer variables $\vertvar_{\rankind}$ equal to the vertex index at rank $\rankind \in [n-1]$. As such, the value of these dual variables are equivalent to the primal variables in \eqref{form:MDCPrank} \cite{Smith01dualmodels}. {\crev We note many of the constraints of this formulation depend on the entries of the adjacency matrix, $\adjmatrix$.}
\bsubeq \label{form:MDCPvertex}
\begin{alignat}{2}
\vertexCP: \min \ & \sum_{\rankind\in [n-1]} \doublevar_\rankind&& \label{MDCPvertex:obj} \\
\text{s.t.} \ & \text{AllDifferent}(\vertvar_0, \vertvar_1,..., \vertvar_{n-1}) \label{MDCPvertex:allDiff}\\ 
& \adjmatrix_{\vertvar_i,\vertvar_j} = 1  &&  \forall \ i \in [ \K-1], \ j \in [i+1, \K] \label{MDCPvertex:firstClique}\\
& \sum_{j=0}^{\rankind-1} \adjmatrix_{\vertvar_\rankind,\vertvar_j} \geq \K +  (1-\doublevar_{\rankind}) \qquad  &&  \forall \ \rankind \in [\K+1, n-1]  \label{MDCPvertex:pred}\\
& \doublevar_\rankind = 0 && \forall \ \rankind \in [\K-1] \label{MDCPvertex:doubleClique}\\
&  \doublevar_\K =1 && \label{MDCPvertex:doubleK}\\
&  \doublevar \in \{0,1\}^n, \vertvar \in [|\vertexset|-1] ^n
\end{alignat}
\esubeq
Objective \eqref{MDCPvertex:obj} minimizes the number of ranks with double variables. Constraint \eqref{MDCPvertex:allDiff} enforces the one to one mapping of ranks to vertices.  Constraints \eqref{MDCPvertex:firstClique} ensure that we have an initial clique of size $\K+ 1$ by  forcing all vertices in the first $\K+ 1$ ranks to all be pairwise adjacent. These constraints use the so-called element constraints, which allow the use of variables as array indices in CP. Constraints \eqref{MDCPvertex:pred} force all vertices in ranks greater than $\K$ to have at least $\K$ adjacent predecessors if the are double and at least $\K+1$ adjacent predecessors if they are not. Constraints \eqref{MDCPvertex:doubleClique} and \eqref{MDCPvertex:firstClique} are the same as constraints \eqref{MDIPvr:cliqueDoubles} and \eqref{MDIPvr:firstDouble} which allow fixing of the variables whose double values are known. 

For $\vertexCP$, $\rankfunc^{-1}(\rankind)= \vertvar_\rankind, \forall \ \rankind \in [n-1] $ characterizes the DVOP order and the double function is $\doublefunc(\rankfunc^{-1}(\rankind))= \doublevar_\rankind, \forall \ \rankind \in [n-1].$
The third CP formulation uses both primal and dual variables and constraints, creating one larger combined formulation which can leverage redundant constraints to make stronger inferences in the CP search.
\bsubeq \label{form:MDCPcomb}
\begin{alignat}{2}
\hspace*{-0.2cm} \combCP: \min \ &  \sum_{\rankind\in [n-1]} \doublevar_\rankind&& \label{MDCPcomb:obj}\\
\text{s.t.} \ & \text{inverse}(\vertvar, \rankvar) \label{MDCPcomb:inverse}\\ 
& \adjmatrix_{\vertvar_i,\vertvar_j} = 1  &&  \forall \ i \in [\K-1], \ j \in [i+1, \K] \label{MDCPcomb:vClique} \\
& \sum_{j=0}^{\rankind-1} \adjmatrix_{\vertvar_\rankind,\vertvar_j} \geq \K +  (1-\doublevar_{\rankind}) \quad  &&  \forall \ \rankind \in [\K+1, n-1]  \label{MDCPcomb:vertexPred} \\
& \doublevar_\rankind = 0 && \forall \ \rankind \in [\K-1] \label{MDCPcomb:doubleClique} \\
&  \doublevar_\K =1 && \label{MDCPcomb:doubleK}  \\
&\rankvar_i \geq \K+1 \lor \rankvar_j \geq \K+1 \quad && \forall \ i,j \in \vertexset : i\neq j, \{i,j\} \notin \edgeset \label{MDCPcomb:rClique} \\
& \rankvar_v \geq \K+1  \Rightarrow \hspace*{-0.15cm} \sum_{u \in \neighb(\vertind)} \hspace*{-0.1cm} \indicate(\rankvar_u \leq \rankvar_\vertind -1 ) \geq \K \quad && \forall \ \vertind \in \vertexset \label{MDCPcomb:rankPred}  \\ 
&  \doublevar \in \{0,1\}^n, \vertvar \in [|\vertexset|-1]^n, \rankvar \in [n-1]^{|\vertexset|}
\end{alignat}
\esubeq
Objective \eqref{MDCPcomb:obj} minimizes the number of ranks which have double vertices. Constraint \eqref{MDCPcomb:inverse} channels the primal and dual variables by enforcing the relationship $(\rankvar_i = j) \equiv (\vertvar_j=i)$. This inverse constraint allows for the elimination of the AllDifferent constraints because the domains of both types of variables are the same. Constraints \eqref{MDCPcomb:rankPred} do not include the double variables as we have chosen to index the doubles in the natural way using ranks instead of vertices. All the other constraints are as defined in formulations \eqref{form:MDCPrank} and \eqref{form:MDCPvertex}.

For $\combCP$, both $\rankfunc(\vertind) =  \rankvar_\vertind, \forall \ v \in \vertexset \text{ and } \rankfunc^{-1}(\rankind)= \vertvar_\rankind, \forall \ \rankind \in [n-1]$ characterize the DVOP order, while $\doublefunc(\rankfunc^{-1}(\rankind))= \doublevar_\rankind,  \forall \ \rankind \in [n-1]$ is the double function.

\section{Hybrid Decomposition Approaches and Valid Inequalities} \label{sec:hybrid}
\medskip
\subsection{A Naive Decomposition}

We present a decomposition of $\vertexCP$, where the master problem will be solved using IP. The master problem, $\naiveMP$, provided in \eqref{naiveMP}, fixes which positions of the order have double vertices. As such, since the double variables are rank-indexed we are able to fix the values of $\K + 1$ variables (which is not possible with vertex-indexed double variables).
\begin{equation}
\label{naiveMP}
\naiveMP: \ \min \left \{  \sum_{\rankind \in [n-1]} \doublevar_\rankind  :
 \doublevar_\rankind = 0 \ \forall \ \rankind \in [\K-1], \ 
\doublevar_\K =1, \ 
 \doublevar \in  \{0,1\}^n  \right\}
\end{equation}

Given a feasible solution to the master problem, $\hat{\doublevar}$, the subproblem then tries to build a DVOP order with doubles in the correct location.
\bsubeq
\begin{alignat}{2}
\naiveSP: \ & \text{AllDifferent}(\vertvar_0, \vertvar_1,..., \vertvar_{n-1}) \quad \\ 
& \adjmatrix_{\vertvar_i,\vertvar_j} = 1  &&  \forall \ i \in [\K-1], \ j \in [i+1, \K], \\
& \sum_{j=0}^{\rankind-1} \adjmatrix_{\vertvar_\rankind,\vertvar_j} \geq \K  &&  \forall \ \rankind \in [\K+1,n-1] \text{ with } \hat{ \doublevar}_\rankind = 1\\
& \sum_{j=0}^{\rankind-1} \adjmatrix_{\vertvar_\rankind,\vertvar_j} \geq \K+1  &&  \forall \ \rankind \in [\K+1,n-1] \text{ with } \hat{ \doublevar}_\rankind = 0 \\
&\vertvar \in [|\vertexset|-1]^n
\end{alignat}
\esubeq

A feasible subproblem means we have found an optimal solution since $\naiveMP$ will give the smallest number of doubles and $\naiveSP$ ensures that there is in fact an order with this number of doubles. If $\naiveSP$ is infeasible, that means there is no order with the $\hat{\doublevar}$ sequence and we must cut off the candidate $\hat{\doublevar}$ solution. 

The simplest cut for $\naiveMP$ is a no-good cut which will force it to choose a different position for at least one double vertex or single vertex:
\[ \sum_{\mathclap{\substack{\rankind\in [n-1]:\\ \hat{ \doublevar}_\rankind = 1}}}  (1-\doublevar_\rankind) + \sum_{\mathclap{\substack{\rankind\in [n-1]:\\ \hat{ \doublevar}_\rankind = 0}}}  \doublevar_\rankind \geq 1.\]
However, as this cut eliminates only the  $\hat{\doublevar}$ solution, it is very weak.

Fortunately, this naive decomposition exactly fits into the definition of \emph{combinatorial Benders decomposition} \cite{combinatorialBenders}, since we have one set of constraints in the extensive form $\vertexCP$, namely \eqref{MDCPvertex:pred} that links the $\doublevar$ and $\rankvar$ variables, and there is exactly one $\doublevar$ in each of the constraints. If $\naiveSP$ is infeasible we can solve for an irreducible infeasible subsystem (IIS), let $\iis \subseteq [\K+1,n-1]$,  and pass the following cut to $\naiveMP$ instead of the no-good cut:
{\crev 
\begin{equation} \label{CombBendCutNaive}
 \sum_{\mathclap{\substack{\rankind\in [n-1]:\\ \hat{ \doublevar}_\rankind = 0\\\rankind  \in \iis}}}  \doublevar_\rankind \geq 1
\end{equation}
As we are minimizing the number of doubles, we need to only consider those positions in which the master problem selects non-doubles, namely where $\hat{ \doublevar}_\rankind = 0$.}
\vspace*{-0.5cm}
\begin{algorithm} 
\small
	\SetAlgoLined
	optimal := false \\
	\While{not optimal}{
		solve $\naiveMP$ and get solution $(\hat{\doublevar})$, 
		solve $\naiveSP$ with $(\hat{\doublevar})$ \\
		\eIf{$\naiveSP$ infeasible}{
			find $\iis$, 
			add \eqref{CombBendCutNaive} to $\naiveMP$		
		}{
			let $\hat{ \vertvar}$ be an optimal solution of $\naiveSP$ \\
			accept $(\hat{\doublevar}, \hat{ \vertvar})$, 
			optimal := true 
		}
	}
	\caption{ Naive Decomposition}
	\label{alg:naive}
\end{algorithm}
\vspace*{-0.5cm}
\subsection{Valid Inequalities}
Analysis of the preliminary computational results showed that the doubles in an order were likely to be close to the first clique, which is logical as $\K$ is small compared to the number of vertices and these are the vertices with the least predecessors. The valid inequalities obtained via Algorithm \ref{alg:naiveVI} are the result of analyzing the structure of graphs with various double sequences. 
Given an instance $\instance$ and the set of cliques, $\clique$, in $\graph$, Algorithm \ref{alg:naiveVI} attempts to find structures for which positions $\K+1,\K+2, $ and $\K+3$ in the order can be non-doubles. {\crev We note that these valid inequalities can be added to any other formulation which use rank-based double variables.} We next explain the procedure in Algorithm \ref{alg:naiveVI}, step by step.

\begin{algorithm}[h]
\small
	\SetAlgoLined
	candidates := $\emptyset$ \\ 
	\eIf{$\clique_{\K+2} = \emptyset$}
	{
		$\doublevar_{\K+1} = 1$\\
		found := false\\
		\ForEach{ $\vertexset^C = \vertexset^{C_1} \cup \vertexset^{C_2}$, $C_1, C_2 \in \clique_{\K+1}$ s.t. $|\vertexset^{C_1} \cap \vertexset^{C_2}| = \K$}  	{ \label{foreach1}
			\If{ $\exists \  \vertind \in \vertexset \setminus \vertexset^{C} $ s.t. $|\neighb(\vertind) \cap \vertexset^{C}| \geq \K+1$}
			{ 
				candidates.add($(\vertexset^{C}, v)$) \\
				found := true
			}
		}
		\eIf{not found}{
			$\doublevar_{\K+2} = 1$
		}{
			found := false \\
			\ForEach{ $(\vertexset^{C}, v) \in $ candidates }
			{
				\If{ $\exists \  \vertind' \in \vertexset \setminus (\vertexset^{C}\cup\{v\}) $ s.t. $|\neighb(\vertind') \cap ( \vertexset^{C} \cup\{v\})| \geq \K+1$}
				{ 
					found := true \\
					\textbf{break}
				}
			}
			\If{not found}{
				$\doublevar_{\K+2} + \doublevar_{\K+3} \geq 1$ \label{cut1}
			}
			
		}
	}{ \label{elsestatement}
		found := false\\
		\ForEach{$C \in \clique_{\K+2}$}{
			\If{ $\exists \  \vertind \in \vertexset \setminus \vertexset^{C} $ s.t. $|\neighb(\vertind) \cap \vertexset^{C}| \geq \K+1$}
			{ 
				candidates.add($(\vertexset^{C}, v)$) \\
				found := true
			}
		}
		\eIf{not found}{
			$\doublevar_{\K+1} + \doublevar_{\K+2} \geq 1$ \label{cut2}
		}{
			found := false \\
			\ForEach{ $(\vertexset^{C}, v) \in $ candidates }
			{
				\If{ $\exists \  \vertind' \in \vertexset \setminus (\vertexset^{C}\cup\{v\}) $ s.t. $|\neighb(\vertind') \cap ( \vertexset^{C}\cup\{v\})| \geq \K+1$}
				{ 
					found := true \\
					\textbf{break}
				}
			}
			\If{not found}{
				$\doublevar_{\K+1} + \doublevar_{\K+2} + \doublevar_{\K+3} \geq 1$ \label{cut3}
			}
		}
	}
	\caption{Valid inequalities for $\MinD$}
	\label{alg:naiveVI}
\end{algorithm}

\begin{figure}[h]
	\centering
	\begin{subfigure}[t]{0.3\textwidth}
		\begin{tikzpicture}[scale=0.9,main_node/.style={scale=0.9,circle,fill=white!80,draw,inner sep=0pt, minimum size=16pt},
		line width=1.2pt]
		\node[main_node] (v0) at (0,0) {$v_0$};
		\node[main_node] (v1) at (1.25,0.75) {$v_1$};
		\node[main_node] (v2) at (2.5,0) {$v_2$};
		\node[main_node] (v3) at (1,2.5) {$v_3$};
		\node[main_node] (v4) at (3.5,1) {$v_4$};
		\draw[-, blue] (v0) -- (v1) -- (v2) -- (v3) -- (v2) -- (v0) --(v3) -- (v1);
		\draw[-, ForestGreen] (v3)-- (v4) -- (v1) -- (v4) -- (v2) ;
		\draw[dashed, gray] (v0) edge[out=270, in=270] node [left] {} (v4);
		\end{tikzpicture}
		\caption{Structure for the first $\K+1$ ranks, there is no $(\K+2)$-clique if the dashed edge is not in the input graph.} \label{fig:5clique}
	\end{subfigure}
	~
	\begin{subfigure}[t]{0.3\textwidth}
		\begin{tikzpicture}[scale=0.9,main_node/.style={scale=0.9,circle,fill=white!80,draw,inner sep=0pt, minimum size=16pt},
		line width=1.2pt]
		\node[main_node] (v0) at (0,0) {$v_0$};
		\node[main_node] (v1) at (1.25,0.75) {$v_1$};
		\node[main_node] (v2) at (2.5,0) {$v_2$};
		\node[main_node] (v3) at (1,2.5) {$v_3$};
		\node[main_node] (v4) at (3.5,1) {$v_4$};
		\draw[-, blue] (v0) -- (v1) -- (v2) -- (v3) -- (v2) -- (v0) --(v3) -- (v1);
		\draw[-, ForestGreen] (v3)-- (v4) -- (v1) -- (v4) -- (v2) ;
		\draw[dashed, white] (v0) edge[out=270, in=270] node [left] {} (v4);
		\end{tikzpicture}
		\caption{Two overlapping $(\K+1)$-cliques with intersection of exactly $\K$ vertices. Here $\doublevar_{\K+1}$ is a double.} \label{fig:overlap}
	\end{subfigure}
	~
	\begin{subfigure}[t]{0.3\textwidth}
		\begin{tikzpicture}[scale=0.9,main_node/.style={scale=0.9,circle,fill=white!80,draw,inner sep=0pt, minimum size=16pt},
		line width=1.2pt]
		\node[main_node] (v0) at (0,0) {$v_0$};
		\node[main_node] (v1) at (1.25,0.75) {$v_1$};
		\node[main_node] (v2) at (2.5,0) {$v_2$};
		\node[main_node] (v3) at (1,2.5) {$v_3$};
		\node[main_node] (v4) at (3.5,1) {$v_4$};
		\node[main_node] (v5) at (-0.5,0.75) {$v_5$};
		\draw[-, blue] (v0) -- (v1) -- (v2) -- (v3) -- (v2) -- (v0) --(v3) -- (v1);
		\draw[-, ForestGreen] (v3)-- (v4) -- (v1) -- (v4) -- (v2) ;
		\draw[-, red] (v0) -- (v5) -- (v2) (v5)-- (v1);
		\draw[dashed, white] (v0) edge[out=270, in=270] node [left] {} (v4);
		\end{tikzpicture}
		\caption{Three overlapping $(\K+1)$-cliques with intersection of exactly $\K$ vertices. Here $\doublevar_{\K+1}$ and $\doublevar_{\K+2}$ are doubles.} \label{fig:overlap3}
	\end{subfigure}
	\caption{Graph substructures which lead to valid inequalities for $\naiveMP$, with $\K=3$.}
	\label{fig:naiveVI}
\end{figure}
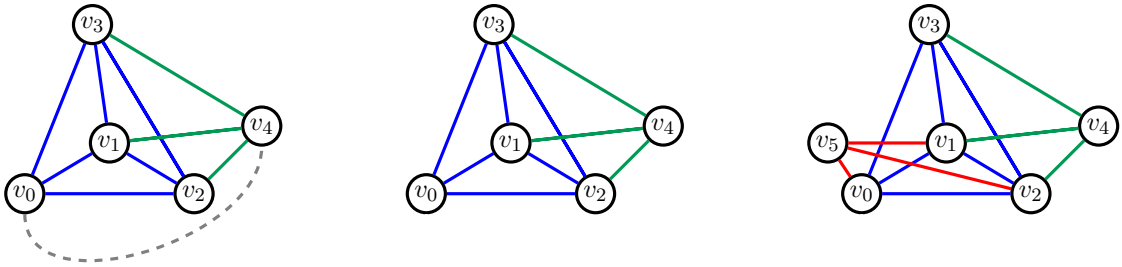

The algorithm begins by checking if there is a clique of size $\K+2$ in $\graph$. If such a clique does not exist, the vertex at position $\K+1$ must be a double since we are unable to extended the initial $(\K+1)$-clique to a  $(\K+2)$-clique. Thus the vertex at position $\K+1$ cannot be adjacent to all  $\K+1$ of its predecessors and must have exactly $\K$ predecessors instead. In this case we obtain $\doublevar_{\K+1} = 1$ as a valid inequality. As seen in Figure \ref{fig:5clique}, if the dashed edge does not  exists, there is no $(\K+2)$-clique, so the initial clique can be $v_0, v_1, v_2, v_3$ and $v_4$ is a double or the initial clique can be $v_1, v_2, v_3, v_4$ and $v_0$ is double, thus the vertex at position $\K+1$ is always a double.

If $\doublevar_{\K+1}$ has been fixed to be a double, we search further for a structure in the graph that would allow position $\K+2$ to be a non-double. On line \ref{foreach1}, we consider all possible graph sub-structures which can fill positions $0$ to $\K+1$ in the order. Since there is no $(\K+2)$-clique, these structures are two $(\K+1)$-cliques, whose intersection is exactly $\K$ vertices,  as seen in Figure \ref{fig:overlap},  one clique is the initial clique and the second contains the vertex at position $\K+1$ and its $\K$ adjacent predecessors. We would like to extend these structures by another vertex which is adjacent to at least $\K+1$ of the vertices in the structure. If such a vertex exists then we have a candidate structure for which position $\K+2$  can be a non-double, if no such candidate is found, we know $\doublevar_{\K+2}$ is a double, this is the case in Figure \ref{fig:overlap3}, where we have three cliques $(\K+1)$-cliques, whose intersection is exactly $\K$ vertices but no $(\K+2)$-clique as a  substructure. Finally if we have at least one candidate we try to extend it in the same manner as before by looking for a vertex which could be a non-double in position $\K+3$, if one exits we cannot say anything about the whether  position $\K+3$ has a double. However, if no such structure exists, we can add the valid inequality $\doublevar_{\K+2} + \doublevar_{\K+3} \geq 1$ (line \ref{cut1}), since we know positions $\K+2$ and $\K+3$ are not both non-double.

If there is a  $(\K+2)$-clique (line \ref{elsestatement}), it is possible for position $\K+1$ to be a non-double, so we similarly search for a vertex which can extend one of the $(\K+2)$-cliques so that the vertex is a non-double. If no such vertex exists, we know positions $\K+1$ and $\K+2$ are not both non-double, so we can add the valid inequality $\doublevar_{\K+1} + \doublevar_{\K+2} \geq 1$ (line \ref{cut2}). Finally, if we have one such vertex clique candidate, we can try to extend it a final time to find a structure which would allow all three positions, $\K+1,\K+2, $ and $\K+3$, to be non-doubles. If we are unable to find such a candidate structure, we add the valid inequality $ \doublevar_{\K+1} +\doublevar_{\K+2} + \doublevar_{\K+3} \geq 1$ (line \ref{cut3}) as we need at least one double vertex in these positions.

We are also able to fix the double value of positions starting from the end of the order. We remark that if there is no vertex with degree exactly $\K$, then the last position in the order cannot be a double because it must have at least $\K+1$ adjacent predecessors. Similarly, if there is no vertex with degree greater than or equal to $\K+1$, positions $n-1$ and $n-2$ cannot be doubles because the vertex with smallest degree has at least $\K+2$ neighbours so in the best case if it is in position $n-1$, it has $\K+2$ adjacent predecessors, and if it is in position $n-2$, it has at least $\K+1$ adjacent predecessors. We extend this notion to a general rule based on the minimum degree of a vertex in the input graph.

\begin{fix}
	Given, $\instance$, let $m$ be the minimum degree of $v \in \vertexset$. Then, we can fix 
	$\doublevar_{n-i} =0 \quad \forall \  i= [1, m-\K].
	$
\end{fix}

Preliminary computational results show that the performance of the naive decomposition is weak. The master problem solves in less than a second,  and using combinatorial Benders, we were able to observe that the problem converges in relatively few iterations, however finding an IIS requires a considerable amount of time.  Thus this naive decomposition approach is unbalanced, it has a very weak master problem and a very strong subproblem. This motivates the next decomposition approach. We remark, however, that the addition of valid inequalities and variable fixing strengthens the naive decomposition significantly.

\subsection{A Witness-based Decomposition and an Extended Formulation}
\label{subsec:witness}
 
Ideally, our decomposition would be more balanced; the master problem would make decisions about the doubles but also about some aspects of the vertex ordering. However, since our current problem space only has two types of decisions, namely the order and the doubles, we are unable to decompose any further. Thus, to help improve the decomposition we add more decisions. First, we allow the problem to reason directly about which vertices will be in the initial clique. We would also like to add some decisions that restrict the space of vertex orders but is less restrictive than the linear ordering constraints used in $\OGcycles$, to do so we introduce the idea of a \textit{witness}. We say a vertex $u$ witnesses vertex $\vertind$ in a given DVOP order if $u$ supports the validity of the position of $\vertind$. More specifically, the witnesses of $\vertind$ are a minimal set of adjacent predecessors of $\vertind$ so that each double vertex has exactly $\K$ witnesses, and each non-double vertex has exactly $\K+1$ witnesses. We also make the convention that all the vertices of the initial $(\K+1)$-clique witness each other and all of their neighbours. The latter means that if a vertex, $\vertind$, outside the initial clique is adjacent to a vertex $u$ inside the initial clique, then $u$ will always be a witness to $\vertind$. 

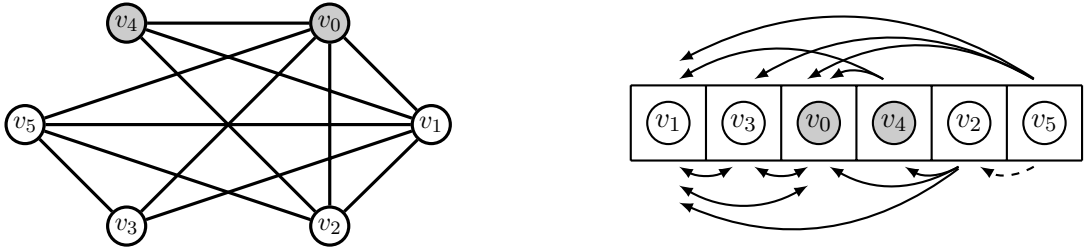
\begin{figure}[h]
	\centering
	\begin{subfigure}[b]{0.6\textwidth}
		\begin{center}
			\begin{tikzpicture}[scale=0.9,main_node/.style={scale=0.9,circle,fill=white!80,draw,inner sep=0pt, minimum size=16pt},
			line width=1.2pt]
			\node[main_node,fill = gray!40] (v0) at (3,3) {$v_0$};
			\node[main_node] (v1) at (4.5,1.5) {$v_1$};
			\node[main_node] (v2) at (3,0) {$v_2$};
			\node[main_node] (v3) at (0,0) {$v_3$};
			\node[main_node,fill = gray!40] (v4) at (0, 3) {$v_4$};
			\node[main_node] (v5) at (-1.5, 1.5) {$v_5$};
			\draw[-] (v5) -- (v0) -- (v1) -- (v2) -- (v0) -- (v4) -- (v1) -- (v3) -- (v0);
			\draw[-] (v4) -- (v2) --(v5);
			\draw[-] (v1) -- (v5) --(v3);
			\end{tikzpicture}
		\end{center}
		\caption{A feasible DVOP order $\Lorder v_1, v_3, v_0, v_4, v_2, v_5\Rorder$, with $\K =2$.}
		\label{fig:DVOPgraph2} 
	\end{subfigure}
	~
	\begin{subfigure}[b]{0.35\textwidth}
		\begin{center}
			\begin{tikzpicture}[main_node/.style={circle,fill=white!80,draw,inner sep=0pt, minimum size=16pt}, thick]
			\draw (0,0) grid (6,1);
			\node[main_node] (v1) at (0.5,0.5) {$v_1$};
			\node[main_node] (v3) at (1.5,0.5) {$v_3$};
			\node[main_node,fill = gray!40] (v0) at (2.5,0.5) {$v_0$};
			\node[main_node,fill = gray!40] (v4) at (3.5,0.5) {$v_4$};
			\node[main_node] (v2) at (4.5,0.5) {$v_2$};
			\node[main_node] (v5) at (5.5,0.5) {$v_5$};
			
			\node[gray] (a1) at (0.5,0) {};
			\node[gray] (b1) at (0.5,-.25) {};
			\node[gray] (c1) at (0.5,1) {};
			\node[gray] (d1) at (0.5,-.5) {};
			\node[gray] (e1) at (0.5,1.25) {};
			\node[gray] (a3) at (1.5,0) {};
			\node[gray] (b3) at (1.5,1) {};
			\node[gray] (a0) at (2.5,0) {};
			\node[gray] (b0) at (2.5,-.25) {};
			\node[gray] (c0) at (2.5,1) {};
			\node[gray] (d0) at (2.2,1) {};
			\node[gray] (a4) at (3.5,1) {};
			\node[gray] (b4) at (3.5,0) {};
			
			\node[gray] (a2) at (4.5,0) {};
			
			\node[gray] (a5) at (5.5,1) {};
			\node[gray] (b5) at (5.5,0) {};
			
			\draw[<->, >=latex] (a1) edge [bend right] (a3) ;
			\draw[<->, >=latex] (a3) edge [bend right] (a0) ;
			\draw[<->, >=latex] (b1) edge [bend right] (b0) ;
			\draw[->, >=latex] (a4) edge [bend right] (c0) ;
			\draw[->, >=latex] (a4) edge [bend right] (c1) ;
			\draw[->, >=latex] (a2) edge [bend left] (b4) ;
			\draw[->, >=latex] (a2) edge [bend left] (a0) ;
			\draw[->, >=latex] (a2) edge [bend left] (d1) ;
			\draw[->, >=latex] (a5) edge [bend right] (d0) ;
			\draw[->, >=latex] (a5) edge [bend right] (e1) ;
			\draw[->, >=latex] (a5) edge [bend right] (b3) ;
			\draw[->, >=latex, dashed] (b5) edge [bend left] (a2) ;
			\end{tikzpicture}
		\end{center}
		\caption{The witness graph.}
		\label{fig:DVOPwitnessGraph} 
	\end{subfigure}
	\caption{An instance with a DVOP order, and its witness graph where there is a solid arc $(v,u)$ if $u$ witnesses $v$, while there is a dashed arc $(v,u)$ if $u$ is adjacent to $v$ in the input graph but is not selected to witness $v$. (The double vertices in the order are highlighted in gray.)}
	\label{fig:DVOPgraphWitness}
\end{figure}

For example, when $\K=2$, the graph in Figure \ref{fig:DVOPgraph2} has DVOP order $\Lorder v_1, v_3, v_0, v_4, v_2, v_5\Rorder$  with two doubles, namely $v_0$ and $v_4$. Since the vertices in the first clique witness each other we have  $v_1$ is witnessed by $v_3$ and $v_0$, $v_3$ is witnessed by $v_1$ and $v_0$, and $v_0$ is witnessed by $v_3$ and $v_1$. Vertex $v_4$ is a double outside the initial clique and thus has $\K = 2$ adjacent predecessors, meaning it needs $\K$ witnesses, $v_4$ is witnessed by $v_0$ and $v_1$. Vertex $v_2$ is not a double, so it must have at least $\K +1$ adjacent predecessors, meaning it needs $\K+1$ witnesses. Since $v_2$ has exactly $\K+1$ adjacent predecessors, it can only be witnessed by them, i.e., $v_0, v_1$ and $v_4$ are all witnesses for $v_2$. Finally, vertex $v_5$ is not a double so it needs  $\K+1$ witnesses. However, $v_5$ has four adjacent predecessors, due to our convention, we choose the initial clique vertices $v_1, v_3,$ and $ v_0$ to witness $v_5$, although any combination of three vertices from  $v_1, v_3, v_0 $ and $ v_2$ would be a valid choice for the witness set of $v_5$. The relationship between vertices and their witnesses can be seen in Figure \ref{fig:DVOPwitnessGraph}.

In our new decomposition, the master problem decides which vertices are doubles and which are in the initial clique, as well as an appropriate number of witnesses for each vertex. The subproblem then looks for an order that respects the clique and witness solutions. For the master problem, in addition to double variables, $\doublevar_\vertind$, we introduce binary variables $\cliquevar_\vertind = 1$ if  vertex $\vertind \in \vertexset$ is in the initial $(\K+1)$-clique, and $\witnessvar_{\vertind u} = 1 \ \forall \  \vertind \in \vertexset, u \in \neighb(\vertind)$ if $v$ uses $u$ as a witness, i.e., if $u$ is a witness for $\vertind$. The  $\cliquevar$ fix the vertices in the initial clique, and the $\witnessvar$ enforce the required number of witnesses for each vertex, depending on $\doublevar$. This formulation incorporates some of the ideas from linear ordering seen in \cite{Omer2017}, {\crev including defining ($2 |\edgeset|$-many) witness variables on the edges of the graph. The key difference between the two formulations is the treatment of the initial clique. Whereas \cite{Omer2017} enumerate the possible first cliques of the order a priori, we constrain the vertices in these positions to form a clique. Further comparisons of the similarities and differences are presented in Appendix \ref{app:ccg_wb}.}  
The new master problem is defined below:
\bsubeq \label{form:EFMP}
\begin{alignat}{2}
\efMP: \min \ &  \sum_{\vertind \in \vertexset} \doublevar_\vertind + 1&& \label{EFMP:obj} \\
\text{s.t.} \ & \sum_{\vertind \in \vertexset} \cliquevar_\vertind = \K + 1 &&  \label{EFMP:1stclique}\\
& \cliquevar_\vertind + \cliquevar_u  \leq 1 && \forall \ \vertind, u \in \vertexset : u \neq \vertind, \{\vertind, u\} \notin \edgeset \label{EFMP:1stcliqueadj}  \\
& \cliquevar_\vertind \leq \witnessvar_{u \vertind}&& \forall \ \vertind \in \vertexset, u \in \neighb(\vertind) \label{EFMP:cliquewitness} \\
& \sum_{u \in \neighb(\vertind)} \witnessvar_{\vertind u} = (\K+1) (1-\cliquevar_\vertind) - \doublevar_\vertind + K\cliquevar_\vertind \quad && \forall  \ \vertind \in \vertexset  \label{EFMP:linking} \\
&  \doublevar \in  \{0,1\}^{|\vertexset|}, \cliquevar \in \{0,1\}^{|\vertexset|}, \witnessvar \in \{0,1\}^{2 |\edgeset|} \label{EFMP:dom} 
\end{alignat}
\esubeq
Constraint \eqref{EFMP:1stclique} ensures that exactly $\K+1$ vertices are selected for the initial $(\K+1)$-clique and constraints \eqref{EFMP:1stcliqueadj} ensure if two vertices are not adjacent, they cannot both be in the initial clique. Constraints \eqref{EFMP:cliquewitness} ensure a vertex will always be witnessed by its neighbours in the initial clique. Constraints \eqref{EFMP:linking} link the double variables, the clique variables and the witness variables. If a vertex $\vertind$ is a double, $\doublevar_{\vertind} = 1$, and outside the initial clique, $\cliquevar_\vertind = 0 $, it must have $\K$ witnesses over all its neighbours, given by the first and second terms of the right-hand-side. If $\vertind$ is not a double, $\doublevar_{\vertind} = 0$, but still outside the initial clique, $\cliquevar_\vertind = 0 $, the number of witnesses it requires is $\K+1$, since all but the first term cancel. If $\vertind$ is in the initial clique, $\cliquevar_\vertind =1$ it must be witnessed by all other vertices in the initial clique, by  \eqref{EFMP:cliquewitness}, which has size $\K+1$ so $\vertind$ has $\K$ witnesses, implying  $\doublevar_{\vertind} = 0$. Finally, objective function \eqref{EFMP:obj} minimizes the number of doubles vertices but also accounts for the linking constraint setting the vertex in position $\K$ to be a non-double when it is always a double by adding  $1$ to the objective.

The subproblem can be any $\rankfunc(\cdot)$ as long as it is properly linked to $\efMP$ variables. We use $\rankCP$ constraints in the subproblem. Given an $\efMP$ candidate solution $(\hat{\cliquevar}, \hat{\witnessvar}, \hat{\doublevar})$, the subproblem looks for a consistent DVOP order:
\bsubeq
\begin{alignat}{2}
\efSP: & \text{AllDifferent}(\rankvar_0,\rankvar_1, ..., \rankvar_{n-1}) \label{cp2}\\ 
& \hat{ \cliquevar}_\vertind =1 \implies \rankvar_\vertind \leq \K &&\forall \ \vertind \in \vertexset \label{cp3}\\
&\hat{ \cliquevar}_\vertind + \hat{ \cliquevar}_u  = 0 \text{ and } \hat{\witnessvar}_{ \vertind u}= 1 \implies \rankvar_u \leq  \rankvar_\vertind -1 \quad && \forall \ \vertind \in \vertexset, u \in \neighb(\vertind) \label{cp4}\\
& \rankvar \in [n-1]^{|\vertexset|}  \label{cp1}
\end{alignat}
\esubeq
Constraints \eqref{cp2} and \eqref{cp1} ensure there is a linear order. Constraints  \eqref{cp3} force all those vertices which have been chosen to be in the initial clique to take a rank in $ [\K]$, and since there are $\K+1$ vertices selected in the clique all these positions will be filled. Note that the actual order of the clique does not matter, as we are merely looking for the existence of an order and any permutation of the initial clique will give a valid DVOP order.  Constraints \eqref{cp4} ensure that for all vertices outside the initial clique, if a vertex $u$ witnesses vertex $\vertind$, then $u$ cannot come after $\vertind$ in the order as otherwise it would not be a valid witness. 

Given a feasible master problem solution $(\hat{\cliquevar}, \hat{\witnessvar}, \hat{\doublevar})$, if $\efSP$ is feasible, we say that there exists a subproblem feasible solution, or a $\rankfunc(\cdot)$, \emph{respecting $(\hat{\cliquevar}, \hat{\witnessvar})$}. Similarly, if there exists a subproblem solution satisfying \eqref{cp2}, \eqref{cp3}, and \eqref{cp1}, but not necessarily \eqref{cp4}, we say that there exists a subproblem solution, or a $\rankfunc(\cdot)$, \emph{respecting $\hat{\cliquevar}$}.

If the master problem solution does not lead to a subproblem feasible solution respecting it, thus a feasible DVOP, we would like to cut it off. To develop such cuts, we consider the possible ``sources of infeasibility" in a master problem solution $(\hat{\cliquevar}, \hat{\witnessvar}, \hat{\doublevar})$. Clearly, since $\hat{\doublevar}$ does not appear in our subproblem it is not a source of infeasibility. The $ \hat{\witnessvar}$, assuming $\hat{\cliquevar}_v=0$ for the vertices in question may lead to an infeasible solution due to the logical constraints \eqref{cp4}, which ensure the proper ranks of the vertices, contradicting with the constraints \eqref{cp2} and  \eqref{cp1} which ensure a linear order. As we will show later, this will indeed arise when the given master problem solution has ``cycles on the $\witnessvar$". For example, recall the graph in Figure \ref{fig:DVOPgraph2} with order $\Lorder v_1, v_3, v_0, v_4, v_2, v_5\Rorder$, a solution which assigns $\hat{\cliquevar}_{v_2}= \hat{\cliquevar}_{v_5} = 0$ and $\hat{\witnessvar}_{ v_2 v_5} =  \hat{\witnessvar}_{ v_5 v_2} =1$ is feasible to the master problem. However, by \eqref{cp2} and \eqref{cp1}, the constraints \eqref{cp4} give a contradiction as we would like to have $v_2$ precedes $v_5$ in the order and $v_5$ precedes $v_2$ simultaneously. This can be interpreted as having a 2-cycle, consisting of nodes $v_2$ and $v_5$, in the $\hat{\witnessvar}$ solution.

As intuitively explained above, we can think of infeasibilities caused by $\hat{\witnessvar}$ as \textit{cycles}. Formally, let $\digraph = (\vertexset, \arcset)$ be the directed graph equivalent of $\graph = (\vertexset, \edgeset)$, where $\arcset = \{(i,j),  (j,i) \ : \  \forall \ \{i, j\} \in \edgeset\}$ where $(i,j)$ denotes a directed arc from $i$ to $j$. If vertex $u$ witnesses vertex $\vertind$, i.e., $\hat{\witnessvar}_{ \vertind u} =1$, this is equivalent to  the arc $(v,u)$ being selected in $\digraph$, as in the solid arcs in the witness graph in Figure \ref{fig:DVOPwitnessGraph}. So, if  vertex $u$ witnesses vertex $\vertind$, and $\vertind$ witnesses  $u$, i.e., $\hat{\witnessvar}_{ \vertind u}= \hat{\witnessvar}_{u \vertind} =1$, we have selected $(v,u)$ and $(u,v)$, which forms a 2-cycle in $\digraph$.  We define \emph{cycles in  $\witnessvar$} as $C = (\vertexset^C, \arcset^C)$, where $C$ is a cycle subgraph of $\digraph$ and let $\cycles$ denote the set of all cycles in $\digraph$.

The only other potential source of subproblem infeasibility in a master problem solution $(\hat{\cliquevar}, \hat{\witnessvar}, \hat{\doublevar})$ is the $\hat{\cliquevar}$ vector. However, we will next prove that the cycles in $\witnessvar$ are indeed the only source of infeasibility. In that regard, the first observation is that $\hat{\cliquevar}$ provides a feasible initial $(\K+1)$-clique for a DVOP order that the subproblem is searching for, which is stated in the following lemma. This is illustrated by the solid arcs in Figure \ref{fig:WitnessProof}.

\begin{lemma} \label{lemma:feasMP}
	If $(\hat{\cliquevar}, \hat{\witnessvar}, \hat{\doublevar})$ is feasible to $\efMP$ then there exists $\hat{\rankvar} \in \Z^{|\vertexset|}$ such that \eqref{cp2}, \eqref{cp3}, and \eqref{cp1} are satisfied at $(\hat{\cliquevar}, \hat{\rankvar})$. 
\end{lemma}

\proof{Proof.}  
	Let $\hat{\vertexset}^1 = \{ \vertind \in \vertexset : \hat{\cliquevar}_\vertind = 1\}$ and $\hat{\vertexset}^0 = \{ \vertind \in \vertexset : \hat{\cliquevar}_\vertind = 0\}$. Note that $|\hat{\vertexset}^1| = \K+1$ and $|\hat{\vertexset}^0| = |\vertexset| - (\K+1)$ due to \eqref{EFMP:1stclique} and \eqref{EFMP:dom}. Arbitrarily ordering $\hat{\vertexset}^1$ and $\hat{\vertexset}^0$ as $\hat{\vertind}^1_{1},\hat{\vertind}^1_{2},\hdots,\hat{\vertind}^1_{\K+1}$ and $\hat{\vertind}^2_{1},\hat{\vertind}^2_{2},\hdots,\hat{\vertind}^2_{|\vertexset| - (\K+1)}$, respectively, construct  $\hat{\rankvar}$ as 
	\[ \hat{\rankvar}_\vertind = \begin{cases} 
	i-1, & \text{if} \ \vertind = \hat{\vertind}^1_{i} \ \text{for some} \ i \in [1,\K+1] \\
	i+\K, & \text{if} \ \vertind = \hat{\vertind}^2_{i} \ \text{for some} \ i \in [1,|\vertexset| - (\K+1)] 
	\end{cases} \quad \forall \ \vertind \in \vertexset.
	\]
	Then, it is easy to confirm that \eqref{cp2}, \eqref{cp3}, and \eqref{cp1} hold at $(\hat{\cliquevar}, \hat{\rankvar})$. 
\endproof

The next result shows that in a master {\crev feasible solution} $\hat{\witnessvar}$, every vertex in the initial clique is witnessed only by vertices in the initial clique, i.e., there cannot be the dotted back-edges illustrated in Figure \ref{fig:WitnessProof}.
\vspace*{-0.3cm}
\begin{figure}[h]
\centering
\vspace*{-0.3cm}
		\begin{tikzpicture}[main_node/.style={circle,fill=white!80,draw,inner sep=0pt, minimum size=16pt}, box/.style={rectangle,draw=black,thick, minimum size=1cm}, thick]
		\draw[step=1cm] (0,0) grid (6,1);
		\node[box,fill=gray!40] at (0.5,0.5) {};
		\node[box,fill=gray!40] at (1.5,0.5) {};
		\node[box,fill=gray!40] at (2.5,0.5) {};
		\node[] (v1) at (0.5,0.5) {\footnotesize $\hat{\cliquevar}_\vertind = 1$};
		\node[] (v3) at (1.5,0.5) {\footnotesize $\hat{\cliquevar}_\vertind = 1$};
		\node[] (v0) at (2.5,0.5) {\footnotesize $\hat{\cliquevar}_\vertind = 1$};
		\node[] (v4) at (3.5,0.5) {\footnotesize $\hat{\cliquevar}_\vertind = 0$};
		\node[] (v2) at (4.5,0.5) {\footnotesize $\hat{\cliquevar}_\vertind = 0$};
		\node[] (v5) at (5.5,0.5) {\footnotesize $\hat{\cliquevar}_\vertind = 0$};
		
		\node[gray] (a1) at (0.5,0) {};
		\node[gray] (b1) at (0.5,-.25) {};
		\node[gray] (c1) at (0.5,1) {};
		\node[gray] (d1) at (0.5,-.5) {};
		\node[gray] (e1) at (0.5,1.25) {};
		\node[gray] (a3) at (1.5,0) {};
		\node[gray] (b3) at (1.5,1) {};
		\node[gray] (a0) at (2.5,0) {};
		\node[gray] (b0) at (2.5,-.25) {};
		\node[gray] (c0) at (2.5,1) {};
		\node[gray] (d0) at (2.2,1) {};
		\node[gray] (a4) at (3.5,1) {};
		\node[gray] (b4) at (3.5,0) {};
		\node[gray] (a2) at (4.5,0) {};
		\node[gray] (a5) at (5.5,1) {};
		\node[gray] (b5) at (5.5,0) {};
		
		\draw[<->, >=latex] (a1) edge [bend right] (a3) ;
		\draw[<->, >=latex] (a3) edge [bend right] (a0) ;
		\draw[<->, >=latex] (b1) edge [bend right] (b0) ;
		\node[] (lemma1) at (3.3,-0.4) {\small Lemma \ref{lemma:feasMP}};
		
		\node[gray] (dot1) at (1.5,1) {};
		\node[gray] (dot2) at (4.5,1) {};
		\draw[->, >=latex, dotted] (dot1) edge [bend left] (dot2);
		\node (mental) at (3,1.5) {\large $\times$};
		
		\node[gray] (dot3) at (2.2,1) {};
		\node[gray] (dot4) at (3.5,1) {};
		\draw[->, >=latex, dotted] (dot3) edge [bend left] (dot4);
		\node (mental) at (2.8,1.2) {\large $\times$};
		
		\node[] (lemma2) at (4.5,1.6) {\small Lemma \ref{lemma:feasMPclique}};
		
		\end{tikzpicture}			
		\caption{Witness graph for a given master candidate solution $(\hat{\cliquevar}, \hat{\witnessvar}, \hat{\doublevar})$ for a case with $\K =2$. The initial $(\K+1)$-clique positions are in gray. The arcs correspond to the witness relationships provided by $\hat{\witnessvar}$.}
		\label{fig:WitnessProof} 
\end{figure}
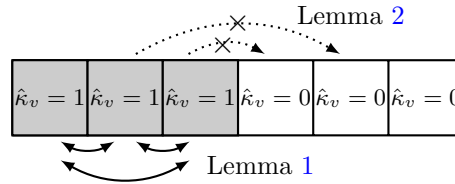
\begin{lemma} \label{lemma:feasMPclique}
	If $(\hat{\cliquevar}, \hat{\witnessvar}, \hat{\doublevar})$ is feasible to $\efMP$, then $\hat{\witnessvar}_{\vertind u} = 0 \ \forall \ \{ \vertind,u\} \in \edgeset \ \text{with} \ \hat{\cliquevar}_\vertind = 1, \hat{\cliquevar}_u =0$. 
\end{lemma}

\proof{Proof.} 
Given  $(\hat{\cliquevar}, \hat{\witnessvar}, \hat{\doublevar})$ feasible to $\efMP$, let $v \in \vertexset$ with $\hat{\cliquevar}_{v} =1$, that is $v$ has been selected for the initial clique. By constraints \eqref{EFMP:1stclique}, \eqref{EFMP:1stcliqueadj}, and \eqref{EFMP:cliquewitness} and the $\witnessvar$ domains, $\hat{\witnessvar}_{\vertind i} = 1$ for all $i \in \vertexset$, such that $v \neq i$ and $\hat{\cliquevar}_i =1$, i.e., $v$ is witnessed at least by $\K$ other vertices in the first clique. By \eqref{EFMP:linking}, we have $\sum_{j \in \neighb(\vertind)} \hat{\witnessvar}_{\vertind j} = \K - \hat{\doublevar}_\vertind $ but since $\hat{\doublevar} \geq 0$, $\sum_{j \in \neighb(\vertind)} \hat{\witnessvar}_{\vertind j} \leq \K$.
		Thus, $v$ is witnessed by exactly $\K$ vertices, and those are exactly the other vertices in the initial clique. As such,  $\hat{\witnessvar}_{\vertind u} = 0 \ \forall \ \{ \vertind,u\} \in \edgeset \ \text{with} \ \hat{\cliquevar}_\vertind = 1, \hat{\cliquevar}_u =0$ as required.
\endproof

Lemma \ref{lemma:feasMPclique} implies that there is no cycle in $\hat{\witnessvar}$ including vertices both from inside and outside of the first $(\K+1)$-clique described by $\hat{\cliquevar}$ since completing such a cycle would require the use of a back-edge from the initial clique to outside. 

\begin{corollary}
	\label{cor:cycle}
	Let $(\hat{\cliquevar}, \hat{\witnessvar}, \hat{\doublevar})$ be feasible to $\efMP$. Define $\hat{\vertexset}^1 = \{ \vertind \in \vertexset : \hat{\cliquevar}_\vertind = 1\}$ and $\hat{\vertexset}^0 = \{ \vertind \in \vertexset : \hat{\cliquevar}_\vertind = 0\}$. For any $C = (\vertexset^C, \arcset^C) \in \cycles$ with $\hat{\witnessvar}_{\vertind u} = 1$ for all $(\vertind , u) \in \arcset^C$, either $\vertexset^C \cap \hat{\vertexset}^1 = \emptyset$ or $\vertexset^C \cap \hat{\vertexset}^0 = \emptyset$.
\end{corollary}

Given an $\efMP$ solution $(\hat{\cliquevar}, \hat{\witnessvar}, \hat{\doublevar})$, as there always exists a feasible initial $(\K+1)$-clique by Lemma \ref{lemma:feasMP}, if the subproblem is infeasible, the issue is failing to build the rest of the order. As mentioned before, a reason could be the existence of a cycle in the $\hat{\witnessvar}$ solution. Due to Corollary \ref{cor:cycle}, such a cycle can only occur outside of the initial clique. We propose the following  inequality to remove infeasibilities caused by a cycle $C = (\vertexset^C, \arcset^C)$  in the $\witnessvar$,
\begin{equation} \label{cycleBreakingCut}
\sum_{(\vertind , u) \in \arcset^C} \witnessvar_{\vertind u} \leq |\vertexset^C|- 1 + \indicate(|\vertexset^C| \leq \K+1) \cliquevar_\iota
\end{equation}
where $\iota$ is any (e.g., the smallest) index of a vertex in $\vertexset^C$. When there is a cycle $C$ with $|\vertexset^C| > \K+1$ we have the classical cycle breaking cut. However, if  $|\vertexset^C| \leq \K+1$ it is possible we have detected a cycle in the initial $(\K+1)$-clique, which is valid. In this case, we need to \emph{lift} the classical cycle breaking into the $(\cliquevar,\witnessvar)$-space to be valid for the master problem, which is accomplished by adding $\cliquevar_\iota$ to the right-hand-side of the cut, so that if the cycle is in the initial clique, $\cliquevar_\iota =1$ and the inequality is redundant. By Corollary \ref{cor:cycle}, $\cliquevar_\iota$ is only $1$ if the cycle is in the initial clique, and otherwise the cycle breaking inequality will hold. The following lemma formalizes the validity of the proposed cuts.

%
%

\begin{lemma} 
	The inequality \eqref{cycleBreakingCut} is valid for $\{({\cliquevar}, {\witnessvar}) : \exists \ \doublevar, \rankvar \text{ s.t. } ({\cliquevar}, {\witnessvar}, \doublevar)$  is feasible to  $\efMP$  and $\rankvar$  respects $({\cliquevar}, {\witnessvar})  \}$.
\end{lemma}

\proof{Proof.} 
 \textit{Case 1:}  If $|\vertexset^C| > \K+1$, $\hat{\witnessvar}_{ij} = 1$ for all $(i,j)  \in \arcset^C$. By Lemma \ref{lemma:feasMPclique} we know $\hat{\cliquevar}_\vertind = 0 $ for all $v \in \vertexset^C$. So we are outside the initial clique and would like to break cycles of any length, and the cut reduces to the standard cycle breaking inequality. \\
\textit{Case 2:}  If $|\vertexset^C| \leq \K+1$ and $\cliquevar_\iota = 1$ the inequality is redundant and thus valid since we only want to break cycles outside the initial clique. If $|\vertexset^C| \leq \K+1$ and $\cliquevar_\iota = 0$ then the cycle is not in the  initial clique and we would like to break it, and again the inequality is the standard cycle breaking inequality.
\endproof

In order to show $\efMP$ and the inequalities \eqref{cycleBreakingCut} for all $C\in \cycles$ are sufficient to prove a solution is optimal to $\MinD$, we introduce the extended formulation obtained by combining $\efMP$ and $\efSP$.  
\bsubeq \label{form:EF}
\begin{alignat}{2}
\hspace*{-0.25cm} \efExtensive: \min \ &  \sum_{\vertind \in \vertexset} \doublevar_\vertind + 1&& \label{EF:obj} \\
\text{s.t.} \ & \sum_{\vertind \in \vertexset} \cliquevar_\vertind = \K + 1 && \label{EF:sizeClique}\\
& \cliquevar_\vertind + \cliquevar_u  \leq 1 && \hspace*{-0.55cm} \forall \ \vertind, u \in \vertexset : u \neq \vertind, \{\vertind, u\} \notin \edgeset\label{EF:adjClique}\\
& \cliquevar_\vertind \leq \witnessvar_{u \vertind}&& \hspace*{-0.55cm} \forall \ \vertind \in \vertexset, u \in \neighb(\vertind)\label{EF:cliqueWitness} \\
& \sum_{u \in \neighb(\vertind)} \witnessvar_{\vertind u} = (\K+1)(1-\cliquevar_\vertind) - \doublevar_\vertind + \K\cliquevar_\vertind \quad &&  \hspace*{-0.55cm}\forall  \ \vertind \in \vertexset \label{EF:linking}\\
& \text{AllDifferent}(\rankvar_0,\rankvar_1, ..., \rankvar_{n-1})\label{EF:allDiff}\\ 
&  \cliquevar_\vertind =1 \Rightarrow \rankvar_\vertind \leq \K && \hspace*{-0.55cm} \forall \ \vertind \in \vertexset\label{EF:rankClique}\\
& \cliquevar_\vertind + \cliquevar_u  = 0 \text{ and } \witnessvar_{ \vertind u } = 1 \Rightarrow \rankvar_u \leq  \rankvar_\vertind -1 \quad &&\hspace*{-0.55cm}  \forall \ \vertind \in \vertexset, u \in \neighb(\vertind) \label{EF:rankPred}\\
&  \doublevar \in  \{0,1\}^{|\vertexset|}, \witnessvar \in \{0,1\}^{|\arcset|}, \cliquevar \in  \{0,1\}^{|\vertexset|}, \rankvar \in  [n-1]^{|\vertexset|}
\end{alignat}
\esubeq
%
Let $\EFfeas$ denote the feasible region of  $\efExtensive$. Finally, we will show that $\efMP$ and the cycle breaking cuts \eqref{cycleBreakingCut} for all $C\in \mathcal{C}$ are sufficient to prove optimality.

\begin{thm} \label{thm:proj}
	If $(\hat{\cliquevar}, \hat{\witnessvar}, \hat{\doublevar})$ is feasible to $\efMP$ and \eqref{cycleBreakingCut} for all $C\in \mathcal{C}$, then $(\hat{\cliquevar}, \hat{\witnessvar}, \hat{\doublevar}) \in proj_{\cliquevar, \witnessvar, \doublevar}(\EFfeas)$, i.e., there exists $\hat{\rankvar}$ such that $(\hat{\cliquevar}, \hat{\witnessvar}, \hat{\doublevar}, \hat{\rankvar}) \in \EFfeas$.
\end{thm}

\proof{Proof.} 
By induction on $p \in \N$ with $\K \leq p \leq n-1$ denoting the position for which there exists a DVOP order respecting $(\hat{\cliquevar}, \hat{\witnessvar})$ in the $[p]$ subset of positions.
	
	\noindent	\textbf{Base Case:} When $p = \K$ there is a DVOP order respecting the $(\hat{\cliquevar}, \hat{\witnessvar})$ since we have a $(\K+1)$-clique from the $\hat{\cliquevar}$ which all witness each other.
	
	\noindent	\textbf{Hypothesis:} Assume we have a DVOP order respecting $(\hat{\cliquevar}, \hat{\witnessvar})$  in the $[p]$ subset of positions.
	
	\noindent	\textbf{Inductive Step:} We show there exists a DVOP order in positions $[p+1]$ respecting $(\hat{\cliquevar}, \hat{\witnessvar})$. Fix a DVOP order respecting  $(\hat{\cliquevar}, \hat{\witnessvar})$ in positions $[p]$. Denote the vertices in the order as $v_1, v_2, \hdots, v_p$, denote, without loss of generality, all other vertices as $q_1, q_2, \hdots, q_{n-p}$.  \\
\underline{\textit{Case 1:}} If there is a $q$ such that $\sum_{i= 1}^{p} \hat{\witnessvar}_{qv_i} = K$, we have an order with $p+1$ where $q$ is the $(p+1)^{\text{th}}$ vertex in the order. So we assume no such vertex exists.
\\
\underline{\textit{Case 2:}} Since all $q$ need at least $\K$ witnesses and we force each vertex to be witnessed by all its neighbours in the initial clique, there is no $q$ with all its witnesses in the initial clique. Thus every $q$ has at least one witness in the set of $q$ vertices or in $v_{\K+1},v_{\K+2}, \hdots, v_p$, however not all $q$ can have such witnesses in the $q$ vertices or there will be a cycle. So there exists a $\hat{\witnessvar}_{q_jv_i} = 1$, that is a $q_j$ that is witnessed by some $v_i \in \{v_{\K+2}, v_{\K+3}, \hdots, v_p\}$.	Without loss of generality, let such a vertex be $q_1$, so $\hat{\witnessvar}_{q_1v_i} = 1$, but it must have $\hat{\witnessvar}_{q_1q_j} = 1$ otherwise it would fall into Case 1. Without loss of generality, let $q_2$ witness $q_1$, so $\hat{\witnessvar}_{q_2q_1} = 0$. If $q_2$ has all witnesses in $v$, we have a DVOP where $q_2$ is the $(p+1)^{\text{th}}$ vertex in the order. So we assume $q_2$ has a witness in $q$, say $q_3$, if $q_3$ has all witnesses in $v$, we have a DVOP where $q_3$ is the $(p+1)^{\text{th}}$ vertex in the order. So we assume $q_3$ has a witness in $q$, say $q_4$, and so on. Otherwise, any $\hat{\witnessvar}_{q_iq_j} = 0 \ \forall \ j < i$, in the worst case $q_i = q_{n-p}$ so it must have all witnesses in $v$. Thus we have a vertex by which we can extend the DVOP order. 
\endproof

We provide a valid inequality for $\efExtensive$ the proof of which is provided in the Appendix \ref{app:proof}.
\begin{prop} 
\label{eq:witnessVI}
	For any  $\vertind \in \vertexset$, the inequality $\doublevar_\vertind \leq 1 - \cliquevar_\vertind$ is valid for $\efExtensive$.
\end{prop}


As the inequalities \eqref{cycleBreakingCut} for all $C\in \mathcal{C}$ are sufficient for proving optimality, we would like to generate them in a cutting plane fashion using $\efSP$ instead of adding all the cycle breaking inequalities to $\efMP$ initially. 
The procedure used to solve the witnessed-based hybrid decomposition is shown in a cutting plane fashion in Algorithm \ref{alg:hybrid} for ease of presentation, however in practice it is implemented as a branch-and-cut algorithm. We note that if the subproblem is infeasible, we generate a cycle $C$ in $\hat{\witnessvar}$, but do not specify the method. As previous stated, the cycle generation can be via a combinatorial Benders procedure wherein we solve an IIS to find the cut. However, we may also use the methods  in \cite{Tarjan1976} to search for a topological ordering of $\hat{\witnessvar}$ using depth-first search. If there is no topological order, we can find a cycle in the digraph since there will be some pair of vertices $u,v$ in the order with $\hat{ \witnessvar}_{vu} =1 $ and and such that the search returns $\rankvar_v \leq \rankvar_u$ which is a contradiction. In this case, $(v, u)$ is in the cycle, with the remaining arcs found on the shortest path from $u$ to $v$.
\begin{algorithm}
\small
	\SetAlgoLined
	optimal := false \\
	\While{not optimal}{
		solve $\efMP$ and get solution $(\hat{\cliquevar}, \hat{\witnessvar}, \hat{\doublevar})$, 
		solve $\efSP$ with $ (\hat{\cliquevar}, \hat{\witnessvar}, \hat{\doublevar})$ \\
		\eIf{$\efSP$ infeasible}{
			generate a cycle $C$ in $\hat{\witnessvar}$, 
			add \eqref{cycleBreakingCut} to $\efMP$		
		}{
			let $\hat{\rankvar}$ be an optimal solution of $\efSP$ \\
			accept  $(\hat{\cliquevar}, \hat{\witnessvar}, \hat{\doublevar}, \hat{\rankvar})$,
			optimal := true 
		}
	}	\caption{Witnessed-based hybrid decomposition. } 
	\label{alg:hybrid}
\end{algorithm}

\section{Computational Results}  \label{sec:results}

In this section we present a computational study which compares the mathematical models presented in Sections \ref{sec:models} and \ref{sec:hybrid}, as well as the existing models from the literature. We first review the data used in the experiments, then the experimental set up, and finally we outline and discuss the results. Full results tables, {\crev with the exception of the pseudo-protein instances}, can be found in the Online Supplement \ref{app:MinD_results}; in this section, we only summarize our main observations from the tables. 

\subsection{Instances}
Our data set consists of {\crev three} kinds of graph instances which are used to test the performance of the mathematical models; random instances, synthetic instances\footnote{Random and synthetic instances can be found at \url{https://sites.google.com/site/mervebodr/home/MINDOUBLE_Instances.zip}.}, and {\crev pseudo-protein instances\footnote{Pseudo-protein instances can be found at \url{https://gitlab.insa-rennes.fr/Jeremy.Omer/MinDouble}.}.}
\paragraph{Random Instances.}
The test data set consisting of randomly generated graphs, with  $n \in \{20,25,30,35\}$ and expected edge density (measured as $\density = \frac{2|\edgeset|}{ n (n-1)} $) in $\{0.3,0.4,0.5\}$. 
The number of doubles for these instances ranged from $1$ to $14$, with an average of $2$ doubles per instance. A summary of the random instances can be seen in Table \ref{table:randInstSize}. Of these instances, only three were infeasible for DVOP with $\K=3$, we left these instances in the test set to ensure our algorithms can also detect infeasible DVOP instances {\crev and observed that for some algorithms, this time is as fast as the greedy algorithm for DVOP.} However, we do not include {\crev these instances} in the performance profiles. 

\begin{table}[h]
\begin{minipage}[t]{0.65\textwidth}
\small
	\caption{Random instances.}
	\label{table:randInstSize}
	\centering
	\begin{tabular}{c c c c c c c} 
		\toprule
		& & & & & Mean  \# \\
		$n$ & $D$ & \# Instances & \# Feasible & \# Infeasible &  Doubles \\
		\midrule
		20& 0.3 & 3 & 1 & 2 & 14 \\
		& 0.4 & 3 & 3 & 0 & 3 \\
		& 0.5 & 3 & 3 & 0 & 1 \\[0.5ex]
		25 & 0.3 & 3 & 2 & 1 & 6 \\
		& 0.4 & 3 & 3 & 0 & 1 \\
		& 0.5 & 3 & 3 & 0 & 1 \\[0.5ex]
		30& 0.3 & 3 & 3 & 0 & 4\\
		& 0.4 & 3 & 3 & 0 & 1 \\
		& 0.5 & 3 & 3 & 0 & 1 \\[0.5ex]
		35& 0.3 & 3 & 3 & 0 &3 \\
		& 0.4 & 3 & 3 & 0 & 1 \\
		& 0.5 & 3 & 3 & 0 & 1 \\[0.5ex]
		\midrule
		Total & & 36 & 33 & 3 & \\
		\bottomrule
	\end{tabular}
\end{minipage}
\begin{minipage}[t]{0.35\textwidth}
\small
	\caption{Synthetic instances.}
	\label{table:SynSize}
	\centering
	\begin{tabular}{c c c } 
		\toprule
		& &  Mean  \# \\
		$n$ & \# Instances & Doubles \\
		\midrule
		25 & 9& 2 \\[0.5ex]
		30& 9 & 3\\[0.5ex]
		35&9 &  4 \\[0.5ex]
		\midrule
		Total & 27 & \\
		\bottomrule
	\end{tabular}
\end{minipage}
\end{table}
	
\paragraph{Synthetic Instances.} Synthetic instances for $\MinD$ were created to simulate the \emph{sparsest} possible graph instances that would still have a DVOP order. Given $\K$, a number of doubles, and a \emph{noise} factor, a synthetic instance is constructed by first building a $(\K+1)$-clique. We then randomly select vertices outside the initial clique to be doubles, next we add edges between vertices ensuring that the doubles and non-doubles have the appropriate number of predecessors in the order. Finally, we add in some noise edges. Notice that the addition of the extra edges may cause there to be less double vertices, however there will always be at least one double vertex by definition.

We generate synthetic instances with $\K=3$ and fix $n \in  \{25, 30, 35 \}$, the number of doubles in  $\in \{\lceil 0.1\times n\rceil, \lceil 0.1\times n\rceil+1, \lceil 0.1\times n \rceil+2 \}$, and a noise factor in $\{0.1,0.15,0.2\}$. For each $n$, doubles, noise triple we generate a minimal graph as described above, where the number of extra edges added is $\lceil n\times \text{noise} \rceil$. {\crev The full algorithm for generating synthetic instances can be found in the Online Supplement \ref{appendix:syn}.} Table \ref{table:SynSize} summarizes the synthetic instances. We remark that the synthetic instances are always feasible.

{\crev \paragraph{Pseudo-Protein  Instances.} The pseudo-protein instances come from \cite{Omer2017}, although they are not the same instances used in their publication. The instances are generated by modifying existing Protein Data Bank instances {\cRev by first trimming the instance to a desired number of nodes, and then randomly removing edges and checking for feasibility until the desired density is met.} In total this test set has $399$ instances with nodes in $\{30,40,\hdots, 100\}$ and expected edge density in $[0.03,0.12]$. We remark that these pseudo-protein instances are always feasible.}

\subsection{Implementation Details}
IPs are solved using IBM ILOG CPLEX version 12.8.0 and CPs are solved using IBM ILOG CP Optimizer version 12.8.0, which we  implemented in C++. All experiments were run on MacOS with 16GB RAM and a 2.3 GHz Intel Core i5 processor, using a single thread. We used the lazy constraint callback function to implement the branch-and-cut procedures. We used $\K=3$ for all experiments and set the time limit to $1000$ seconds. {\crev (We also tested the random instances with $\K = 4$ and  $\K = 5$. As the findings were similar to the $\K=3$ case, their results are provided in Tables \ref{table:K4} and \ref{table:K5} of Online Supplement \ref{app:MinD_results}, respectively.)
	
The three existing formulations of  \cite{Omer2017}, $\OGcycles$, $\OGranks$, and $\CCG$, are re-implemented using the details provided in their publication.  Also, we provided the greedy solutions of \cite{Omer2017} to warm start all the algorithms for random and pseudo-protein instances.}

\subsection{Computational Results and Discussion}
\begin{figure}[h]
	\centering
	\includegraphics[width=0.75\textwidth]{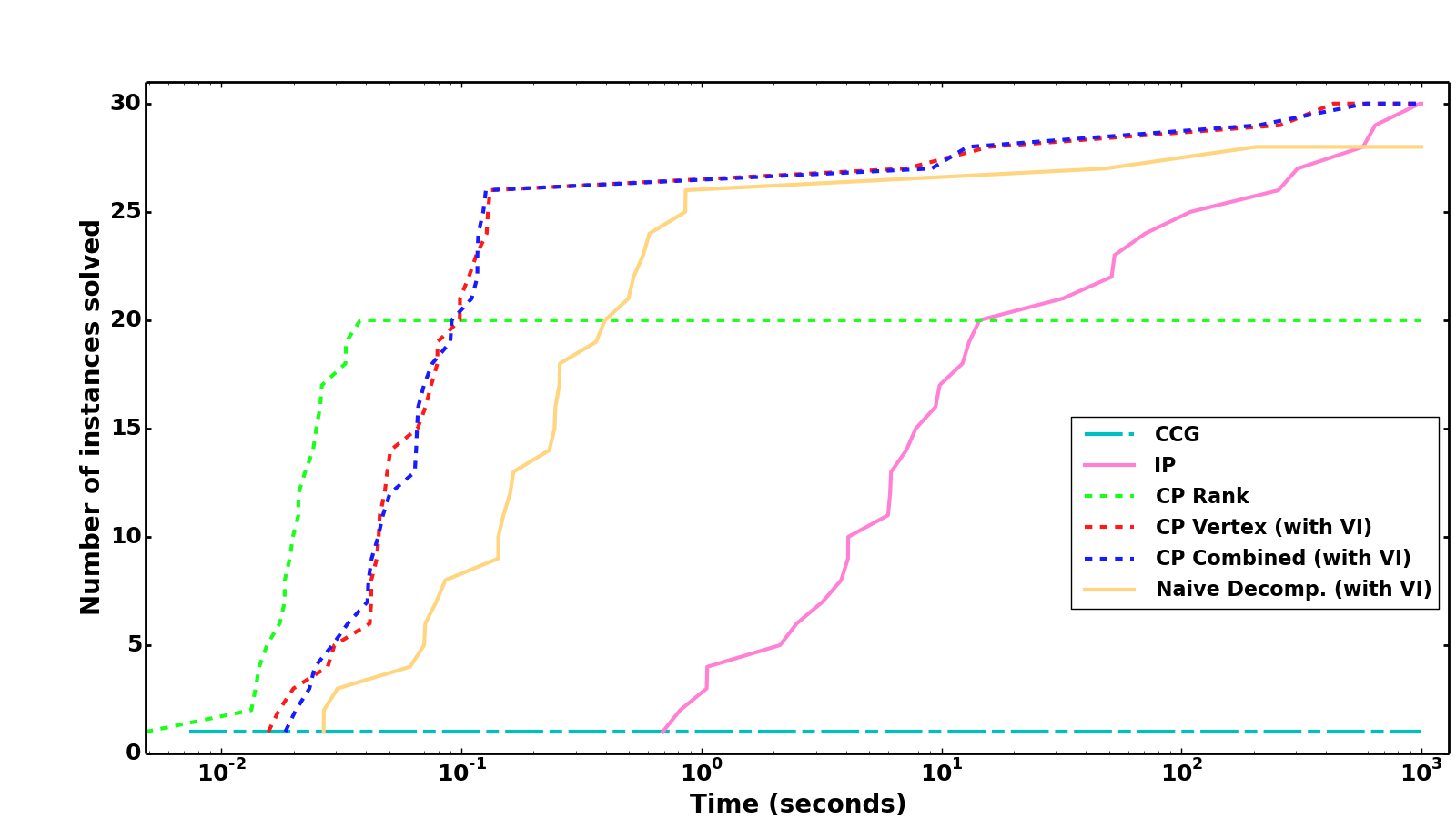}
	\caption{\cRev Solution times of the models on random instances.}
	\label{fig:minDrandom}
\end{figure}

We begin by comparing the performance of selected formulations on the random instances. The performance profile in Figure \ref{fig:minDrandom} shows that many of the methods have similar performance. Overall, $\vertexCP$, $\combCP$, and the naive decomposition {\crev all} with valid inequalities perform the best, although they are only able to solve $30$ of the $33$ feasible instances within the time limit. We remark that $\CCG$ is only able to solve one of the {\cRev feasible} random instances, and all the formulations of this work outperform those the literature by both time and number of instances solved.

The three existing formulations from the literature  $\OGcycles$, $\OGranks$, and $\CCG$ are all outperformed by $\IP$ (details can be seen in Table \ref{table:MinDIP} of Online Supplement \ref{app:MinD_results}).
$\IP$ is able to solve more than half of the instances with $D = 0.4$ and $D = 0.5$ in less than a second. However, solving instances with  $D = 0.3$ is more difficult for $\IP$ as the time limit is reached on $4$ of these instances, and the time to a solution exceeds $100$ seconds on $4$ others. We remark that $\CCG$ is only able to solve $4$ out of $36$ instances, three of which are infeasible, and hits the time limit on all others. When $\CCG$ reaches the time limit on instances with $n \geq 30$ {\crev and $\density \geq 0.4$} it has applied at least {\cRev $14,000$} cycle separation cuts. Similarly, $\OGcycles$ is only able to solve {\cRev $9$} instances, three of which are infeasible, and $\OGranks$ solves $4$ instances, three of which are infeasible. {\cRev Despite solving fewer instances than $\OGcycles$, $\CCG$ has better solution times on the instances that both formulations are able to solve. For this reason, and because \cite{Omer2017} conclude that of the existing approaches, $\CCG$ performs the best, we proceed to compare with $\CCG$.}

The  CP formulations have better performance than all of the IP formulations. Their detailed comparison is given in Table \ref{table:MinDCP} of Online Supplement \ref{app:MinD_results}. %
$\rankCP$ performs best on instances with $D = 0.4$ or $D = 0.5$, solving most in a fraction of a second. However, if $\rankCP$ cannot solve an instance almost immediately, it cannot solve it in the time limit. As $\rankCP$ solves only $20$ instances in total, we deem its performance worse overall than that of $\vertexCP$ and $\combCP$. {\crev When valid inequalities and variable fixing are added to $\vertexCP$ the time required to solve $25$ of the total of $32$ solved instances decreases. Similarly, when  valid inequalities and variable fixing are added to $\combCP$ the time decreases for $24$ of $33$ instances. Thus we conclude that the $\vertexCP$ and $\combCP$ formulations with valid inequalities outperform those without.} 

As $n$ increases, $\vertexCP$ {\crev with valid inequalities} is able to solve instances with $D = 0.4$ and $D = 0.5$ in a fraction of a second, however, it takes significantly longer to solve instances with $D = 0.3$ reaching the time limit for {\crev $4$} of these lowest density instances. $\combCP$ {\crev with valid inequalities} reaches the time limit when solving $3$ instances all of which have $D = 0.3$, but for all other instances with $n\geq 25$ and $D = 0.4$ or $D = 0.5$ the solution times are competitive. Finally, we compare the extended formulation $\efExtensive$ which is also a CP, it is outperformed by the other CPs on all but $2$ instances and hits the time limit on $12$ instances. We conclude that  $\vertexCP$ and $\combCP$ {\crev with valid inequalities} have the best CP performance. 

\begin{figure}[h]
	\centering
	\includegraphics[width=0.75\textwidth]{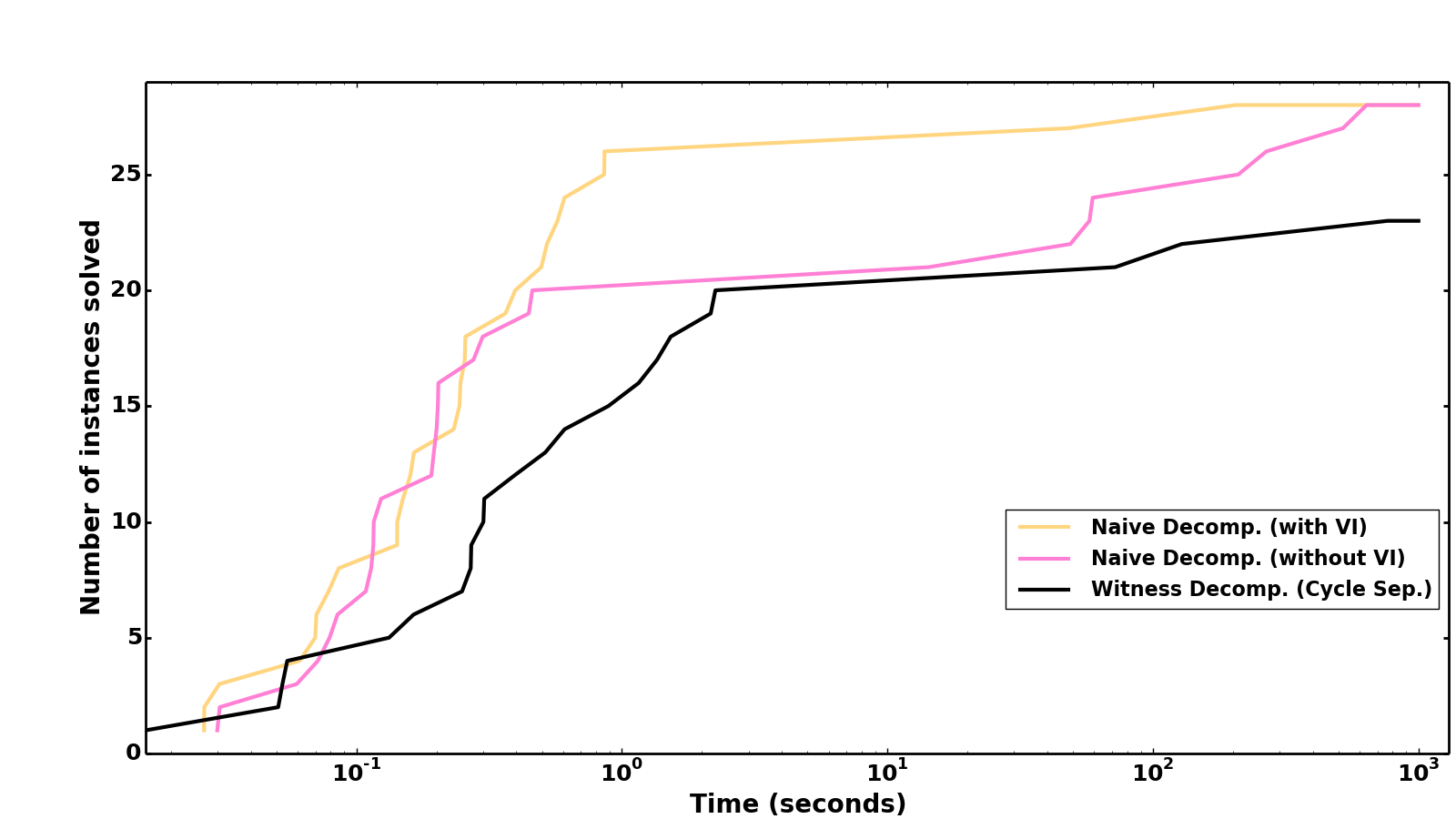}
	\caption{Solution times of the decompositions on random instances.}
	\label{fig:minDdecomp}
\end{figure}

To show the strength of the valid inequalities and variable fixing for the naive decomposition, we compare the naive decomposition with and without these enhancements and the best performing witness-based decomposition in Figure \ref{fig:minDdecomp}. We were able to detect at least one valid inequality or apply variable fixing for all random instances. The results (provided in detail in Table \ref{table:MinDVI} of Online Supplement \ref{app:MinD_results})
show that the naive decomposition reaches the time limit on the same {\crev$7$} instances with or without  the addition of the valid inequalities. We also remark that {\crev for most instances} the naive decomposition with valid inequalities has fewer cuts, and because the majority of the time to find a solution is spent finding an IIS to generate cuts, the  naive decomposition with valid inequalities is faster than without.

We consider a variety of options with which to solve the witness-based decomposition, first we note that because we would like to break cycles in $\witnessvar$, we can add cycle breaking constraints to $\efExtensive$ before beginning the search. In Table  \ref{table:MinDWitness}  of Online Supplement \ref{app:MinD_results}, 
 we compare  the witness-based decomposition using an IIS with no cycle breaking in $\efExtensive$, breaking only $2$-cycles, and breaking both $2$ and $3$-cycles.  Of these the best performance is from breaking both $2$ and $3$-cycles  which we also implement with the cycle separation. We add the valid inequality from Proposition \ref{eq:witnessVI} to all formulations. Naturally, finding an IIS is slower than separating a cycle using depth-first search, our results indicate that the cycle separation is several orders of magnitude faster than finding an IIS (these negligible times have been omitted from Table \ref{table:MinDWitness}).
   In terms of the number of cuts, neither the IIS nor the cycle separation version dominates, but the cycle separation version  is able to solve $3$ more instances than the IIS version within the time limit. We also remark that the cycle separation version has the fastest solution time on $20$ instances.

\begin{figure}[h]
	\centering
	\includegraphics[width=0.75\textwidth]{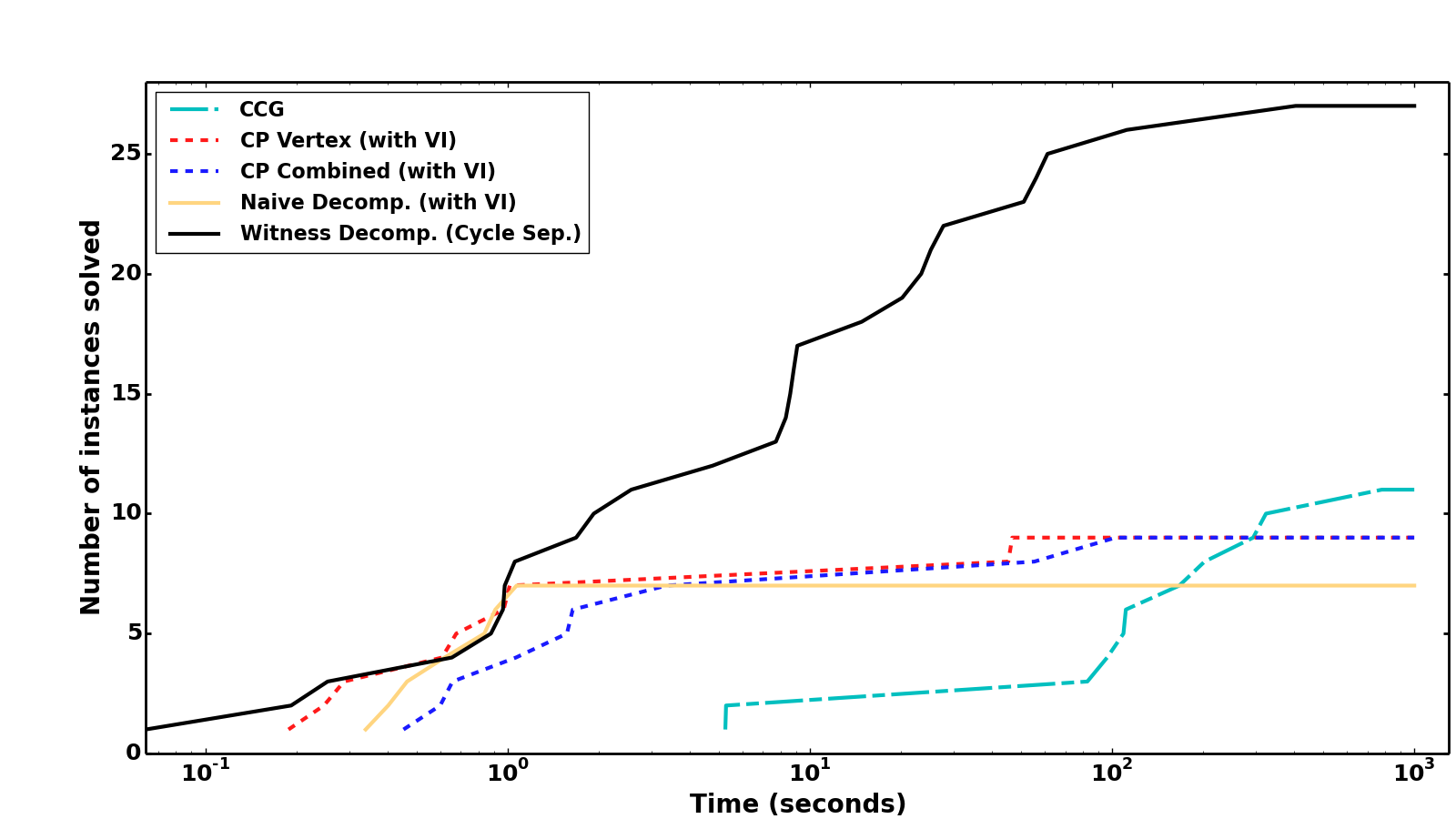}
	\caption{\cRev Solution times of the best performing models on synthetic instances.}
	\label{fig:MinDsyn}
\end{figure}

Figure \ref{fig:MinDsyn} gives the solution times for the best performing methods on the synthetic instances, as well as the best method from the literature,  $\CCG$. We observe that for synthetic instances the witness-based decomposition with cycle separation outperforms the other formulations both in speed and in terms of the number of instances solved.
Table \ref{table:Syn} of Online Supplement \ref{app:MinD_results}
shows the results for synthetic instances. For this test set, $\CCG$ reaches the time limit on {\cRev $16$} instances. It struggles in particular when the number of nodes increases, as it only solves $4$ instances with $n\geq30$. The two CP formulations $\vertexCP$ and $\combCP$  {\crev with valid inequalities}, have very similar performance, reaching the time limit on $18$ instances; those with both a large number of nodes and of double vertices. On synthetic instances, the naive decomposition with valid inequalities reaches the time limit on the most instances, $20$ in total. This decomposition is however, faster than the witness-based decomposition on instances that they both solve. The naive decomposition stays at the root node for all instances as they are hitting the time limit while trying to find an IIS to generate a cut. Finally, the witness-based decomposition is the best performing overall, it is the fastest to a solution for all {\crev but four instances and only reaches the time limit on a single instance}. 

{\crev
Figure \ref{fig:MinDprot} gives the solution times for the best performing methods on pseudo-protein instances. Similarly to the synthetic instances, the best performing method is the witness-based decomposition, it is able to solve $222$ of the $399$ instances within the time limit. $\CCG$ is able to solve {\cRev $148$} instances within the time limit. Similar to its performance on the random instances, $\rankCP$ dominates all other methods in the beginning as it is able to solve $57$ pseudo-protein instances very quickly, but is unable to solve any more instances after that. $\vertexCP$ and $\combCP$ with valid inequalities have similar performance, solving $58$  and $57$ instances within the time limit, respectively. Both $\IP$ and the naive decomposition with valid inequalities were outperformed by all other methods, and were only able to solve $57$ and $46$ instances within the time limit, respectively.

\begin{figure}[h]
	\centering
	\includegraphics[width=0.75\textwidth]{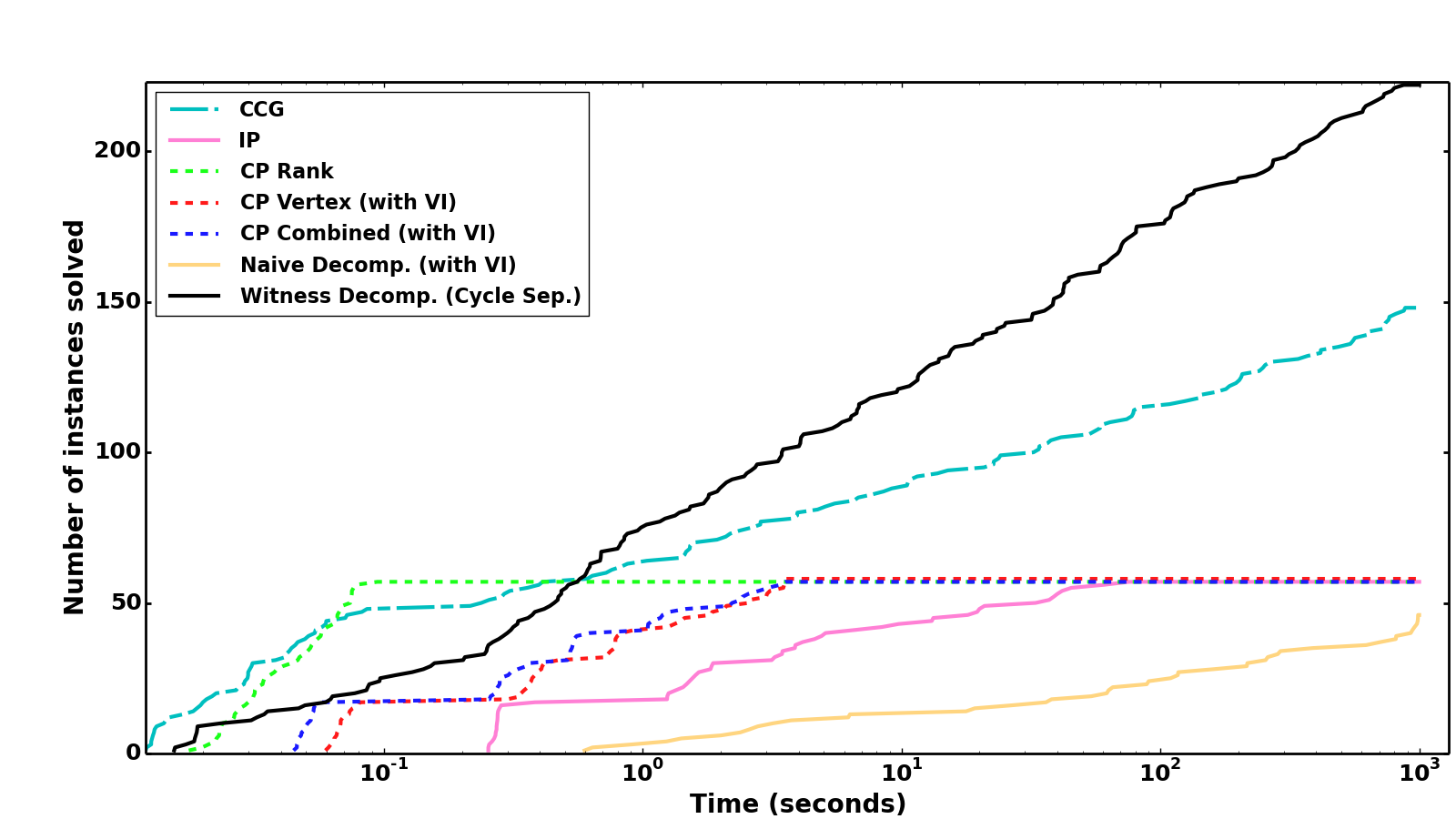}
	\caption{\cRev Solution times of the best performing models on pseudo-protein instances.}
	\label{fig:MinDprot}
\end{figure}}
\section{Conclusion} \label{sec:conclusion}

We propose the first CP formulations and hybrid decomposition methods for $\MinD$ as well as a novel IP. We provide the first valid inequalities for $\MinD$ and use polyhedral theory to prove the correctness of the witness-based decomposition.  

Our computational results show $\vertexCP$ and $\combCP$ with valid inequalities perform best on random instances, but still struggle with the lowest density instances in this data set. For both synthetic {\crev and pseudo-protein} instances, we observe that the witness-based decomposition has superior performance to all other formulations. We conclude {\crev scalability is still an issue for this problem as the most difficult instances to solve are those in the pseudo-protein dataset with large number of nodes. We also remark that for all instance sets studied there exists a novel formulation in this work that outperforms those of the literature. }

\section*{Acknowledgments}
\crev We are grateful to Leo Liberti for introducing us to the distance geometry area and this particular topic. Also, we would like to thank J\'{e}r\'{e}my Omer for kindly providing the pseudo-protein instances similar to the ones used in Omer and Gon\c{c}alves (2017). 

\bibliographystyle{abbrv} 
\bibliography{bibliography} 

\newpage
\appendix
	{\cRev \section{Integer Programming Model for $\MinN$} \label{appendix:IPminNodes}
	
This model is an extension of the model $\IP$ for $\MinD$, with the addition of constraints \eqref{MNIPvr:cliqueNodes}-\eqref{MNIPvr:NodesDouble} and with a different objective function.
Define binary variables $\doublevar_\rankind= 1$ if the vertex at position $\rankind \in [n-1]$ is a double, $0$ otherwise. Define binary variables $ \vrvar_{\vertind \rankind}=1$ if the vertex $\vertind \in \vertexset$ is at rank $\rankind \in [n-1]$, $0$ otherwise. Finally, define positive integer variables $m_r$ equal to the number of nodes at level $r$ of the BP tree. We also introduce binary indicator variables $\ipindicator_{\vertind\rankind}$ for all $\vertind \in \vertexset, \ \rankind \in [n-1]$ to express logical constraints. 

\bsubeq 
\begin{alignat}{2}
\min \ &  \sum_{\rankind \in [n-1]} m_\rankind && \label{MNIPvr:obj} \\
\text{s.t.} \ &\sum_{\rankind \in [n-1]}\vrvar_{\vertind\rankind} = 1 && \forall \  \vertind \in \vertexset \label{MNIPvr:vertex1rank}\\
&\sum_{\vertind \in V} \vrvar_{\vertind\rankind} = 1 && \forall \ \rankind \in [n-1] \label{MNIPvr:rank1vertex}\\
& \sum_{u \in \neighb(\vertind)} \sum_{j \in [\rankind -1]} \vrvar_{uj} \geq \rankind\vrvar_{\vertind\rankind} \quad &&  \forall \ \vertind \in \vertexset, \rankind \in [1,\K] \label{MNIPvr:initalClique}\\
&\sum_{u \in \neighb(\vertind)} \sum_{j \in [\rankind -1]} \vrvar_{uj} \geq \K \vrvar_{\vertind \rankind} &&  \forall \ \vertind \in \vertexset, \rankind \in [\K+1, n-1] \label{MNIPvr:Kpred}\\
& \doublevar_\rankind = 0 && \forall \ \rankind \in [\K-1] \label{MNIPvr:cliqueDoubles} \\
&  \doublevar_\K =1 && \label{MNIPvr:firstDouble}\\
& \sum_{u \in \neighb(\vertvar)} \sum_{j \in [\rankind -1]}\vrvar_{uj}    \geq (K+1)\ipindicator_{\vertind\rankind} \quad &&   \forall \ \vertind \in \vertexset, \rankind \in [\K, n-1] \label{MNIPvr:logical1}\\ 
&  \vrvar_{\vertind\rankind} -  \doublevar_\rankind    \leq \ipindicator_{\vertind\rankind} \quad &&   \forall \ \vertind \in \vertexset, \rankind \in [\K, n-1]\label{MNIPvr:logical2}\\ 
& m_\rankind = 1 && \forall \ \rankind \in [\K-1] \label{MNIPvr:cliqueNodes} \\
&  m_\rankind \geq m_{\rankind-1}  && \forall \ \rankind \in [\K, n-1]  \label{MNIPvr:NodesIncrease}\\
&  m_\rankind - 2m_{\rankind-1} \geq -2 ^{\rankind - \K}(1-\doublevar_{\rankind-1})  && \forall \ \rankind \in [\K, n-1]  \label{MNIPvr:NodesDouble}\\
&\vrvar \in \{0,1\}^{|\vertexset| \times n},  \doublevar   \in \{0,1\}^n, \ipindicator \in  \{0,1\}^{|\vertexset| \times n}, m \in \Z_+^n && 
\end{alignat}
\esubeq
Constraints \eqref{MNIPvr:vertex1rank} and \eqref{MNIPvr:rank1vertex} ensure that there is a bijection from the vertices to the ranks, that is each vertex has exactly one rank and vice versa. Constraints \eqref{MNIPvr:initalClique} ensure that we have an initial clique of size $\K +1$, since each vertex in the clique will be adjacent to all its predecessors. Constraints \eqref{MNIPvr:Kpred} ensure that each vertex after the initial clique has at least $\K$ adjacent predecessors. Constraints \eqref{MNIPvr:cliqueDoubles} and \eqref{MNIPvr:firstDouble} fix the double value of positions $[\K]$. Constraints \eqref{MNIPvr:logical1} and  \eqref{MNIPvr:logical2} link the $\vrvar_{\vertind\rankind} $ and $ \doublevar_\vertind$ variables, where we want to enforce that if vertex $\vertind$ is a non-double, then for any rank $\rankind$, either $\vertind$ is not in position $\rankind$ or it has at least $\K + 1$ adjacent predecessors. Constraints \eqref{MNIPvr:cliqueNodes} fix the first $\K-1$ layers to have $1$ node as follows from the problem definition. Constraints \eqref{MNIPvr:NodesIncrease} ensure the number of nodes at each level is non-decreasing. Finally, the constraints \eqref{MNIPvr:NodesDouble} enforce the logical constraint $\doublevar_{\rankind-1} = 1 \implies m_\rankind \geq 2m_{\rankind-1}$. Since we are minimizing the sum over the number of nodes at each level, we do not need to also enforce $\doublevar_{\rankind-1} = 1 \implies m_\rankind \leq 2m_{\rankind-1}$, as at an optimal solution this constraint is satisfied. The big-M is active in these constraints when $\doublevar_{\rankind-1} = 0 $, the maximum difference is given when every position before $\rankvar$ is a double, thus $2^{\rankvar-\K}$. 
}

\section{Details of the Existing Models} \label{appendix:existing}

Prior to this work, \cite{Omer2017} present two IP formulations and one branch-and-cut procedure for $\MinD$.

\subsection{The cycles IP}
Define binary precedence variables $\precvar_{ij} = 1$ if and only if vertex $i \in \vertexset$ precedes vertex $j \in \vertexset$ in the order, double variables  $\doublevar_\vertind =1$ if and only if vertex $\vertind \in \vertexset$ is double or is in the first clique. {\crev Let  $\clique$ denote the set of all possible $\K$-cliques that can be extended to a $(\K + 1)$-clique, {\cRev note these cliques are not ordered}. Define binary the initial clique selection variables $\cliquevar_c = 1$ if and only if clique $c \in \clique$,  is the initial clique. } Let $R_i^c $ be the rank from $[1,K]$ of vertex $i$ in clique $c$. Then, $\MinD$ can be formulated as follows:
\bsubeq \label{form:OGcycles}
\begin{alignat}{2}
\hspace*{-0.3cm} \OGcycles: \min \ &  \sum_{\vertind \in \vertexset} \doublevar_\vertind - \K&& \label{OGcycles:obj}\\
\text{s.t.} \ &\precvar_{ij} + \precvar_{ji} = 1&& \forall \  i, j \in \vertexset, i \neq j  \label{OGcycles:2cycle}  \\
& \precvar_{ij} + \precvar_{jk} + \precvar_{ki} \leq 2 \quad && \forall \ \{i,j\} \in \edgeset, k \in \vertexset,  k \neq i,j    \label{OGcycles:3cycle} \\
& \sum_{c \in \mathcal{K}}  \cliquevar_c =1   \label{OGcycles:sumClique} \\
& \sum_{\{i,j\} \in \edgeset} \precvar_{ji} \footnotemark + \hspace*{-0.1cm} \sum_{c \in \mathcal{K}: i \in c} (\K-R_i^c+1)\cliquevar_c \geq \K + (1 - \doublevar_i) \quad&& \forall \ i \in \vertexset  \label{OGcycles:linking}  \\ 
& \precvar_{ij} \in \{0,1\} && \forall  \  i, j \in \vertexset, i \neq j  \label{OGcycles:domPrec} \\
& \cliquevar_c  \in \{0,1\} && \forall \  c \in \mathcal{K}  \label{OGcycles:domClique} \\
&\doublevar_\vertind \in \{0,1\} && \forall \ \vertind \in \vertexset \label{OGcycles:domDouble} 
\end{alignat}
\footnotetext{We believe there is a typo in the original paper \cite{Omer2017}, which incorrectly writes $ \precvar_{ij} $.}
\esubeq
Objective \eqref{OGcycles:obj} minimizes the number of double vertices in the order, the constant term is due to  \eqref{OGcycles:linking} and will be explained shortly. Constraints \eqref{OGcycles:2cycle}, \eqref{OGcycles:3cycle} are the usual linear ordering constraints used to break cycles of size $2$ and $3$ \cite{Grotschel1984} and sufficient to determine a vertex ordering \cite{coudert2016}. Constraint \eqref{OGcycles:sumClique} ensures we choose exactly one initial $\K$-clique in the order. Constraints \eqref{OGcycles:linking} are logical constraints which enforce the appropriate number of adjacent predecessors for double or non-double vertices. Fixing a vertex $i \in \vertexset$, $\sum_{\{i,j\} \in \edgeset} \precvar_{ji}$ gives the number of adjacent predecessors of $i$. If $i$ is not in the chosen initial clique, the second left hand side summation is eliminated, and $i$ has exactly $\K$ adjacent predecessors if it is a double and,  $\doublevar_i = 1$ and otherwise it has at least $\K + 1$, as mentioned in the discussion of (2c). 
If vertex $i$ is inside the selected initial clique, the second left hand side summation becomes $(\K-R_i^c+1)$ and rearranging the terms gives that $i$ has $R_i^c - \doublevar_i$ adjacent predecessors. Since each vertex in the initial clique must be adjacent to all the others in the clique and, $i$ is the $(R_i^c)^{\text{th}}$ vertex in the clique, $i$ must have $R_i^c -1$ adjacent predecessors, meaning $\doublevar_i = 1$ for all $i$ in the chosen initial clique. In order to be consistent across all formulations, we subtract these $\K$ doubles from the objective \eqref{OGcycles:obj}, since by definition the vertices in the first clique are not double. Finally, constraints \eqref{OGcycles:domPrec}, \eqref{OGcycles:domClique}, \eqref{OGcycles:domDouble} give binary domains for all the  decision variables.  

For $\OGcycles$ we have \[\rankfunc(\vertind) = \sum_{\mathclap{\substack{u \in \vertexset \\ u \neq v}}} \precvar_{uv}\] \[\doublefunc(\vertind) = \doublevar_\vertind, \forall \ \vertind \in \vertexset.\]

\subsection{The rank IP}
Define integer variables $\rankvar_i \in [n-1]$ for the rank of vertex $i \in V $ in the order.
\bsubeq \label{form:OGranks}
\begin{alignat}{2}
\OGranks: \min \ &  \sum_{\vertind \in \vertexset} \doublevar_\vertind -\K && \label{OGranks:obj}\\
\text{s.t.} \  &n \, \precvar_{ij} + (\rankvar_i+1) - (\rankvar_j +1) \leq n-1 && \forall \ \{i,j\} \in \edgeset  \label{OGranks:cycle}\\
& \sum_{c \in \mathcal{K}}  \cliquevar_c =1  \label{OGranks:sumClique}\\
& \sum_{\{i,j\} \in \edgeset} \precvar_{ji} + \sum_{c \in \mathcal{K}: i \in c} (\K-R_i^c+1)\cliquevar_c \geq \K + (1 - \doublevar_i) \quad&& \forall \ i \in \vertexset  \label{OGranks:linking}\\ 
& \precvar_{ij} \in \{0,1\} && \forall  i, j \in \vertexset, i \neq j \label{OGranks:domPrec}\\
& \cliquevar_c  \in \{0,1\} && \forall c \in \mathcal{K} \label{OGranks:domClique}\\
&\doublevar_\vertind \in \{0,1\} && \forall \ \vertind \in \vertexset  \label{OGranks:domDouble}\\
&\rankvar_\vertind \in [n-1] && \forall \ \vertind \in \vertexset  \label{OGranks:domRank}
\end{alignat}
\esubeq
Objective function \eqref{OGranks:obj}  and constraints \eqref{OGranks:sumClique}, \eqref{OGranks:linking}, \eqref{OGranks:domPrec}, \eqref{OGranks:domClique}, and \eqref{OGranks:domDouble} are as in \eqref{form:OGcycles}. Constraints \eqref{OGranks:cycle} ensure that if vertex $i$ precedes vertex $j$, then the rank of $i$ is strictly less than the rank of $j$. Constraints \eqref{OGranks:domRank} enforce the domain of the rank variables to be one rank for each vertex in the order. Thus, constraints \eqref{OGranks:cycle} and  \eqref{OGranks:domRank} replace the cycle breaking constraints in  $\OGcycles$.

For any feasible solution to $\OGranks$ we have $\rankfunc(\vertind) = \rankvar_\vertind $, and $\doublefunc(\vertind) = \doublevar_\vertind$.

\subsection{Cycle Constraint Generation $ \CCG $}

The cutting-plane version of $\CCG$ is given in Algorithm \ref{alg:CCG} for ease of exposition, but has been implemented in a branch-and-cut fashion. Algorithm \ref{alg:CCG} takes as input integer $q \geq 2$, it begins by generating  all cycles of size $\leq q$, $\cycles_q$. A master problem MP is then created with the objective \eqref{OGranks:obj} , constraints \eqref{OGranks:sumClique}, \eqref{OGranks:linking}, and cycle breaking constraints which will break all cycles in $\cycles_q$. It then passes the master problem solution to a function which will detect if there is cycle in $\bar{\precvar}$. If the cycle is found, we can add a cycle cut to MP, otherwise we have an optimal solution and are done. We denote the set of cycles of a graph as $\cycles$, where a cycle $C \in \cycles$ has vertex set $\vertexset^C$ and edge set $\edgeset^C$, i.e., $C = (\vertexset^C, \edgeset^C)$.

\begin{algorithm} 
\footnotesize
	\SetAlgoLined
	\KwIn{$q \in \N$ with $q \geq 2$}
	$\cycles_q = \{C \in \cycles : |\vertexset^C| \leq q\}$ \\
	MP := $\bigg\{ \min \displaystyle\sum_{\vertind \in \vertexset} \doublevar_\vertind  - \K: $ \eqref{OGranks:sumClique}, \eqref{OGranks:linking}, and 
	$ \displaystyle\sum_{\{i,j\} \in \edgeset^C} \precvar_{ji} \leq |\vertexset^C| - 1 \ \forall \ C \in \cycles_q	\bigg\}$ \\
	optimal := false \\
	\While{not optimal}{
		solve MP and get solution$(\bar{\precvar}, \bar{\doublevar}, \bar{\cliquevar})$\\
		$C$ := SeparateCycle$(\bar{\precvar})$\\
		\eIf{$C$ is found}{
			add $\displaystyle\sum_{\{i,j\} \in \edgeset^C} \precvar_{ji} \leq |\vertexset^C| - 1$ to MP		
		}{
			accept $(\bar{\precvar}, \bar{\doublevar}, \bar{\cliquevar})$\\
			optimal := true 
		}
	}
	\caption{Cycle Constraint Generation $\CCG$ Procedure}
	\label{alg:CCG}
\end{algorithm}
\noindent As for $\mathbb{(CYCLES)}$,  \[\rankfunc(\vertind) = \sum_{\mathclap{\substack{u \in \vertexset : \\ u \neq v}}} \precvar_{uv}\]  for $\CCG$. 

\cite{Omer2017}  show that Cycle Cut Generation outperforms both IP formulations. Nevertheless, we will compare our formulations against all three of their approaches for completeness. 

	\newpage
\begin{landscape}
\section{Summary of the paper}
\label{app:summary}
	\renewcommand\arraystretch{0.8}
	\setlength{\LTleft}{-11pt}	
	\small 
\begin{longtable}{lllllll}
		\caption{Summary of  the formulations for $\MinD$.}
	\label{table:MinDsummary}\\
	\toprule
\multicolumn{7}{l}{\textbf{$\MinD$}: Minimize the number of vertices with exactly $\K$ adjacent predecessors in the order.} \\
\midrule
{\cellcolor{gray!20!white}\textsc{Literature}} & Variables & Domain & Number & Constraints & Number & Comments \\
Cycles $\OGcycles$& $\doublevar_v$ & $\{0,1\}$ & $|\vertexset|$ & linear ordering & $|\vertexset|^2 + |\vertexset|^2\times|\edgeset|$ \\
 & $\cliquevar_v$ & $\{0,1\}$ & $|\vertexset|$ & clique selection & 1 \\
 & $\precvar_{uv}$ & $\{0,1\}$ & $|\vertexset|^2$ & linking & $|\vertexset|$ \\
Ranks $\OGranks$& $\doublevar_v$ & $\{0,1\}$ & $|\vertexset|$ & linear ordering & $|\edgeset|$ \\
 & $\cliquevar_v$ & $\{0,1\}$ & $|\vertexset|$ & clique selection & 1 \\
 & $\precvar_{uv}$ & $\{0,1\}$ & {\crev $2|\edgeset|$} & linking & $|\vertexset|$ \\
 & $\rankvar_v$ & $[n-1]$ & $n$ \\
Cycle cut generation $\CCG$ & $\doublevar_v$ & $\{0,1\}$ & $|\vertexset|$ & clique selection & 1 & Generate cycles in a  \\
 & $\cliquevar_v$ & $\{0,1\}$ & $|\vertexset|$ & linking & $|\vertexset|$ & branch-and-cut \\
 & $\precvar_{uv}$ & $\{0,1\}$ & {\crev $2|\edgeset|$} & cycle breaking cuts & & procedure \\
{\cellcolor{gray!20!white}\textsc{New}} \\
Vertex-rank IP $\IPVR$  & $\doublevar_r$ & $\{0,1\}$ & $n$ & 1-1 assignment & $2n$ \\
 & $\vrvar_{\vertind \rankind}$ & $\{0,1\}$ & $|\vertexset|^2$ & clique  & $n^2$ \\
 &  &  &  & fixing & $\K+1$ \\
 &  &  &  & linking & $2(|\vertexset|\times n-\K-1)$ \\
CP Rank $\rankCP$& $\doublevar_v$ & $\{0,1\}$ & $n$ & AllDifferent & 1 \\
 & $\rankvar_\vertind$ & $[n-1]$ & $n$ & clique & $|\vertexset|(|\vertexset|-1)/2-|\edgeset|$ \\
 &  &  &  & logical & $|\vertexset|$ \\
CP Vertex $\vertexCP$& $\doublevar_r$ & $\{0,1\}$ & $n$ & AllDifferent & 1 \\
 & $\vertvar_\rankind$ & $[|V|-1]$ & $n$ & clique &  $\K(\K+1)/2$ \\
 &  &  &  & fixing & $\K+1$ \\
 &  &  &  & logical & $n-\K-1$ \\
CP Combined $\combCP$ & $\doublevar_r$ & $\{0,1\}$ & $n$ & inverse & 1 & Combines CP Rank  \\
 & $\rankvar_\vertind$ & $[n-1]$ & $n$ & clique & $(|\vertexset|(|\vertexset|-1) + \K(\K+1))/2-|\edgeset|$ \\
 &$\vertvar_\rankind$ & $[|V|-1]$ & $n$ & logical & $|\vertexset|+ n-\K-1$ \\
 &  &  &  & fixing & $\K+1$ \\
Naive Decomposition & $\doublevar_r$ & $\{0,1\}$ & $n$ & fixing & $\K+1$ & MP (IP): fixes doubles \\
 & $\vertvar_\rankind$ & $[|V|-1]$ & $n$ & AllDifferent & 1 & SP (CP): finds an order \\
 &  &  &  & clique & $\K(\K+1)/2$ \\
 &  &  &  & fixing and logical & $n$ \\
Witness-based Decomposition & $\doublevar_v$ & $\{0,1\}$ & $n$ & clique selection & 1 & MP (IP): fixes doubles,\\
 & $\cliquevar_v$ & $\{0,1\}$ & $|\vertexset|$ & clique witness & $|\vertexset|(|\vertexset|-1)/2+|\edgeset|$ &first clique, witnesses \\
 & $\witnessvar_{uv}$ & $\{0,1\}$ & $2|\edgeset|$ & witness & $|\vertexset|$ & SP (CP): finds an order\\
 & $\rankvar_\vertind$ & $[n-1]$ & $n$ & AllDifferent & 1 &  $\efExtensive$ combines MP and\\
 &  &  &  & rank & $|\vertexset|+ 2|\edgeset|$ &  SP into one CP \\
\bottomrule
\end{longtable}
\end{landscape}
\clearpage 
	{\crev \section{Numerical Comparison of $\CCG$ and Witness-Based Decomposition} \label{app:ccg_wb}
In this section, we present a brief numerical comparison of $\CCG$ and the witness-based decomposition on select random instances with more than one double  vertex in the optimal solution. As mentioned in Section 4.3, the formulations are similar in that they both incorporate ideas from linear ordering formulations by defining variables on the graph edges and adding cycle breaking cuts in an iterative manner. The key difference in the formulations is in their treatment of the initial clique, while $\CCG$ enumerates all possible ordered $\K$-cliques that can be extended to a $(\K + 1)$-clique {\cRev and creates a variable for each clique}, the witness-based decomposition defines indicator variables for the initial clique and constrains the vertices in the first $(\K + 1)$ positions to form a $(\K + 1)$-clique. In particular, the exact position of these vertices within the clique is not assigned. Table \ref{table:CCGWBfixed} presents the results of running both formulations when the initial $\K$-clique has been fixed to an optimal one. 

The columns of Table \ref{table:CCGWBfixed} are:
\BI
\I``Time": The solution time in seconds if the instance is solved in the given time limit, ``TL" if the instance hit the time limit.
\I``BB Nodes": The number of branch-and-bound nodes explored (for the IP formulations); exact if it is less than one million, lower bound rounded to the closest million otherwise where a single decimal point is used up to ten million for better accuracy.
\I`` \# Cuts": The number of cuts added in the branch-and-cut procedure.
\I ``$n$" $ \in \{ 20, 25, 30, 35 \}$: The number of vertices in the input graph.
\I ``$\density$"$ \in \{ 0.3, 0.4, 0.5\}$: The edge density of the input graph.
\I ``Inst.": The assigned instance number from $\{1,2,3\}$ for each $(n,\density)$ combination.
\I ``Obj": Optimal objective  value of the instance; $+\infty$ if it is not solved by any method.
\EI
As expected, the witness-based decomposition is able to solve four of the six instances faster than without the clique fixing, $\CCG$ however is unable to solve any instance within the time limit.

\begin{longtable}{cccr rrr rrr}
	\caption{\small Results of $\CCG$ and the witness-based decomposition on selected random instances with an optimal first clique fixed.}\\
	\label{table:CCGWBfixed} \\
	\toprule 
	&  &  &  & \multicolumn{3}{c}{$\CCG$}   &\multicolumn{3}{c}{Witness-based Decomp.}   \\
	\cmidrule(lr){5-7}
	\cmidrule(lr){8-10}
	$n$ & $\density$ & Inst. & Obj. & Time & BB Nodes & \# Cuts & Time & BB Nodes & \# Cuts \\
	\midrule
	20 & 0.4 & 1 & 2 & TL & $>$1.5M & 827  & \textBF{1.53} & 6951 & 245 \\
	&  & 2 & 3 & TL &  $>$1.1M	& 1337 & \textBF{34.43} & 80141 & 533 \\
	&  & 3 & 4 & TL & $>$1M & 1605 & \textBF{270.00} & 370348 & 1522 \\
	25 & 0.3 & 2 & 5 & TL & $>$1.7M & 485 & \textBF{10.51} & 28141 & 210 \\
	&  & 3 & 7 & TL & $>$1.1M & 1015 & TL & $>$1.2M & 985 \\
	& 0.4 & 2 & 2 & TL & 184450 & 10813 & TL & 344436 & 9140 \\ 
	\bottomrule
\end{longtable}

{\cRev We also compare the number of initial cliques considered by the master problem during the witness-based decomposition and the number of cliques enumerated for $\CCG$. That is, the unique initial clique candidates the master problem of the witness-based decomposition passes to the subproblem and the number of clique variables, $\kappa_c$, created in $\CCG$. } We focus on the three instances for which the witness-based decomposition was able to find an optimal solution within the time limit. The results are presented in Table \ref{table:CCGWBcliques}, which has similar columns to Table \ref{table:CCGWBfixed} except gives the number of cliques generated or enumerated in each method.

\begin{longtable}{cccr rr}
	\caption{\small Comparison of the number of cliques considered by $\CCG$ and  the witness-based decomposition.}\\
	\label{table:CCGWBcliques} \\
	\toprule 
	$n$ & $\density$ & Inst. & Obj. & \multicolumn{1}{l}{Cliques} & \multicolumn{1}{l}{Cliques} \\
	& & & & \multicolumn{1}{l}{$\CCG$}  & \multicolumn{1}{l}{Witness-based Decomp.}\\
	\midrule
	20 & 0.4 & 1 & 2 &56 &18  \\ \
	&  & 2 & 3 &38 & 10 \\
	25 & 0.3 & 2 & 5 & 27& 8\\
	\bottomrule
\end{longtable}

We observe that the witness-based decomposition considers far fewer cliques than $\CCG$. This may help to explain the success of witness-based decomposition over $\CCG$.
}
	\vspace*{0.5cm}
	\section{Proof of Proposition 1} \label{app:proof}

We prove the validity of the inequality $\doublevar_\vertind \leq 1 - \cliquevar_\vertind$ for $\efExtensive$.

Without loss of generality, we fix $\vertind \in \vertexset$. If $\cliquevar_\vertind = 0$, the inequality is redundant, since $\doublevar_\vertind$ is binary. If $\cliquevar_\vertind = 1$, we would like $\doublevar_\vertind =0$, since a vertex in the initial clique is not a double. By constraint 
(14e), when $\cliquevar_\vertind = 1$, $\vertind$ will have $\doublevar_\vertind + \K$ witnesses. Thus, since we are minimizing the number of doubles, the objective function will force $\doublevar_\vertind =0$. Meaning we have $\doublevar_\vertind =0$ when $\cliquevar_\vertind = 1$, and so the inequality holds.

	\vspace*{0.5cm}
	\section{$\MinD$ Results} \label{app:MinD_results}

We provide the following tables\footnote{The results tables for the pseudo-protein instances are available upon request. }:
\BI
\I Table \ref{table:MinDIP} compares the IP formulations, including the newly proposed 
formulation $\IP$, as well as the existing formulations from the literature the cycle cut generation procedure {\crev $\CCG$, the cycles formulation $\OGcycles$, and the ranks formulation $\OGranks$} on random instances.
\I Table \ref{table:MinDCP} compares the four novel CP formulations {\crev $\vertexCP$ with and without valid inequalities}, $\rankCP$, {\crev $\combCP$ with and without valid inequalities}, and $\efExtensive$ on random instances.
\I Table \ref{table:MinDVI} compares the {\crev naive decomposition with and without the valid inequalities.} 
\I Table \ref{table:MinDWitness} compares the witness-based decomposition under different options, initializing $\efMP$ with 2- and 3-cycle breaking, and generating cuts using an IIS or the cycle separator. 
\I Table \ref{table:MinDBest} compares the best performing formulations from the previous tables, namely  {$\vertexCP$ with valid inequities, $\combCP$ with valid inequities}, the naive decomposition with valid inequalities, and the witness-based decomposition with 2-cycle breaking and 3-cycle breaking added to $\efMP$  and using cycle separator. 
\I Table \ref{table:Syn} compares the {\crev best performing formulations} on synthetic instances.
{\crev \I Table \ref{table:K4} compares the best performing formulations on feasible random instances with $\K=4$.
\I  Table \ref{table:K5} compares the best performing formulations on feasible random instances with $\K=5$.}
\EI

For the random instances, we have:
\BI
\I ``$n$" $ \in \{ 20, 25, 30, 35 \}$: The number of vertices in the input graph.
\I ``$\density$"$ \in \{ 0.3, 0.4, 0.5\}$: The edge density of the input graph.
\I ``Inst.": The assigned instance number from $\{1,2,3\}$ for each $(n,\density)$ combination.
\I ``Obj": Optimal objective  value of the instance; $+\infty$ if it is not solved by any method.
\EI 
For the synthetic instances, we have:
\BI
\I ``$n$" $ \in \{ 20, 25, 30, 35 \}$: The number of vertices in the input graph.
\I  ``Doubles"  $\in \{\lceil 0.1\times n\rceil, \lceil 0.1\times n\rceil+1, \lceil 0.1\times n \rceil+2 \}$: The upper bound on the number of doubles in the input graph, there may be less depending on the noise edges.
\I ``Noise" : The number of extra edges added to  the input graph.
\I ``Obj": Optimal objective  value of the instance; $+\infty$ if it is not solved by any method.
\EI 
We also have:
\BI
\I ``Time": Solution time in seconds if the instance is solved in the given time limit, ``TL" if the instance hit the time limit.
\I ``IIS Time": The time require for the conflict refiner in CP Optimizer to find an IIS.
\I ``BB Nodes": The number of branch-and-bound nodes explored (for the IP formulations); exact if it is less than one million, lower bound rounded to the closest million otherwise where a single decimal point is used up to ten million for a better accuracy.
\I ``Ch.Pts.": The number of choice points (for the CP formulations); the number convention is the same as the ``BB Nodes".
\I `` \# Cuts": The number of cuts added in the branch-and-cut procedure.
\EI 

\afterpage{
\begin{landscape}
\centering\renewcommand\arraystretch{0.8}
	\small
\begin{longtable}{cccr rr rrr rr rr}
	\caption{\cRev Results for IP formulations on random instances.}
	\label{table:MinDIP} \\
\toprule 
 &  &  &  & \multicolumn{2}{c}{$\IP$} & \multicolumn{3}{c}{$\CCG$}   &\multicolumn{2}{c}{$\OGcycles$}  & \multicolumn{2}{c}{$\OGranks$} \\
 \cmidrule(lr){5-6}
 \cmidrule(lr){7-9}
 \cmidrule(lr){10-11}
  \cmidrule(lr){12-13}
$n$ & $\density$ & Inst. &  Obj. &  Time &  BB Nodes &  Time &  BB Nodes &  \# cuts &  Time &  BB Nodes &  Time &  BB Nodes \\
  \midrule
  20 & 0.3 & 1 & 14 & TL & 453417 & \textBF{0.01} & 0 & 0 & 0.03 & 0 & 0.01 & 0 \\
  &  & 2 & $+\infty$ & 0.06 & 0 & \textBF{0.00} & 0 & 0 & 0.00 & 0 & 0.00 & 0 \\
  &  & 3 & $+\infty$ & 43.16 & 15086 & \textBF{0.08} & 258 & 43 & 0.14 & 20 & 0.15 & 398 \\*[0.18cm]
  & 0.4 & 1 & 2 & \textBF{5.74} & 3028 & TL & 934500 & 921 & 89.22 & 20984 & TL & $>$1.1M \\
  &  & 2 & 3 & \textBF{25.53} & 10890 & TL & 787149 & 1172 & 160.21 & 43842 & TL & $>$1.1M \\
  &  & 3 & 4 & \textBF{30.28} & 12493 & TL & 764550 & 1885 & 279.12 & 91610 & TL & $>$1M \\*[0.18cm]
  & 0.5 & 1 & 1 & \textBF{0.09} & 0 & TL & 385745 & 4450 & TL & 217051 & TL & $>$1.3M \\
  &  & 2 & 1 & \textBF{0.09} & 0 & TL & 402470 & 4108 & TL & 205948 & TL & $>$1.5M \\
  &  & 3 & 1 & \textBF{0.09} & 0 & TL & 426618 & 4108 & TL & 208035 & TL & $>$1.2M \\
  \midrule
  25 & 0.3 & 1 & $+\infty$ & 0.11 & 0 & 0.00 & 0 & 0 & 0.00 & 0 & \textBF{0.00} & 0 \\
  &  & 2 & 5 & \textBF{80.57} & 11650 & TL & $>$1.2M & 629 & 549.19 & 57399 & TL & 735923 \\
  &  & 3 & 7 & 581.19 & 95937 & TL & 792761 & 1307 & \textBF{363.11} & 36138 & TL & 690563 \\*[0.18cm]
  & 0.4 & 1 & 1 & \textBF{0.17} & 0 & TL & 328277 & 5389 & TL & 74699 & TL & 998600 \\
  &  & 2 & 2 & \textBF{35.83} & 8983 & TL & 292200 & 6502 & TL & 75560 & TL & 953026 \\
  &  & 3 & 1 & \textBF{0.18} & 0 & TL & 342172 & 5299 & TL & 69395 & TL & 878376 \\*[0.18cm]
  & 0.5 & 1 & 1 & \textBF{0.21} & 0 & TL & 143260 & 16695 & TL & 80320 & TL & $>$1.1M \\
  &  & 2 & 1 & \textBF{0.21} & 0 & TL & 127055 & 17654 & TL & 81223 & TL & $>$1.1M \\
  &  & 3 & 1 & \textBF{0.21} & 0 & TL & 142009 & 17796 & TL & 87235 & TL & 993935 \\
  \midrule
  30 & 0.3 & 1 & 3 & \textBF{160.00} & 20239 & TL & 369653 & 4186 & TL & 31633 & TL & 895019 \\
  &  & 2 & 4 & \textBF{164.67} & 20843 & TL & 327500 & 4050 & TL & 26437 & TL & 756744 \\
  &  & 3 & 4 & TL & 74413 & TL & 262484 & 4662 & TL & 24985 & TL & 745650 \\*[0.18cm]
  & 0.4 & 1 & 1 & \textBF{0.49} & 0 & TL & 149648 & 14657 & TL & 26089 & TL & $>$1M \\
  &  & 2 & 1 & \textBF{0.45} & 0 & TL & 113121 & 19578 & TL & 29055 & TL & 824766 \\
  &  & 3 & 1 & \textBF{0.51} & 0 & TL & 117574 & 17796 & TL & 42971 & TL & 871401 \\*[0.18cm]
  & 0.5 & 1 & 1 & \textBF{0.65} & 0 & TL & 66740 & 31222 & TL & 29502 & TL & 936367 \\
  &  & 2 & 1 & \textBF{0.57} & 0 & TL & 55952 & 32607 & TL & 42067 & TL & $>$1M \\
  &  & 3 & 1 & \textBF{0.59} & 0 & TL & 68639 & 34761 & TL & 29835 & TL & $>$1M \\
  \midrule
  35 & 0.3 & 1 & 3 & TL & 70945 & TL & 133890 & 13011 & TL & 12328 & TL & 652591 \\
  &  & 2 & 3 & TL & 52051 & TL & 149991 & 11774 & TL & 11495 & TL & 713065 \\
  &  & 3 & 3 & \textBF{655.51} & 46969 & TL & 137145 & 12610 & TL & 9883 & TL & 724207 \\*[0.18cm]
  & 0.4 & 1 & 1 & \textBF{1.00} & 0 & TL & 50158 & 36286 & TL & 12080 & TL & 797966 \\
  &  & 2 & 1 & \textBF{0.93} & 0 & TL & 69639 & 34443 & TL & 16788 & TL & 943108 \\
  &  & 3 & 1 & \textBF{0.92} & 0 & TL & 61730 & 35345 & TL & 11998 & TL & 902507 \\*[0.18cm]
  & 0.5 & 1 & 1 & \textBF{1.16} & 0 & TL & 47893 & 43789 & TL & 19223 & TL & $>$1M \\
  &  & 2 & 1 & \textBF{1.13} & 0 & TL & 48441 & 43845 & TL & 17275 & TL & 761910 \\
  &  & 3 & 1 & \textBF{1.12} & 0 & TL & 48151 & 43090 & TL & 18154 & TL & 745675 \\
  \bottomrule
\end{longtable}
\end{landscape}
}
\afterpage{
\begin{landscape}
\centering
	\renewcommand\arraystretch{0.78}
		\small
\begin{longtable}{cccr rr rr rr rr rr rr}
				\caption{Results for the CP formulations on random instances.}
	\label{table:MinDCP} \\
\toprule
 &  &  &  & \multicolumn{2}{c}{$\vertexCP$}&  \multicolumn{2}{c}{$\vertexCP$}&  \multicolumn{2}{c}{$\rankCP$}&  \multicolumn{2}{c}{$\combCP$} &  \multicolumn{2}{c}{$\combCP$ }&  \multicolumn{2}{c}{$\efExtensive$} \\
  &  &  &  & \multicolumn{2}{c}{}&  \multicolumn{2}{c}{(with VI)}&  \multicolumn{2}{c}{}&  \multicolumn{2}{c}{} &  \multicolumn{2}{c}{(with VI)}&  \multicolumn{2}{c}{} \\
\cmidrule(lr){5-6}
\cmidrule(lr){7-8}
\cmidrule(lr){9-10}
\cmidrule(lr){11-12}
\cmidrule(lr){13-14}
\cmidrule(lr){15-16}
$n$ & $\density$ & Inst. & Obj. & Time & Ch.Pts. & Time & Ch.Pts. & Time & Ch.Pts. & Time & Ch.Pts.  & Time & Ch.Pts.  & Time & Ch.Pts. \\
\midrule
20 & 0.3 & 1 & 14 & TL & $>$19M & TL & $>$17M & TL & $>$34M & TL & $>$13M & TL & $>$14M & \textBF{1.80} & 114233 \\
 &  & 2 & $+\infty$ & 1.16 & 26095 & 1.39 & 28232 & \textBF{0.00} & 0 & 0.00 & 0 & \textBF{0.00} & 0 & 0.00 & 0 \\
 &  & 3 & $+\infty$ & 39.82 & 856311 & 34.94 & 739690 & TL & $>$23M & 8.46 & 200036 & 8.27 & 209521 & \textBF{2.88} & 0 \\*[0.18cm]
 & 0.4 & 1 & 2 & 1.92 & 37797 & 0.03 & 451 & TL & $>$30M & 2.36 & 50296 & \textBF{0.02} & 455 & TL & $>$35M \\
 &  & 2 & 3 & 8.95 & 161665 & 0.03 & 436 & TL & $>$39M & 3.55 & 77317 & \textBF{0.02} & 446 & TL & $>$42M \\
 &  & 3 & 4 & 15.78 & 332463 & \textBF{7.12} & 146090 & TL & $>$45M & 10.08 & 223677 & 9.10 & 201655 & TL & $>$39M \\*[0.18cm]
 & 0.5 & 1 & 1 & 0.02 & 186 & \textBF{0.02} & 186 & 0.02 & 1028 & 0.02 & 453 & 0.02 & 453 & 0.61 & 30411 \\
 &  & 2 & 1 & 0.02 & 160 & 0.02 & 160 & \textBF{0.01} & 978 & 0.03 & 463 & 0.03 & 463 & 0.16 & 8222 \\
 &  & 3 & 1 & \textBF{0.02} & 228 & 0.02 & 228 & 0.02 & 1543 & 0.02 & 445 & 0.02 & 445 & 0.14 & 7412 \\
 \midrule
25 & 0.3 & 1 & $+\infty$ & TL & $>$13M & TL & $>$13M & \textBF{0.00} & 0 & 0.02 & 6 & 0.01 & 12 & 0.00 & 0 \\
 &  & 2 & 5 & TL & $>$13M & TL & $>$12M & TL & $>$38M & TL & $>$12M & TL & $>$12M & TL & $>$34M \\
 &  & 3 & 7 & TL & $>$13M & TL & $>$13M & TL & $>$42M & TL & $>$13M & TL & $>$12M & TL & $>$37M \\*[0.18cm]
 & 0.4 & 1 & 1 & 0.09 & 1192 & 0.05 & 447 & \textBF{0.01} & 1122 & 0.04 & 470 & 0.05 & 470 & 0.54 & 27543 \\
 &  & 2 & 2 & 3.39 & 49464 & 0.05 & 465 & 993.79 & $>$45M & 6.93 & 100309 & \textBF{0.03} & 458 & TL & $>$35M \\
 &  & 3 & 1 & 0.24 & 3814 & 0.04 & 451 & \textBF{0.02} & 1613 & 0.04 & 450 & 0.04 & 450 & 5.17 & 213898 \\*[0.18cm]
 & 0.5 & 1 & 1 & 0.08 & 1198 & 0.04 & 467 & \textBF{0.00} & 238 & 0.04 & 471 & 0.05 & 471 & 0.37 & 16623 \\
 &  & 2 & 1 & 0.05 & 529 & 0.04 & 459 & \textBF{0.01} & 1123 & 0.05 & 469 & 0.04 & 469 & 0.87 & 37055 \\
 &  & 3 & 1 & 0.04 & 511 & 0.04 & 462 & \textBF{0.02} & 1502 & 0.05 & 460 & 0.04 & 460 & 0.37 & 15621 \\
 \midrule
30 & 0.3 & 1 & 3 & 31.34 & 358648 & 0.07 & 485 & TL & $>$51M & 25.72 & 296148 & \textBF{0.06} & 469 & TL & $>$30M \\
 &  & 2 & 4 & 35.81 & 423625 & 15.80 & 138956 & TL & $>$60M & 48.10 & 643116 & \textBF{12.84} & 164305 & TL & $>$34M \\
 &  & 3 & 4 & 353.13 & $>$3.7M & 430.65 & $>$4.5M & TL & $>$56M & \textBF{258.87} & $>$3M & 588.08 & $>$6.9M & TL & $>$27M \\*[0.18cm]
 & 0.4 & 1 & 1 & 0.32 & 3185 & 0.08 & 478 & \textBF{0.02} & 1138 & 0.07 & 478 & 0.07 & 478 & 1.05 & 49638 \\
 &  & 2 & 1 & 0.75 & 7566 & 0.08 & 468 & \textBF{0.02} & 1719 & 0.06 & 458 & 0.06 & 458 & 0.55 & 25189 \\
 &  & 3 & 1 & 0.15 & 1680 & 0.07 & 467 & \textBF{0.02} & 1720 & 0.06 & 476 & 0.07 & 476 & 4.43 & 177270 \\*[0.18cm]
 & 0.5 & 1 & 1 & 0.05 & 169 & 0.05 & 169 & \textBF{0.02} & 1138 & 0.06 & 475 & 0.08 & 475 & 0.36 & 16984 \\
 &  & 2 & 1 & 0.06 & 207 & 0.05 & 207 & \textBF{0.02} & 1148 & 0.04 & 181 & 0.04 & 181 & 0.64 & 26581 \\
 &  & 3 & 1 & 0.06 & 203 & 0.05 & 203 & \textBF{0.02} & 1740 & 0.08 & 471 & 0.07 & 471 & 1.23 & 57488 \\
 \midrule
35 & 0.3 & 1 & 3 & 388.29 & $>$3M & 258.45 & $>$1.6M & TL & $>$58M & 985.83 & $>$9.3M & \textBF{209.41} & $>$2.1M & TL & $>$35M \\
 &  & 2 & 3 & 100.39 & 748258 & \textBF{0.12} & 492 & TL & $>$62M & 145.04 & $>$1.4M & 0.12 & 503 & TL & $>$32M \\
 &  & 3 & 3 & 118.75 & 934246 & 0.13 & 480 & TL & $>$58M & 175.48 & $>$1.5M & \textBF{0.12} & 490 & TL & $>$27M \\*[0.18cm]
 & 0.4 & 1 & 1 & 0.13 & 616 & 0.13 & 487 & \textBF{0.03} & 1852 & 0.12 & 489 & 0.12 & 489 & 1.38 & 54887 \\
 &  & 2 & 1 & 0.27 & 2208 & 0.13 & 477 & \textBF{0.03} & 1233 & 0.12 & 492 & 0.11 & 492 & 1.38 & 57876 \\
 &  & 3 & 1 & 0.09 & 208 & 0.11 & 208 & \textBF{0.03} & 1855 & 0.12 & 499 & 0.12 & 499 & 3.99 & 138625 \\*[0.18cm]
 & 0.5 & 1 & 1 & 0.10 & 205 & 0.10 & 205 & \textBF{0.02} & 1123 & 0.09 & 171 & 0.09 & 171 & 0.53 & 21071 \\
 &  & 2 & 1 & 0.08 & 88 & 0.07 & 88 & \textBF{0.04} & 1842 & 0.08 & 168 & 0.09 & 168 & 0.42 & 18364 \\
 &  & 3 & 1 & 0.09 & 170 & 0.10 & 170 & \textBF{0.03} & 1776 & 0.11 & 501 & 0.13 & 501 & 0.72 & 31121 \\
\bottomrule
\end{longtable}
\end{landscape}
}

\afterpage{
\begin{landscape}
\centering
	\renewcommand\arraystretch{0.8}
	\small
	\begin{longtable}{cccr rr rr rr rr l}
			\caption{Results for the naive decomposition with and without valid inequalities (VI) on random instances.}
		\label{table:MinDVI} \\
		\toprule
 &  &  &  & \multicolumn{4}{c}{Naive Decomp. (without VI)}& \multicolumn{5}{c}{Naive Decomp.  (with VI)} \\
\cmidrule(lr){5-8}
\cmidrule(lr){9-13}
$n$ & $\density$ & Inst. & Obj. & Time &  IIS time &  BB Nodes &  \# Cuts & Time &  IIS time &  BB Nodes &  \# Cuts & VI added \\
\midrule
20 & 0.3 & 1 & 14 & TL & 992.91 & 0 & 6 & TL & 991.49 & 0 & 6 & $y_{\K+1} + y_{\K+2}  \geq 1 $ \\
 &  & 2 & $+\infty$ & \textBF{302.18} & 285.99 & 1 & 56 & 597.65 & 549.67 & 1 & 106 & $y_{\K+1} = 1$, $y_{\K+2} = 1$   \\
 &  & 3 & $+\infty$ & TL & 994.63 & 11 & 23 & TL & 970.16 & 0 & 22 & $y_{\K+1} = 1$, $y_{\K+2} = 1$   \\*[0.18cm]
 & 0.4 & 1 & 2 & 14.28 & 12.87 & 0 & 1 & \textBF{0.23} & 0 & 0 & 0 & $y_{\K+1} + y_{\K+2}  \geq 1 $ \\
 &  & 2 & 3 & 48.87 & 45.35& 0 & 2 & \textBF{0.03} & 0 & 0 & 0 & $y_{\K+1} = 1$, $y_{\K+2} = 1$   \\
 &  & 3 & 4 & 57.65 & 55.11 & 0 & 3 & \textBF{48.10} & 44.18 & 0 & 1 & $y_{\K+1} = 1$, $y_{\K+2} = 1$   \\*[0.18cm]
 & 0.5 & 1 & 1 & 0.08 & 0 & 0 & 0 & \textBF{0.08} & 0 & 0 & 0 \\
 &  & 2 & 1 & 0.03 & 0 & 0 & 0 & \textBF{0.03} & 0 & 0 & 0 \\
 &  & 3 & 1 & 0.03 & 0 & 0 & 0 & \textBF{0.03} & 0 & 0 & 0 \\ \midrule
25 & 0.3 & 1 & $+\infty$ & TL & 983.07 & 0 & 1 & TL & 980.37 & 0 & 1 & $y_{\K+1} + y_{\K+2}  \geq 1 $ \\
 &  & 2 & 5 & TL & 734.242 & 0 & 4 & TL & 655.09 & 0 & 2 & $y_{\K+1} = 1$, $y_{\K+2} + y_{\K+3}\geq 1$   \\
 &  & 3 & 7 & TL & 930.99 & 0 & 5 & TL & 927.70 & 0 & 4 & $y_{\K+1} = 1$, $y_{\K+2} = 1$   \\*[0.18cm]
 & 0.4 & 1 & 1 & 0.19 & 0 & 0 & 0 & \textBF{0.16} & 0 & 0 & 0 \\
 &  & 2 & 2 & 59.29 & 56.01 & 0 & 1 & \textBF{0.50} & 0 & 0 & 0 & $y_{\K+1} + y_{\K+2}  \geq 1 $ \\
 &  & 3 & 1 & 0.11 & 0 & 0 & 0 & \textBF{0.09} & 0 & 0 & 0 \\*[0.18cm]
 & 0.5 & 1 & 1 & 0.08 & 0 & 0 & 0 & \textBF{0.07} & 0 & 0 & 0 \\
 &  & 2 & 1 & 0.07 & 0 & 0 & 0 & \textBF{0.07} & 0 & 0 & 0 \\
 &  & 3 & 1 & \textBF{0.06} & 0 & 0 & 0 & 0.06 & 0 & 0 & 0 \\
 \midrule
30 & 0.3 & 1 & 3 & 209.49 & 206.41 & 0 & 2 & \textBF{0.86} & 0 & 0 & 0 & $y_{\K+1} = 1$, $y_{\K+2} + y_{\K+3}\geq 1$   \\
 &  & 2 & 4 & 267.44 & 260.88 & 0 & 3 & \textBF{204.33} & 187.15 & 0 & 1 & $y_{\K+1} = 1$, $y_{\K+2} = 1$   \\
 &  & 3 & 4 &TL & 967.79 & 0 & 5 & TL & 937.77 & 0 & 2 & $y_{\K+1} + y_{\K+2}  \geq 1 $ \\*[0.18cm]
 & 0.4 & 1 & 1 & \textBF{0.45} & 0 & 0 & 0 & 0.57 & 0 & 0 & 0 \\
 &  & 2 & 1 & \textBF{0.30} & 0 & 0 & 0 & 0.40 & 0 & 0 & 0 \\
 &  & 3 & 1 & \textBF{0.12} & 0 & 0 & 0 & 0.14 & 0 & 0 & 0 \\*[0.18cm]
 & 0.5 & 1 & 1 & \textBF{0.12} & 0 & 0 & 0 & 0.16 & 0 & 0 & 0 \\
 &  & 2 & 1 & \textBF{0.11} & 0 & 0 & 0 & 0.15 & 0 & 0 & 0 \\
 &  & 3 & 1 & \textBF{0.12} & 0 & 0 & 0 & 0.14 & 0 & 0 & 0 \\

 \midrule
35 & 0.3 & 1 & 3 & TL & 910.55 & 0 & 2 & TL & 944.54 & 0 & 1 & $y_{\K+1} + y_{\K+2}  \geq 1 $ \\
 &  & 2 & 3 & 637.82 & 596.94 & 0 & 2 & \textBF{0.61} & 0 & 0 & 0 & $y_{\K+1} = 1$, $y_{\K+2} + y_{\K+3}\geq 1$   \\
 &  & 3 & 3 & 519.48 & 469.57 & 0 & 2 & \textBF{0.86} & 0 & 0 & 0 & $y_{\K+1} = 1$, $y_{\K+2} = 1$   \\*[0.18cm]
 & 0.4 & 1 & 1 & \textBF{0.46} & 0 & 0 & 0 & 0.52 & 0 & 0 & 0 \\
 &  & 2 & 1 & \textBF{0.20} & 0 & 0 & 0 & 0.24 & 0 & 0 & 0 \\
 &  & 3 & 1 & \textBF{0.28} & 0 & 0 & 0 & 0.36 & 0 & 0 & 0 \\*[0.18cm]
 & 0.5 & 1 & 1 & \textBF{0.20} & 0 & 0 & 0 & 0.26 & 0 & 0 & 0 \\
 &  & 2 & 1 & \textBF{0.20} & 0 & 0 & 0 & 0.26 & 0 & 0 & 0 \\
 &  & 3 & 1 & \textBF{0.20} & 0 & 0 & 0 & 0.25 & 0 & 0 & 0 \\
\bottomrule
\end{longtable}
\end{landscape}
}

\begin{landscape}
\renewcommand\arraystretch{0.8}
\vspace*{-1.5cm}
\centering
	\small
	\begin{longtable}{cccr rr rr rr  rrr}
		\caption{Results for the witness-based decomposition on random instances, the column headers define the options used to solve each instance. We include only the times for those methods which are clearly outperformed.}
		\label{table:MinDWitness} \\
		\toprule
  &  &  &  & \multicolumn{1}{c}{IIS} &  \multicolumn{1}{c}{IIS} &  \multicolumn{4}{c}{IIS} &\multicolumn{3}{c}{cycle separation} \\
  &  &  &  & \multicolumn{1}{c}{no cycle breaking } &  \multicolumn{1}{c}{ 2-cycle breaking} &  \multicolumn{4}{c}{2-cycle \& 3-cycle breaking} &\multicolumn{3}{c}{2-cycle \& 3-cycle breaking} \\
 \cmidrule(lr){5-5}
 \cmidrule(lr){6-6}
 \cmidrule(lr){7-10}
 \cmidrule(lr){11-13}
$n$ & $\density$ & Inst. & Obj. & Time & Time & Time &  IIS time &  BB Nodes &  \# Cuts & Time &  BB Nodes &  \# Cuts \\
\midrule
20 & 0.3 & 1 & 14 & 0.82 & 0.22 & 0.11 & 0.04 & 29 & 2 & \textBF{0.05} & 30 & 1 \\
 &  & 2 & $+\infty$ & 0.00 & 0.00 & \textBF{0.00} & 0.00 & 0 & 0 & 0.00 & 0 & 0 \\
 &  & 3 & $+\infty$ & 1.27 & 1.11 & 0.49 & 0.39 & 393 & 18 & \textBF{0.09} & 250 & 30 \\*[0.18cm]
 & 0.4 & 1 & 2 & 503.49 & 322.40 & 188.39 & 23.55 & 241864 & 764 & \textBF{72.06} & 136729 & 1086 \\
 &  & 2 & 3 & TL & TL & 826.28 & 44.20 & 868805 & 1309 & \textBF{765.68} & $>$1M & 1858 \\
 &  & 3 & 4 & TL & TL & TL & 56.34 & 819500 & 1557 & TL & 850346 & 2865 \\*[0.18cm]
 & 0.5 & 1 & 1 & TL & 19.35 & 0.02 & 0.00 & 0 & 0 & \textBF{0.02} & 0 & 0 \\
 &  & 2 & 1 & 0.86 & 2.51 & \textBF{0.46} &0.38 & 9 & 12 & 0.51 & 534 & 290 \\
 &  & 3 & 1 & 67.42 & TL & 899.36 & 167.57 & 483400 & 4726 & \textBF{0.30} & 250 & 191 \\
 \midrule
25 & 0.3 & 1 & $+\infty$ & 0.00 & \textBF{0.00} & 0.00 & 0.00 & 0 & 0 & 0.00 & 0 & 0 \\
 &  & 2 & 5 & 265.95 & 182.43 & \textBF{102.67} & 15.28 & 147275 & 448 & 128.44 & 214519 & 428 \\
 &  & 3 & 7 & TL & TL & TL & 17.55 & $>$1.2M & 653 & TL & 814100 & 1561 \\*[0.18cm]
 & 0.4 & 1 & 1 & 49.47 & 22.69 & 1.32 & 1.18 & 19 & 28 & \textBF{0.05} & 13 & 12 \\
 &  & 2 & 2 & TL & TL & TL & 390.42 & 76940 & 10783 & TL & 435200 & 13778 \\
 &  & 3 & 1 & TL & 822.20 & 74.26 & 63.96 & 12997 & 1917 & \textBF{2.25} & 3508 & 1213 \\*[0.18cm]
 & 0.5 & 1 & 1 & 14.83 & 5.61 & 33.98 & 31.50 & 1285 & 997 & \textBF{1.53} & 1128 & 900 \\
 &  & 2 & 1 & 144.43 & \textBF{0.02} & 0.10 & 0.05 & 0 & 2 & 0.05 & 0 & 2 \\
 &  & 3 & 1 & 1.69 & 74.92 & 7.33 & 6.78 & 140 & 231 & \textBF{0.13} & 0 & 3 \\
 \midrule
30 & 0.3 & 1 & 3 & TL & TL & TL & 395.67 & 96929 & 8766 & TL & 110687 & 16086 \\
 &  & 2 & 4 & TL & TL & TL & 227.83 & 164700 & 5215 & TL & 120326 & 14070 \\
 &  & 3 & 4 & TL & TL & TL & 129.15 & 412610 & 3121 & TL & 271708 & 6599 \\*[0.18cm]
 & 0.4 & 1 & 1 & TL & 2.44 & TL & 431.60 & 60001 & 10376 & \textBF{1.16} & 983 & 874 \\
 &  & 2 & 1 & TL & TL & 111.46 & 101.08 & 3962 & 2389 & \textBF{2.16} & 1076 & 1398 \\
 &  & 3 & 1 & 47.21 & 7.55 & 0.48 & 0.34 & 3 & 10 & \textBF{0.30} & 162 & 226 \\*[0.18cm]
 & 0.5 & 1 & 1 & 11.14 & 1.97 & TL & 784.87 & 37527 & 16374 & \textBF{0.25} & 30 & 57 \\
 &  & 2 & 1 & 4.93 & 2.63 & 0.88 & 0.65 & 7 & 19 & \textBF{0.27} & 19 & 41 \\
 &  & 3 & 1 & \textBF{0.19} & 40.20 & 5.71 & 5.23 & 105 & 130 & 0.27 & 63 & 69 \\
 \midrule
35 & 0.3 & 1 & 3 & TL & TL & TL & 389.40 & 81254 & 8553 & TL & 200157 & 10093 \\
 &  & 2 & 3 & TL & TL & TL & 330.50 & 102679 & 6448 & TL & 144100 & 12350 \\
 &  & 3 & 3 & TL & TL & TL & 384.89 & 102852 & 8311 & TL & 106642 & 23345 \\*[0.18cm]
 & 0.4 & 1 & 1 & TL & TL & TL & 670.87 & 48007 & 13694 & TL & 104470 & 26293 \\
 &  & 2 & 1 & TL & TL & 1.24 & 0.95 & 14 & 21 & \textBF{0.39} & 173 & 227 \\
 &  & 3 & 1 & TL & 25.52 & TL & 734.19 & 56599 & 13849 & \textBF{0.61} & 230 & 339 \\*[0.18cm]
 & 0.5 & 1 & 1 & 1.10 & 19.22 & 6.28 & 5.64 & 116 & 133 & \textBF{0.16} & 0 & 3 \\
 &  & 2 & 1 & 11.84 & 20.52 & 159.28 & 144.22 & 3981 & 2899 & \textBF{0.89} & 233 & 389 \\
 &  & 3 & 1 & 2.51 & 25.54 & 3.20 & 2.74 & 41 & 61 & \textBF{1.36} & 389 & 671 \\
\bottomrule
\end{longtable}
\end{landscape}
\begin{landscape}
\renewcommand\arraystretch{0.8}
\vspace*{-1.28cm}
\centering
	\small
	\begin{longtable}{cccr rr rr rr rrr rrr}
		\caption{Results for the best performing formulations on random instances.}
		\label{table:MinDBest} \\
		\toprule
 &  &  &  & \multicolumn{2}{c}{$\IP$} & \multicolumn{2}{c}{$\vertexCP$}&\multicolumn{2}{c}{$\combCP$}  &   \multicolumn{3}{c}{Naive Decomp. (with VI)} & \multicolumn{3}{c}{Witness-based Decomp.} \\
  &  &  &  & \multicolumn{2}{c}{} & \multicolumn{2}{c}{ (with VI)}&\multicolumn{2}{c}{ (with VI)}  &   \multicolumn{3}{c}{} & \multicolumn{3}{c}{cycle separation} \\
    &  &  &  & \multicolumn{2}{c}{} & \multicolumn{2}{c}{}&\multicolumn{2}{c}{}  &   \multicolumn{3}{c}{} & \multicolumn{3}{c}{2-cycle \& 3-cycle breaking} \\
\cmidrule(lr){5-6}
\cmidrule(lr){7-8}
\cmidrule(lr){9-10}
\cmidrule(lr){11-13}
\cmidrule(lr){14-16}
$n$ & $\density$ & Inst. &  Obj. &  Time &  BB Nodes & Time & Ch.Pts. & Time & Ch.Pts. & Time &  BB Nodes &  \#Cuts & Time &  BB Nodes &  \#Cuts \\
\midrule
20 & 0.3 & 1 & 14 & 956.68 & 453417 & 932.80 & $>$17M & TL & $>$14M & TL & 0 & 6 & \textBF{0.05} & 30 & 1 \\
 &  & 2 & $+\infty$ & 0.06 & 0 & 1.39 & 28232 & \textBF{0.00} & 0 & 597.65 & 1 & 106 & 0.00 & 0 & 0 \\
 &  & 3 & $+\infty$ & 43.16 & 15086 & 34.94 & 739690 & 8.27 & 209521 & TL & 0 & 22 & \textBF{0.09} & 250 & 30 \\*[0.18cm]
 & 0.4 & 1 & 2 & 5.74 & 3028 & 0.03 & 451 & \textBF{0.02} & 455 & 0.23 & 0 & 0 & 72.06 & 136729 & 1086 \\
 &  & 2 & 3 & 25.53 & 10890 & 0.03 & 436 & \textBF{0.02} & 446 & 0.03 & 0 & 0 & 765.68 & $>$1M & 1858 \\
 &  & 3 & 4 & 30.28 & 12493 & \textBF{7.12} & 146090 & 9.10 & 201655 & 48.10 & 0 & 1 & TL & 850346 & 2865 \\*[0.18cm]
 & 0.5 & 1 & 1 & 0.09 & 0 & \textBF{0.02} & 186 & 0.02 & 453 & 0.08 & 0 & 0 & 0.02 & 0 & 0 \\
 &  & 2 & 1 & 0.09 & 0 & \textBF{0.02} & 160 & 0.03 & 463 & 0.03 & 0 & 0 & 0.51 & 534 & 290 \\
 &  & 3 & 1 & 0.09 & 0 & \textBF{0.02} & 228 & 0.02 & 445 & 0.03 & 0 & 0 & 0.30 & 250 & 191 \\
 \midrule
25 & 0.3 & 1 & $+\infty$ & 0.11 & 0 & TL & $>$13M & 0.01 & 12 & TL & 0 & 1 & \textBF{0.00} & 0 & 0 \\
 &  & 2 & 5 & \textBF{80.57} & 11650 & TL & $>$12M & TL & $>$12M & TL & 0 & 2 & 128.44 & 214519 & 428 \\
 &  & 3 & 7 & \textBF{581.19} & 95937 & TL & $>$13M & TL & $>$12M & TL & 0 & 4 & TL & 814100 & 1561 \\*[0.18cm]
 & 0.4 & 1 & 1 & 0.17 & 0 & \textBF{0.05} & 447 & 0.05 & 470 & 0.16 & 0 & 0 & 0.05 & 13 & 12 \\
 &  & 2 & 2 & 35.83 & 8983 & 0.05 & 465 & \textBF{0.03} & 458 & 0.50 & 0 & 0 & TL & 435200 & 13778 \\
 &  & 3 & 1 & 0.18 & 0 & 0.04 & 451 & \textBF{0.04} & 450 & 0.09 & 0 & 0 & 2.25 & 3508 & 1213 \\*[0.18cm]
 & 0.5 & 1 & 1 & 0.21 & 0 & \textBF{0.04} & 467 & 0.05 & 471 & 0.07 & 0 & 0 & 1.53 & 1128 & 900 \\
 &  & 2 & 1 & 0.21 & 0 & 0.04 & 459 & \textBF{0.04} & 469 & 0.07 & 0 & 0 & 0.05 & 0 & 2 \\
 &  & 3 & 1 & 0.21 & 0 & \textBF{0.04} & 462 & 0.04 & 460 & 0.06 & 0 & 0 & 0.13 & 0 & 3 \\
 \midrule
30 & 0.3 & 1 & 3 & 160.00 & 20239 & 0.07 & 485 & \textBF{0.06} & 469 & 0.86 & 0 & 0 & TL & 110687 & 16086 \\
 &  & 2 & 4 & 164.67 & 20843 & 15.80 & 138956 & \textBF{12.84} & 164305 & 204.33 & 0 & 1 & TL & 120326 & 14070 \\
 &  & 3 & 4 & TL & 74413 & \textBF{430.65} & $>$4.5M & 588.08 & $>$6.9M & TL & 0 & 2 & TL & 271708 & 6599 \\*[0.18cm]
 & 0.4 & 1 & 1 & 0.49 & 0 & 0.08 & 478 & \textBF{0.07} & 478 & 0.57 & 0 & 0 & 1.16 & 983 & 874 \\
 &  & 2 & 1 & 0.45 & 0 & 0.08 & 468 & \textBF{0.06} & 458 & 0.40 & 0 & 0 & 2.16 & 1076 & 1398 \\
 &  & 3 & 1 & 0.51 & 0 & 0.07 & 467 & \textBF{0.07} & 476 & 0.14 & 0 & 0 & 0.30 & 162 & 226 \\*[0.18cm]
 & 0.5 & 1 & 1 & 0.65 & 0 & \textBF{0.05} & 169 & 0.08 & 475 & 0.16 & 0 & 0 & 0.25 & 30 & 57 \\
 &  & 2 & 1 & 0.57 & 0 & 0.05 & 207 & \textBF{0.04} & 181 & 0.15 & 0 & 0 & 0.27 & 19 & 41 \\
 &  & 3 & 1 & 0.59 & 0 & \textBF{0.05} & 203 & 0.07 & 471 & 0.14 & 0 & 0 & 0.27 & 63 & 69 \\
 \midrule
35 & 0.3 & 1 & 3 & TL & 70945 & 258.45 & $>$1.6M & \textBF{209.41} & $>$2.1M & TL & 0 & 1 & TL & 200157 & 10093 \\
 &  & 2 & 3 & TL & 52051 & \textBF{0.12} & 492 & 0.12 & 503 & 0.61 & 0 & 0 & TL & 144100 & 12350 \\
 &  & 3 & 3 & 655.51 & 46969 & 0.13 & 480 & \textBF{0.12} & 490 & 0.86 & 0 & 0 & TL & 106642 & 23345 \\*[0.18cm]
 & 0.4 & 1 & 1 & 1.00 & 0 & 0.13 & 487 & \textBF{0.12} & 489 & 0.52 & 0 & 0 & TL & 104470 & 26293 \\
 &  & 2 & 1 & 0.93 & 0 & 0.13 & 477 & \textBF{0.11} & 492 & 0.24 & 0 & 0 & 0.39 & 173 & 227 \\
 &  & 3 & 1 & 0.92 & 0 & \textBF{0.11} & 208 & 0.12 & 499 & 0.36 & 0 & 0 & 0.61 & 230 & 339 \\*[0.18cm]
 & 0.5 & 1 & 1 & 1.16 & 0 & 0.10 & 205 & \textBF{0.09} & 171 & 0.26 & 0 & 0 & 0.16 & 0 & 3 \\
 &  & 2 & 1 & 1.13 & 0 & \textBF{0.07} & 88 & 0.09 & 168 & 0.26 & 0 & 0 & 0.89 & 233 & 389 \\
 &  & 3 & 1 & 1.12 & 0 & \textBF{0.10} & 170 & 0.13 & 501 & 0.25 & 0 & 0 & 1.36 & 389 & 671 \\
\bottomrule
\end{longtable}
\end{landscape}

\afterpage{
\begin{landscape}
	\renewcommand\arraystretch{0.8}
	\begin{center}
		\small
		\begin{longtable}{cccr rrr rr rr rrr rrr}
			\caption{\cRev Results for the best performing formulations on synthetic instances.}
			\label{table:Syn} \\
			\toprule
			 &  &  &  & \multicolumn{3}{c}{$\CCG$} & \multicolumn{2}{c}{$\vertexCP$} & \multicolumn{2}{c}{$\combCP$} & \multicolumn{3}{c}{Naive Decomp. (with VI)} & \multicolumn{3}{c}{Witness-based Decomp.}\\
			 			 &  &  &  & \multicolumn{3}{c}{} & \multicolumn{2}{c}{(with VI)} & \multicolumn{2}{c}{(with VI)} & \multicolumn{3}{c}{} & \multicolumn{3}{c}{cycle separation}\\
			 &  &  &  & \multicolumn{3}{c}{} & \multicolumn{2}{c}{} & \multicolumn{2}{c}{} & \multicolumn{3}{c}{} & \multicolumn{3}{c}{2-cycle \& 3-cycle breaking}\\			 
			 \cmidrule(lr){5-7}
			 \cmidrule(lr){8-9}
			 \cmidrule(lr){10-11}
			 \cmidrule(lr){12-14}
			 \cmidrule(lr){15-17}
$n$ & Doubles & Noise & Obj. & Time &   \multicolumn{1}{c}{BB}  &  \# Cuts &  Time & Ch.Pts.&  Total Time & Ch.Pts. &  Time &  \multicolumn{1}{c}{BB}     &  \# Cuts &  Time &   \multicolumn{1}{c}{BB}    &  \# Cuts \\
&  &  &  &  & \multicolumn{1}{c}{Nodes}   &   &  & &   &  &   &   \multicolumn{1}{c}{Nodes}   &   &   &    \multicolumn{1}{c}{Nodes}   &   \\
\midrule
25 & 3 & 3 & 2 & 82.87 & 132209 & 292 & \textBF{0.25} & 5989 & 1.57 & 27179 & 0.47 & 0 & 0 & 1.01 & 484 & 20 \\
&  & 4 & 2 & 781.13 & $>$1M & 560 & \textBF{0.29} & 6676 & 0.45 & 9189 & 0.34 & 0 & 0 & 2.08 & 1298 & 40 \\
&  & 5 & 2 & TL & $>$1M & 977 & \textBF{0.19} & 4492 & 0.66 & 17568 & 0.40 & 0 & 0 & 7.42 & 3750 & 131 \\*[0.18cm]
& 4 & 3 & 2 & 5.28 & 7007 & 180 & 45.37 & 966391 & 55.34 & $>$1M & TL & 0 & 1 & \textBF{1.09} & 257 & 22 \\
&  & 4 & 2 & 96.35 & 125717 & 423 & 46.76 & 973021 & 102.71 & $>$1.7M & TL & 0 & 1 & \textBF{0.59} & 344 & 7 \\
&  & 5 & 1 & 111.11 & 152250 & 250 & 0.61 & 12744 & \textBF{0.60} & 16244 & 0.62 & 0 & 0 & 1.92 & 667 & 38 \\*[0.18cm]
& 5 & 3 & 3 & 5.26 & 6841 & 70 & TL & $>$20M & TL & $>$14M & TL & 0 & 0 & \textBF{0.94} & 555 & 19 \\
&  & 4 & 3 & 202.42 & 282382 & 298 & TL & $>$21M & TL & $>$13M & TL & 0 & 0 & \textBF{9.17} & 5014 & 170 \\
&  & 5 & 3 & TL & $>$1M & 479 & TL & $>$19M & TL & $>$18M & TL & 0 & 0 & \textBF{108.01} & 34378 & 1494 \\
\midrule
30 & 3 & 3 & 2 & 293.07 & 226800 & 1014 & \textBF{0.68} & 10973 & 1.64 & 26819 & 0.91 & 0 & 0 & 2.64 & 1372 & 31 \\
&  & 5 & 2 & TL & 703581 & 1534 & 0.97 & 13645 & 3.36 & 56883 & \textBF{0.84} & 0 & 0 & 29.56 & 25076 & 174 \\
&  & 6 & 2 & TL & 727815 & 974 & \textBF{1.02} & 15450 & 1.07 & 21795 & 1.07 & 0 & 0 & 96.23 & 43878 & 1012 \\*[0.18cm]
& 4 & 3 & 3 & 166.98 & 192681 & 314 & TL & $>$12M & TL & $>$9.5M & TL & 0 & 0 & \textBF{11.52} & 1873 & 235 \\
&  & 5 & 3 & TL & $>$1M & 541 & TL & $>$12M & TL & $>$9.4M & TL & 0 & 0 & \textBF{21.73} & 13673 & 234 \\
&  & 6 & 3 & TL & $>$1M & 641 & TL & $>$12M & TL & $>$13M & TL & 0 & 0 & \textBF{49.62} & 45561 & 376 \\*[0.18cm]
& 5 & 3 & 4 & 323.10 & 382120 & 261 & TL & $>$14M & TL & $>$12M & TL & 0 & 0 & \textBF{9.71} & 4323 & 157 \\
&  & 5 & 4 & TL & $>$1M & 579 & TL & $>$14M & TL & $>$12M & TL & 0 & 0 & \textBF{59.32} & 70357 & 263 \\
&  & 6 & 3 & TL & 982200 & 569 & TL & $>$13M & TL & $>$14M & TL & 0 & 0 & \textBF{23.30} & 19367 & 187 \\
\midrule
35 & 4 & 4 & 4 & TL & 919450 & 441 & TL & $>$10M & TL & $>$9.5M & TL & 0 & 0 & \textBF{17.53} & 12017 & 121 \\
&  & 6 & 4 & TL & 511962 & 1028 & TL & $>$10M & TL & $>$9.2M & TL & 0 & 0 & \textBF{428.53} & 343809 & 286 \\
&  & 7 & 3 & TL & 459000 & 1277 & TL & $>$10M & TL & $>$13M & TL & 0 & 0 & TL & 773333 & 477 \\*[0.18cm]
& 5 & 4 & 4 & TL & 756700 & 402 & TL & $>$9.9M & TL & $>$9M & TL & 0 & 0 & \textBF{11.77} & 9761 & 48 \\
&  & 6 & 4 & TL & 741600 & 895 & TL & $>$9.8M & TL & $>$9.2M & TL & 0 & 0 & \textBF{134.45} & 117720 & 247 \\
&  & 7 & 4 & TL & 479861 & 1098 & TL & $>$10M & TL & $>$11M & TL & 0 & 0 & \textBF{836.95} & 767184 & 387 \\*[0.18cm]
& 6 & 4 & 4 & 109.11 & 102481 & 298 & TL & $>$9.9M & TL & $>$10M & TL & 0 & 0 & \textBF{2.88} & 1348 & 17 \\
&  & 6 & 4 & TL & 794092 & 599 & TL & $>$10M & TL & $>$9.7M & TL & 0 & 0 & \textBF{42.49} & 31154 & 144 \\
&  & 7 & 3 & TL & 482000 & 1333 & TL & $>$10M & TL & $>$13M & TL & 0 & 0 & \textBF{35.25} & 33400 & 172 \\
\bottomrule
\end{longtable}
\end{center}
\end{landscape}
\clearpage 
\restoregeometry}

\begin{landscape}
\renewcommand\arraystretch{0.8}
\vspace*{-1.28cm}
\centering
	\small
	\begin{longtable}{cccr rr rr rrr rr rrr rrr}
		\caption{\cRev Results for the best performing formulations on random instances with $\K=4$.}
		\label{table:K4} \\
		\toprule
 &  &  &  & \multicolumn{2}{c}{$\IP$} & \multicolumn{3}{c}{$\CCG$} & \multicolumn{2}{c}{$\vertexCP$}&\multicolumn{2}{c}{$\combCP$}  &   \multicolumn{3}{c}{Naive Decomp.} & \multicolumn{3}{c}{Witness-based Decomp.} \\
  &  &  &  & \multicolumn{2}{c}{} &  \multicolumn{3}{c}{} & \multicolumn{2}{c}{ (with VI)}&\multicolumn{2}{c}{ (with VI)}  &   \multicolumn{3}{c}{ (with VI)} & \multicolumn{3}{c}{cycle separation} \\
    &  &  &  & \multicolumn{2}{c}{}&   \multicolumn{3}{c}{}  & \multicolumn{2}{c}{}&\multicolumn{2}{c}{}  &   \multicolumn{3}{c}{} & \multicolumn{3}{c}{2-cycle \& 3-cycle breaking} \\
\cmidrule(lr){5-6}
\cmidrule(lr){7-9}
\cmidrule(lr){10-11}
\cmidrule(lr){12-13}
\cmidrule(lr){14-16}
\cmidrule(lr){17-19}
$n$ & $\density$ & Inst. & Obj. &  Time & \multicolumn{1}{l}{BB}  &Time &   \multicolumn{1}{l}{BB}    &  \#Cuts & Time & Ch.Pts. & Time & Ch.Pts. & Time &   \multicolumn{1}{l}{BB}    &  \#Cuts & Time &   \multicolumn{1}{l}{BB}    &  \#Cuts \\
 & &  & &   &    \multicolumn{1}{l}{Nodes} & &   \multicolumn{1}{l}{Nodes} &  &  & &  &  &  &   \multicolumn{1}{l}{Nodes} &   &  &   \multicolumn{1}{l}{Nodes} &  \\
\midrule
20 & 0.5 & 1 & 3 & 198.54 & 83515 & TL & 323383 & 4152 & 97.08 & $>$2.2M & \textBF{46.34} & $>$1.3M & TL & 0 & 4 & TL & 698737 & 4096 \\
&  & 2 & 2 & 24.41 & 12257 & TL & 429900 & 3618 & \textBF{0.13} & 4468 & 0.26 & 9707 & 28.15 & 0 & 2 & TL & $>$1M & 3622 \\
&  & 3 & 4 & 294.51 & 99787 & TL & 443716 & 3565 & 123.19 & $>$2.8M & \textBF{82.99} & $>$2.1M & TL & 0 & 4 & TL & 610254 & 4891 \\*[0.18cm]
& 0.7 & 1 & 2 & 0.38 & 0 & TL & 87819 & 18949 & \textBF{0.01} & 28 & 0.01 & 36 & 0.04 & 0 & 0 & 0.06 & 1 & 3 \\
&  & 2 & 2 & 0.50 & 100 & TL & 90067 & 17007 & \textBF{0.01} & 102 & 0.01 & 104 & 0.02 & 0 & 0 & 32.17 & 23946 & 8097 \\
&  & 3 & 2 & 0.57 & 20 & TL & 90497 & 18369 & \textBF{0.01} & 32 & 0.01 & 30 & 0.02 & 0 & 0 & 0.11 & 69 & 58 \\
\midrule
25 & 0.4 & 1 & 2 & 100.53 & 23942 & TL & 255230 & 5083 & 3.23 & 60795 & \textBF{0.65} & 18463 & 222.89 & 0 & 2 & TL & 944585 & 1953 \\
&  & 3 & 4 & 191.34 & 39444 & TL & 183556 & 7808 & \textBF{37.57} & 778990 & 115.23 & $>$2.5M & 659.51 & 0 & 4 & TL & 754647 & 1544 \\*[0.18cm]
& 0.5 & 1 & 2 & 9.37 & 2171 & TL & 103080 & 17686 & 0.22 & 5032 & \textBF{0.10} & 2482 & 0.11 & 0 & 0 & TL & 452598 & 10506 \\
&  & 2 & 2 & 3.06 & 544 & TL & 98861 & 17595 & \textBF{0.04} & 518 & 0.15 & 3844 & 0.11 & 0 & 0 & 2.47 & 2679 & 1290 \\
&  & 3 & 2 & 117.14 & 27346 & TL & 103104 & 15080 & \textBF{0.07} & 1565 & 0.34 & 9083 & 43.63 &  &  & TL & 335431 & 18925 \\*[0.18cm]
& 0.7 & 1 & 2 & 0.71 & 0 & TL & 43425 & 34277 & \textBF{0.01} & 31 & 0.02 & 57 & 0.05 & 0 & 0 & 0.09 & 4 & 0 \\
&  & 2 & 2 & 1.41 & 30 & TL & 46617 & 35453 & \textBF{0.01} & 31 & 0.02 & 52 & 0.05 & 0 & 0 & 0.19 & 14 & 23 \\
&  & 3 & 2 & 1.09 & 0 & TL & 50974 & 36328 & \textBF{0.01} & 33 & 0.02 & 48 & 0.05 & 0 & 0 & 0.03 & 0 & 0 \\
\midrule
30 & 0.4 & 1 & 2 & TL & 117419 & TL & 111696 & 14955 & \textBF{3.19} & 48203 & 6.74 & 119001 & TL & 0 & 2 & TL & 259546 & 9562 \\
&  & 2 & 3 & 840.23 & 91914 & TL & 110217 & 15295 & 249.28 & $>$4M & \textBF{56.69} & 994726 & TL & 0 & 2 & TL & 128032 & 16278 \\
&  & 3 & 3 & TL & 99842 & TL & 96684 & 13686 & 252.12 & $>$3.5M & \textBF{77.23} & $>$1.3M & TL & 0 & 1 & TL & 193379 & 12298 \\*[0.18cm]
& 0.5 & 1 & 2 & 25.32 & 2743 & TL & 45286 & 32276 & \textBF{0.04} & 71 & 0.09 & 1323 & 0.52 & 0 & 0 & 27.12 & 13384 & 5019 \\
&  & 2 & 2 & 11.98 & 604 & TL & 47654 & 32034 & \textBF{0.05} & 626 & 0.47 & 9244 & 0.81 & 0 & 0 & TL & 115100 & 24688 \\
&  & 3 & 2 & 800.83 & 103399 & TL & 45780 & 28561 & 0.12 & 1596 & \textBF{0.08} & 1396 & TL & 0 & 1 & TL & 172140 & 25786 \\*[0.18cm]
& 0.7 & 1 & 2 & 1.74 & 0 & TL & 36983 & 46019 & \textBF{0.03} & 38 & 0.04 & 284 & 0.13 & 0 & 0 & 0.11 & 0 & 1 \\
&  & 2 & 2 & 1.55 & 0 & TL & 34086 & 45394 & \textBF{0.03} & 37 & 0.03 & 38 & 0.15 & 0 & 0 & 0.19 & 4 & 6 \\
&  & 3 & 2 & 1.18 & 0 & TL & 31620 & 49277 & \textBF{0.03} & 38 & 0.03 & 37 & 0.11 & 0 & 0 & 0.14 & 2 & 5 \\
\midrule
35 & 0.4 & 2 & 2 & 40.95 & 1573 & TL & 56761 & 29474 & 0.49 & 7500 & 1.33 & 18297 & \textBF{0.27} & 0 & 0 & 0.65 & 287 & 308 \\
&  & 3 & 2 & 811.71 & 56844 & TL & 52878 & 28871 & 0.43 & 4930 & \textBF{0.07} & 125 & TL & 0 & 1 & TL & 142021 & 21304 \\*[0.18cm]
& 0.5 & 1 & 2 & 15.03 & 731 & TL & 35316 & 41620 & \textBF{0.10} & 562 & 0.97 & 14074 & 0.23 & 0 & 0 & 0.37 & 22 & 35 \\
&  & 2 & 2 & 12.19 & 164 & TL & 34960 & 45237 & \textBF{0.06} & 115 & 0.40 & 5611 & 1.44 & 0 & 0 & 8.25 & 3354 & 3285 \\
&  & 3 & 2 & 10.40 & 150 & TL & 34937 & 44867 & \textBF{0.08} & 199 & 0.13 & 1308 & 0.80 & 0 & 0 & TL & 214703 & 22183 \\*[0.18cm]
& 0.7 & 1 & 2 & 5.30 & 40 & TL & 20362 & 47286 & \textBF{0.06} & 55 & 0.06 & 62 & 0.20 & 0 & 0 & 0.15 & 0 & 0 \\
&  & 2 & 2 & 2.87 & 0 & TL & 19120 & 47919 & \textBF{0.06} & 57 & 0.06 & 59 & 0.21 & 0 & 0 & 0.24 & 43 & 54 \\
&  & 3 & 2 & 4.17 & 0 & TL & 18270 & 47295 & \textBF{0.06} & 55 & 0.06 & 60 & 0.22 & 0 & 0 & 0.11 & 3 & 3 \\
\bottomrule
\end{longtable}
\end{landscape}

\begin{landscape}
\renewcommand\arraystretch{0.8}
\vspace*{-1.28cm}
\centering
	\small
\begin{longtable}{cccr rr rrr rr rr rrr rrr}
		\caption{\cRev Results for the best performing formulations on random instances with $\K=5$.}
		\label{table:K5} \\
		\toprule
  &  &  &  & \multicolumn{2}{c}{$\IP$} & \multicolumn{3}{c}{$\CCG$} & \multicolumn{2}{c}{$\vertexCP$}&\multicolumn{2}{c}{$\combCP$}  &   \multicolumn{3}{c}{Naive Decomp.} & \multicolumn{3}{c}{Witness-based Decomp.} \\
 &  &  &  & \multicolumn{2}{c}{} &  \multicolumn{3}{c}{} & \multicolumn{2}{c}{ (with VI)}&\multicolumn{2}{c}{ (with VI)}  &   \multicolumn{3}{c}{ (with VI)} & \multicolumn{3}{c}{cycle separation} \\
 &  &  &  & \multicolumn{2}{c}{}&   \multicolumn{3}{c}{}  & \multicolumn{2}{c}{}&\multicolumn{2}{c}{}  &   \multicolumn{3}{c}{} & \multicolumn{3}{c}{2-cycle \& 3-cycle breaking} \\
 \cmidrule(lr){5-6}
 \cmidrule(lr){7-9}
 \cmidrule(lr){10-11}
 \cmidrule(lr){12-13}
 \cmidrule(lr){14-16}
 \cmidrule(lr){17-19}
$n$ & $\density$ & Inst. & Obj. &  Time & \multicolumn{1}{l}{BB}  &Time &   \multicolumn{1}{l}{BB}    &  \#Cuts & Time & Ch.Pts. & Time & Ch.Pts. & Time &   \multicolumn{1}{l}{BB}    &  \#Cuts & Time &   \multicolumn{1}{l}{BB}    &  \#Cuts \\
& &  & &   &    \multicolumn{1}{l}{Nodes} & &   \multicolumn{1}{l}{Nodes} &  &  & &  &  &  &   \multicolumn{1}{l}{Nodes} &   &  &   \multicolumn{1}{l}{Nodes} &  \\
\midrule
20 & 0.7 & 1 & 2 & 0.52 & 0 & TL & 127066 & 13270 & \textBF{0.01} & 28 & 0.02 & 541 & 0.03 & 0 & 0 & 0.23 & 44 & 71 \\
&  & 2 & 2 & 209.33 & 117075 & TL & 120642 & 13778 & 0.05 & 1211 & \textBF{0.05} & 1609 & 969.40 & 0 & 1 & TL & 310055 & 13348 \\
&  & 3 & 2 & 0.86 & 260 & TL & 171232 & 10873 & \textBF{0.01} & 28 & 0.02 & 484 & 0.03 & 0 & 0 & 0.13 & 3 & 2 \\
\midrule
25 & 0.5 & 2 & 4 & TL & 183338 & TL & 133379 & 10894 & \textBF{314.66} & $>$4.5M & 559.60 & $>$9M & TL & 0 & 3 & TL & 329012 & 5492 \\*[0.18cm]
& 0.7 & 1 & 2 & 2.30 & 0 & TL & 58130 & 23335 & 0.02 & 147 & \textBF{0.02} & 56 & TL & 0 & 3 & 0.31 & 81 & 103 \\
&  & 2 & 2 & 1.65 & 80 & TL & 57317 & 28280 & \textBF{0.02} & 158 & 0.07 & 1536 & 0.05 & 0 & 0 & 0.26 & 4 & 12 \\
&  & 3 & 2 & 1.12 & 0 & TL & 48025 & 26474 & \textBF{0.01} & 35 & 0.03 & 533 & 0.05 & 0 & 0 & 0.42 & 32 & 33 \\
\midrule
30 & 0.5 & 1 & 3 & TL & 132475 & TL & 57854 & 21568 & 553.77 & $>$5.8M & 450.33 & $>$6.6M & \textBF{0.05} & 0 & 0 & TL & 158584 & 14454 \\
&  & 2 & 3 & TL & 123439 & TL & 66413 & 20646 & 894.57 & $>$8.2M & \textBF{601.74} & $>$8.1M & TL & 0 & 2 & TL & 125635 & 16910 \\*[0.18cm]
& 0.7 & 1 & 2 & 4.74 & 310 & TL & 32065 & 40585 & \textBF{0.03} & 43 & 0.04 & 44 & TL: & 0 & 2 & 0.17 & 2 & 3 \\
&  & 2 & 2 & 2.43 & 10 & TL & 30354 & 42925 & \textBF{0.03} & 79 & 0.08 & 950 & 0.13 & 0 & 0 & 36.70 & 9460 & 6976 \\
&  & 3 & 2 & 1.57 & 0 & TL & 27713 & 39782 & \textBF{0.03} & 42 & 0.04 & 188 & 0.16 & 0 & 0 & 0.19 & 7 & 23 \\
\midrule
35 & 0.5 & 1 & 2 & TL & 64824 & TL & 39674 & 35014 & 0.97 & 11610 & 0.41 & 5392 & \textBF{0.15} & 0 & 0 & TL & 86717 & 22700 \\
&  & 2 & 2 & TL & 61939 & TL & 43800 & 34937 & 0.92 & 11223 & \textBF{0.20} & 2673 & TL & 0 & 0 & TL & 148236 & 29015 \\
&  & 3 & 2 & T;L & 61253 & TL & 35592 & 37923 & \textBF{1.35} & 15605 & 1.58 & 22016 & TL & 0 & 0 & TL & 113575 & 17829 \\*[0.18cm]
& 0.7 & 1 & 2 & 8.49 & 14 & TL & 18249 & 45102 & \textBF{0.06} & 42 & 0.06 & 55 & 0.26 & 0 & 0 & 0.35 & 11 & 36 \\
&  & 2 & 2 & 4.30 & 0 & TL & 19194 & 44076 & \textBF{0.06} & 44 & 0.06 & 61 & 0.25 & 0 & 0 & 0.29 & 8 & 24 \\
&  & 3 & 2 & 3.28 & 0 & TL & 21799 & 42927 & 0.06 & 42 & \textBF{0.06} & 57 & 0.26 & 0 & 0 & 0.55 & 19 & 21 \\
\bottomrule
\end{longtable}
\end{landscape}

	{\crev \section{Synthetic Instance Generation} \label{appendix:syn}
	Algorithm \ref{alg:syn} gives the procedure to generate synthetic instances. 
	\begin{algorithm} 
		\footnotesize
		\DontPrintSemicolon
		\SetAlgoLined
		\KwIn{$\K$, \textit{numDoubles}, \textit{noise}, $n$}
		$\doublevar := \{0\  | \ \forall \rankind \in [n-1] \} $ \;
		$\doublevar_\K := 1$\;
		\While{ $\sum_{\rankind \in [n-1]} \doublevar_\rankind < $  numDoubles}{
			$r := $ random integer in $[\K+1, n-1]$\;
			\If{$\doublevar_r = 0$}{
				$\doublevar_r := 1$\;
			}
		}
		$\edgeset = \emptyset$ \;
		\For{all $i,j$ pairs $i,j \in [\K]$}{
			add an edge $(i,j)$ to $\edgeset$ \;
		}
		\ForEach{$v\in [\K+1, n-1]$}{
			\eIf{$\doublevar_v = 1$}{
				randomly select a subset of vertices of size $\K$, $U \subseteq [v-1]$, and add edges $(u,v)$ for all $u \in U$ to $\edgeset$\;
			}{
				randomly select a subset of vertices of size $\K+1$, $U \subseteq [v-1]$, and add edges $(u,v)$ for all $u \in U$ to $\edgeset$\;
			}
		}
		$noiseCount := 0$\;
		\While{noiseCount $< \lceil$noise$\times n \rceil$}{
			randomly select $u,v \in [\K+1, n-1]$, $u \neq v$, with $\doublevar_u = 0$ and $ \doublevar_v = 0$\;
			\If{$(u,v) \notin \edgeset$}{
				add edge $(u,v)$ to $\edgeset$ \;
				$noiseCount++$	\;
		}

	}

		\Return  $G := ([n-1],\edgeset)$ \;
		\caption{Synthetic Instance Generation Procedure}
		\label{alg:syn}
	\end{algorithm}
	
}

\end{document}